%% file: 1DiagrCalc-3.tex
\documentclass[a4paper,11pt,reqno]{amsart}

\usepackage{graphicx}
\usepackage{amsmath}
\usepackage{amssymb}
\usepackage{amsfonts}
\usepackage{epsfig}
\usepackage{epigraph}

\usepackage{comment}
\usepackage{color}

\usepackage{mathrsfs}

\usepackage{graphics}
\usepackage{float}

\graphicspath{{Orbits/}
              {Exceptions/}
              {Conj_Classes/}}

\def\colorcom#1#2#3{
 \def#1##1 {
  \begin{color}{#3}
      \bf[\,#2: ##1\,]
  \end{color}}
} \colorcom\RS{RS}{red}

\makeatletter

\newtheorem{theorem}{Theorem}[section]

\@addtoreset{equation}{section}

\def\SC{{\bf c}}
\def\DualSC{{}^t\SC}
\def\InvDualSC{{\SC}^*}

\newtheorem{corollary}[theorem]{Corollary}
\newtheorem{remark}[theorem]{Remark}

\newtheorem{lemma}[theorem]{Lemma}

\def\PerfProof{{\it Proof.\ }}

\begin{document}

\title[Root systems and diagram calculus. III]
{\qquad\qquad Root systems and diagram calculus. \newline
  III. Semi-Coxeter orbits of linkage diagrams and the Carter theorem}
         \author{Rafael Stekolshchik}

\date{}

\begin{abstract}
   A diagram obtained from the Carter diagram $\Gamma$ by adding one root together with its bonds such that the resulting subset of roots is linearly independent is said to be the {\it linkage diagram}.
   Given a linkage diagram, we associate the linkage labels vector, which is introduced like the vector of Dynkin labels.
   Similarly to the dual Weyl group, we introduce the group $W^{\vee}_L$ associated with $\Gamma$, and we call it
   the dual partial Weyl group.  The linkage labels vectors connected under the action of $W^{\vee}_L$ constitute the linkage system $\mathscr{L}(\Gamma)$,  which is similar to the weight system arising in the representation theory of the semisimple Lie algebras.
   The Carter theorem  states that every element of a Weyl group $W$ is expressible as the product of
   two involutions. We give the proof of this theorem based on the description of the linkage system $\mathscr{L}(\Gamma)$ and semi-Coxeter orbits of linkage labels vectors for any Carter diagram $\Gamma$.
   The main idea of the proof is based on the fact that, with a few exceptions,
   in each semi-Coxeter orbit there is a special linkage diagram -- called {\it unicolored}, for which
   the decomposition into the product of two involutions is trivial.
\end{abstract}

\maketitle

~\\
~\\
~\\
~\\
\tableofcontents

\newpage
~\\
~\\

\setlength{\epigraphwidth}{105mm}

\epigraph{Gelfand requested that I review the H. Weyl - Van der Waerden papers on semisimple Lie groups.
I found them very difficult to read, and I tried to find my own ways. It came to my mind that
there is a natural way to select a set of generators for a semisimple Lie algebra by using
simple roots (i.e., roots which cannot be represented as a sum of two positive roots).
Since the angle between any two simple roots can be equal
only to $\pi/2$, $2\pi/3$, $3\pi/4$, $5\pi/6$, a system of simple roots can be represented
by a simple diagram. An article was submitted to {\it Matematicheskii Sbornik}
in October $1944$, \cite{Dy46}. Only a few years later, when recent literature from the
West reached Moscow, I discovered that similar diagrams have been used
by Coxeter for describing crystallographic groups.}
{E. B. Dynkin,
Foreword in \lq\lq Selected papers of E. B. Dynkin with commentary \rq\rq , \cite[p. 2]{Dy00}}

\section{\sc\bf Introduction}
\subsection{The Carter theorem}

In the present paper we give the proof of the Carter theorem on the
decomposition of every element of any Weyl group $W$ into the
product of two involutions.

\begin{theorem}[\cite{Ca72},Theorem C]
 \label{th_Carter}  The following equivalent statements hold for the Weyl group:

  {\rm (i)}  Every element of a Weyl group $W$ is expressible as the product of
two involutions.

  {\rm (ii)} Every element of $W$ is contained in some dihedral subgroup.

  {\rm (iii)} For each element $w \in W$ there is an involution $i \in W$ such that $iwi = w^{-1}$ .
\end{theorem}

\begin{corollary}
  Every element of $W$ is conjugate to its inverse.
\end{corollary}

In \cite{Sp74}, Springer  gives a proof of the Carter theorem for
all finite Coxeter groups including the non-crystallographic cases
$I_n$ (dihedral group), $H_3$ and $H_4$. Springer deduced the proof
from the classification of so-called regular elements in the Coxeter
groups and by inspection  from the known character tables of the
irreducible Weyl groups, \cite[\S8.6, \S8.7]{Sp74}.

 The proof given by Carter in \cite{Ca72} uses the calculation of all conjugacy classes in the Weyl group.
 Our proof uses the classification of linkage systems and semi-Coxeter orbits for every Carter diagrams.
 The definitions of linkage systems and semi-Coxeter orbits will be given below in Section \ref{sec_Carter_diagr}.
 The linkage systems for Carter diagrams from $\mathsf{C4} \coprod \mathsf{DE4}$ are presented in \cite{St10.II},
 where $\mathsf{C4}$ is the class of Carter diagrams, each of which contains $4$-cycle $D_4(a_1)$ as a subdiagram, and
 $\mathsf{DE4}$ is the class of Carter diagrams, each of which contains $D_4$ as a subdiagram,
 see Section \ref{sec_classes_diagr}. In this paper, we give the complete description of semi-Coxeter orbits for
 any $\Gamma \in \mathsf{C4} \coprod \mathsf{DE4}$. The linkage systems and semi-Coxeter orbits for $A_l$ and $B_l$ will be presented in \cite{St11}.

\input 1definition.tex

\input 2CarterTh.tex
\input 3exceptional.tex

\input 4finding.tex

\begin{appendix}
\vspace{15mm}
\input A01matrSemiCox.tex

\input A1orbits.tex

\end{appendix}

\input biblio-3.tex

\end{document}

%% file: 1definition.tex
~\\
\subsection{The Carter diagrams and semi-Coxeter elements}
  \label{sec_Carter_diagr}
\subsubsection{Solid and dotted edges}
  \label{sec_solid_dot}

  The Carter diagram (= admissible diagram) \cite[\S4]{Ca72}
  is the diagram $\Gamma$  satisfying two conditions:
  \vspace{3mm}

  (a) The nodes of $\Gamma$ correspond to a set of linearly independent roots.

  (b) Each subgraph of $\Gamma$ which is a cycle contains even number of vertices.
   \vspace{3mm}

Let $w = w_1 w_2$ be the decomposition of $w$ into the product of
two involutions. By \cite[Lemma 5]{Ca72} each of $w_1$ and $w_2$ can
be expressed as a product of reflections corresponding to mutually orthogonal roots:
\begin{equation}
   \label{two_invol}
       w = w_1 w_2, \quad
       w_1 = s_{\alpha_1} s_{\alpha_2} \dots s_{\alpha_k}, \quad
       w_2 = s_{\beta_1} s_{\beta_2} \dots s_{\beta_h}, 
       \quad \text{where} \quad k + h = l_C(w).
 \end{equation}
For details, see \cite[\S4]{Ca72}, \cite[\S1.1]{St10.I}.  
We denote by {\it $\alpha$-set} (resp. {\it $\beta$-set}) the subset of
 roots corresponding to $w_1$ (resp. $w_2$):
\begin{equation}
   \label{two_sets}
       \alpha\text{-set} = \{ \alpha_1, \alpha_2, \dots, \alpha_k \}, \quad
       \beta\text{-set} =  \{ \beta_1, \beta_2, \dots, \beta_h \}.
 \end{equation}
 Any coordinate from $\alpha$-set (resp. $\beta$-set) of the linkage labels vector
 we call {\it $\alpha$-label} (resp. {\it $\beta$-label}).
 The decomposition \eqref{two_invol} is said to be the {\it bicolored decomposition}.

  For the Dynkin diagrams, a number of bonds for non-orthogonal roots
  describes the angle between roots, and the ratio of lengths of two roots.
  For the Carter diagrams, we add designation distinguishing acute and obtuse angles
  between roots. Recall, that for the Dynkin diagrams, all angles between simple roots
  are obtuse and a special designation is not necessary.
  A {\it solid edge} indicates an obtuse angle between roots exactly as for simple roots
  in the case of Dynkin diagrams. A {\it dotted edge} indicates an acute angle between the
  roots cinsidered, see \cite{St10.I}. For examples of diagrams with dotted and solid edges, see Table \ref{tab_repres_exc}.

\subsubsection{Semi-Coxeter elements}
  \label{sec_semi_cox_elem}
   A conjugacy class of $W$  which can be described by a connected
Carter diagram with number of nodes equal to the rank of $W$ is called
a {\it semi-Coxeter} class, \cite{CE72} (or, a {\it primitive} conjugacy class, \cite{KP85}).
Let us fix some basis of roots corresponding to the given Carter diagram $\Gamma$:
\begin{equation}
  \label{eq_basis}
   \{\alpha_1,\dots,\alpha_k,\beta_1,\dots,\beta_h \},
\end{equation}
where $\alpha_i$, $\beta_j$ are roots (not necessarily simple) corresponding 
to $\Gamma$.
The element
\begin{equation}
  \label{eq_semiCox}
        {\bf c} = w_{\alpha}w_{\beta}, \quad \text{where}  \quad
        w_{\alpha} = \prod\limits_{i=1}^k{s}_{\alpha_i}, \quad
        w_{\beta} = \prod\limits_{j=1}^h{s}_{\beta_j},
\end{equation}
given in the basis \eqref{eq_basis} we call the semi-Coxeter element. It is the
representative of the semi-Coxeter class. The dual semi-Coxeter element
\begin{equation}
  \label{eq_semiCox_2}
        {\bf c}^* = {}^t{w}_{\alpha}{}^t{w}_{\beta}, \quad \text{where}  \quad
        w_{\alpha} = \prod\limits_{i=1}^k{s}^*_{\alpha_i}, \quad
        w_{\beta} = \prod\limits_{j=1}^h{s}^*_{\beta_j},
\end{equation}
is used for the proof of the Carter theorem. Semi-Coxeter elements for diagrams
$D_l$, $D_l(a_k)$, $E_l$, $E_l(a_k)$, where $l \leq 7$, are presented in
Tables \ref{tab_partial Cartan_1}-\ref{tab_partial Cartan_3}.

Note that roots \eqref{eq_basis} are not necessarily simple.
If all roots \eqref{eq_basis} are simple, the Carter diagram $\Gamma$
is a Dynkin diagram and the semi-Coxeter element \eqref{eq_semiCox}
coincides with the corresponding Coxeter element.

\subsubsection{Linkages, linkage diagrams and linkage systems}
\label{sec_linkage}

 Let $w = w_1w_2$ be the bicolored decomposition of some element
 $w \in W$, where $w_1$, $w_2$ are two involutions, associated, respectively,
 with $\alpha$-set $\{ \alpha_1, \dots, \alpha_k \}$ and $\beta$-set
 $\{ \beta_1, \dots, \beta_h \}$  of roots from the root system $\varPhi$, see
 \eqref{two_invol}, \eqref{two_sets}, and $\Gamma$ be the Carter
 diagram associated with this bicolored decomposition.
 We consider the {\it extension} of the root basis $\Pi_w$ by means of the root $\gamma \in \varPhi$,
 such that the set of roots
 \begin{equation}
   \label{alpha_beta}
    \Pi_w(\gamma) = \{ \alpha_1, \dots, \alpha_k, \beta_1, \dots, \beta_h, \gamma \}
 \end{equation}
 is linearly independent. Let us multiply $w$ on the right by the reflection
 $s_{\gamma}$ corresponding to $\gamma$ and consider the diagram
 $\Gamma' = \Gamma \cup \gamma$ together with new edges.
 These edges are
 ~\\
\begin{minipage}{10.8cm}
 \begin{equation*}
   \begin{cases}
  \text{{\it solid}, for $(\gamma, \tau) = -1$}, \\
  \text{{\it dotted}, for $(\gamma, \tau) = 1$}, \\
   \end{cases}
 \end{equation*}
 where $\tau$ one of elements \eqref{alpha_beta}. The diagram $\Gamma'$
 is said to be the {\it linkage diagram}, and the root $\gamma$ is said to be the {\it linkage} or the {\it $\gamma$-linkage}. 
 Consider vectors $\gamma^{\vee}$ belonging to the dual space $L^{\vee}$ and defined by
 \eqref{dual_gamma}.
 Vector \eqref{dual_gamma} is said to be {\it linkage labels vector} or, for brevity, {\it linkage labels}.
  There is, clearly, the one-to-one correspondence between linkage labels vectors $\gamma^{\vee}$
 (with labels $\gamma^{\vee}_i \in \{0, -1, 1\}$) and simply-laced linkage diagrams
 (i.e., such linkage diagrams that  $(\gamma, \tau) \in \{0, -1, 1\}$).
 \end{minipage}
\begin{minipage}{5.5cm}
  \quad
 \begin{equation}
   \label{dual_gamma}
   \gamma^{\vee} :=
   \left (
    \begin{array}{c}
    (\gamma, \alpha_1) \\
    \dots,     \\
    (\gamma, \alpha_k)  \\
    (\gamma, \beta_1)  \\
    \dots,  \\
    (\gamma, \beta_h) \\
    \end{array}
   \right )
 \end{equation}
\end{minipage}
 
   Let $L$ be the linear space spanned by the roots associated with $\Gamma$,
 The linkage labels vector is the element of the dual linear space $L^{\vee}$. We denote
 the linkage labels vector by $\gamma^{\vee}$. A certain group $W^{\vee}_L$ named
 the dual partial Weyl group acts in the dual space $L^{\vee}$.
 This group acts on the linkage label vectors, i.e., on the set of linkage diagrams:
 \begin{equation*}
      (w\gamma)^{\vee} = w^{*}\gamma^{\vee},
 \end{equation*}
 where $w^{*} \in W^{\vee}_L$, see Proposition $2.9$ from \cite{St10.II}. The set of
 linkage diagrams (=linkage labels) under action of $W^{\vee}_L$ constitute the
 diagram called the {\it linkage system} similarly to the weight system
 in the theory of representations  of semisimple Lie algebras,  see \cite[p. 30]{Sl81}, \cite[p. 4]{St10.II}.
\begin{remark}
  \label{rem_abuse}
 \rm{
 By abuse of language, we sometimes say {\it linkages} instead of {\it linkage diagrams}.
 Similarly, remembering only the algebraic nature of the linkage diagram, we use
 the term {\it linkage label vector} or {\it linkage labels}.
 }
 \end{remark}
 \subsubsection{Classes of Carter diagrams}
   \label{sec_classes_diagr}
 We divide all Carter diagrams to the following classes:
 ~\\
 
 {\it Simply-laced Carter diagrams:}

  1. $\mathsf{DE4}$, Dynkin diagrams containing $D_4$ as a subdiagram,

  2. $\mathsf{C4}$, Carter diagrams containing $4$-cycle $D_4(a_1)$ as a subdiagram,

  3. $\mathsf{A}$,  Dynkin diagrams $A_l$,
 ~\\
 
{\it  Multiply-laced Carter diagrams:}

  4. $\mathsf{BC}$, Dynkin diagrams $B_l, C_l$,

  5. $\mathsf{FG}$, Dynkin diagrams $F_4$, $G_2$, and the $4$-cycle with two double bonds $F_4(a_1)$.
~\\

 For $\Gamma \in \mathsf{C4} \coprod \mathsf{DE4}$, the linkage systems are described in \cite{St10.II}.
 In this case, for $l \leq 7$, the linkage systems $\mathscr{L}(\Gamma)$ looks as follows:
 every linkage diagram containing at least one non-zero $\alpha$-label (see Section \ref{sec_solid_dot})
 belongs to a certain $8$-cell "spindle-like" linkage subsystem called
 {\it loctet} (= linkage octet). The loctets are the main construction blocks in every linkage system.
 If all $\alpha$-labels (resp. $\beta$-labels) of $\gamma^{\vee}$ are zeros, the linkage diagram  $\gamma^{\vee}$ 
 is said to be {\it $\beta$-unicolored} (resp. {\it $\alpha$-unicolored}) linkage diagram. 
 Every linkage system is the union of several loctets and several $\beta$-unicolored linkage diagrams, see \cite[\S3]{St10.II}.
 In the case, where $l > 7$, the linkage systems for two infinite series $D_l$ and $D_l(a_k)$ are described as follows: for $D_l(a_k)$ the linkage system looks as {\it wind rose of linkages}, see
 \cite[Fig. B.46-B.47]{St10.II}; for $D_l$ the linkage system looks as the Carter diagram $D_l(a_k)$,
 see \cite[Fig. B.48]{St10.II}.
   
 \begin{remark}[multiply-laced Carter diagrams $\mathsf{FG}$]
 {\rm
 For any $\Gamma \in \mathsf{FG}$, the linkage system is trivial. 
 Really, for the case $G_2$, there are maximum two linearly independent roots. Thus, the linkage system $\mathscr{L}(G_2)$ is trivial,  see \cite[Rem. 2.2]{St10.I}.
 Further, for the multiply-laced $4$-cycle $F_4(a_1)$,
 there is no additional fifth edge, otherwise such a diagram contains an extended Dynkin diagram as a subdiagram,
 see \cite[\S A.3.2]{St10.I}. The simple extension of $F_4$ leads to the subdiagram,
 which is one of extended Dynkin diagrams $\widetilde{F}_{41}$, $\widetilde{F}_{42}$, see \cite[Example A.3]{St10.I}, or $\widetilde{CD}_{n}$, $\widetilde{DD}_{n}$ that can not be. Any triangle extending $F_4$ is also moved to one
 of cases $\widetilde{F}_{41}$, $\widetilde{F}_{42}$, see \cite[\S 4.3]{St10.II}. \qed
 }
 \end{remark}
 In \cite{St11}, we will construct remaining cases of linkage systems for two infinite series $A_l$ and $B_l$.

%% file: 2CarterTh.tex
\section{\sc\bf The proof of the Carter theorem}

\subsection{Linear independency and reduced decomposition}

\subsubsection{Reduced decomposition and the Carter length $l_C(w)$}
  \label{sec_reduced}
 Each element $w \in W$ can be expressed in the form
 ~\\
 \begin{equation}
   \label{any_roots}
    w  = s_{\tau_1} s_{\tau_2} \dots s_{\tau_k}, \quad \tau_i \in \varPhi,
 \end{equation}
 ~\\
 where $\varPhi$ is the root system associated with
 the Weyl group $W$; $s_{\tau_i}$ are reflections in $W$
 corresponding to not necessarily simple roots $\tau_i \in \varPhi$.
 We denote by $l_C(w)$ the smallest value $k$ in any expression like \eqref{any_roots}.
 The Carter length $l_C(w)$ is always less than the classical length $l(w)$.
 The decomposition \eqref{any_roots} is called {\it reduced} if
 $l_C(s_{\tau_1} s_{\tau_2} \dots s_{\tau_k}) = k$, i.e., the number of
 reflections in \eqref{any_roots} can not be decreased.

 \begin{lemma}{\cite[Lemma 3]{Ca72}}
  \label{lem_lin_indep}
   Let $\tau_1, \tau_2, \dots, \tau_k \in \varPhi$.
   Then $s_{\tau_1} s_{\tau_2} \dots s_{\tau_k}$ is reduced
   if and only if $\tau_1, \tau_2, \dots, \tau_k$ are linearly
   independent.
 \end{lemma}  \qed

\subsubsection{The basic conjugacy relation}

\begin{lemma}[on conjugacy]
  \label{lem_indep_gamma}

Let $\{ \tau_1, \dots, \tau_n \}$ be the subset of linearly independent roots
(not necessarily simple), $\tau_i \in \varPhi$, and let $w \in W$ be the element, which is decomposed
into the product of reflections $\{ s_{\tau_1}, \dots, s_{\tau_n} \}$.

 1) If $\gamma$ such a root that $\{\gamma, \tau_1, \dots, \tau_n \}$ are linearly
 independent then $\{ w\gamma, \tau_1, \dots, \tau_n \}$ are also linearly
 independent.

 2) The following conjugacy relation holds for any integer $k$:
~\\
\begin{equation*}
    s_{\gamma}w    \simeq s_{w^k\gamma}w. \\
\end{equation*}
~\\
In particular, for the semi-Coxeter element ${\bf c}$,  we have
~\\
\begin{equation}
  \label{semiCox_basic}
    s_{\gamma}\SC  \simeq s_{\SC^k\gamma}\SC.
\end{equation}
\end{lemma}
~\\

\PerfProof  1) Since $s_{\tau_i}(\tau_j) \in \{ \tau_i, \tau_j \}$, and
$s_{\tau_i}(\gamma) \in \{ \gamma, \tau_i \}$, we have
~\\
\begin{equation*}
  \begin{split}
   & w\{\tau_1, \dots, \tau_n \} \subseteq \{\tau_1, \dots, \tau_n \}, \text{ and } \\
   &  w\gamma \in \{\gamma, \tau_1, \dots, \tau_n \}, \text{ i.e., } 
      w\gamma = \gamma + \sum\limits_{i=1}^n{a_i}{\tau_i}
  \end{split}
\end{equation*}
~\\
for some rational factors $a_i$. If $\{ w\gamma, \tau_1, \dots, \tau_n \}$ are linearly dependent, we have
~\\
\begin{equation*}
     w\gamma = \sum\limits_{i=1}^n{b_i}{\tau_i} \text{ for some rational $b_i$,} \text{ i.e., }
       \gamma + \sum\limits_{i=1}^n{a_i}{\tau_i} = \sum\limits_{i=1}^n{b_i}{\tau_i}
\end{equation*}
~\\
 that contradicts to the linear independency of $\{\gamma, \tau_1, \dots, \tau_n \}$.

 2)  Let
\begin{equation}
    w  = \prod\limits_{i=1}^m{s}_{\tau_i},
\end{equation}
where not necessarily all $\tau_i$ are different. For example, it can be that $m > n$.
We have

\begin{equation}
 \label{semiCox_orbit}
  \begin{split}
   s_{\gamma}w  = & s_{\gamma}s_{\tau_1}\dots{s}_{\tau_m} = s_{\tau_1}s_{s_{\tau_1}(\gamma)}
     s_{\tau_2}\dots{s}_{\tau_m} =
     s_{\tau_1}s_{\tau_2}s_{s_{\tau_2}s_{\tau_1}(\gamma)}s_{\tau_3}\dots{s}_{\tau_m} =
        \dots = \\
   & w{s}_{w^{-1}\gamma}  \simeq {s}_{w^{-1}\gamma}w \simeq \dots \simeq {s}_{w^{-k}\gamma}w,
       \text{ for } k > 0.
   \end{split}
\end{equation}
~\\
By mapping $\gamma' = w^{-k}\gamma$ we obtain from \eqref{semiCox_orbit} also that
\begin{equation*}
    s_{\gamma'}w \simeq {s}_{w^k\gamma'}
\end{equation*}
 for any integer $k > 0$ and every root $\gamma'$. Thus, \eqref{semiCox_basic} is proven.
\qed

Eq. \eqref{semiCox_basic} is the basic relation in our proof of the Carter theorem,
see Section \ref{sec_inducition}.

~\\
 \subsection{The proof of Theorem \ref{th_Carter}}
\subsubsection{The induction step}
  \label{sec_inducition}
 Let $\SC$  be the semi-Coxeter element associated with the Carter diagram $\Gamma$
 such that $\SC$ has bicolored decomposition given by \eqref{eq_semiCox}.
The proof of the theorem is carried out by induction  on the
Carter length $l_C(\SC)$ of the decomposition \eqref{eq_semiCox}, see Section \ref{sec_reduced}.
For details, see \cite{Ca72}, \cite[p. 4]{St10.I}.
Suppose, $\gamma$ is the root such that roots $\{\gamma, \alpha_1, \dots, \beta_h \}$
are linearly independent. According to Lemma \ref{lem_indep_gamma}, heading 1)
$\{\SC^n\gamma, \alpha_1, \dots, \beta_h \}$ are also linearly independent.

We will show that
\begin{equation}
  \label{semiCox_2}
   s_{\gamma}\SC  \text{ \it is also associated with a certain Carter diagram}.
\end{equation}
Of course, it suffices to prove the property \eqref{semiCox_2} for any conjugate of $s_{\gamma}\SC$.
The property \eqref{semiCox_2} gives us the induction step. According to Lemma \ref{lem_indep_gamma},2)
it suffices to find such an integer $n$ that any conjugate of $s_{\SC^n\gamma}\SC$
has the bicolored decomposition.

\subsubsection{Semi-Coxeter orbits}
  \label{sec_semi_Coxeter}
We have $(\SC^n\gamma)^{\vee} = ({\SC}^n)^*\gamma^{\vee} = ({\SC}^*)^n\gamma^{\vee}$ for any $n$, see \cite[Proposition 2.9]{St10.II}. Let us consider the sequence of linkages
\begin{equation}
  \label{eq_seq_linkages}
   (\SC^n\gamma)^{\vee} = ({\SC}^*)^n\gamma^{\vee}, \quad n = 0, \pm{1}, \pm{2}, \dots
\end{equation}
It is clear that \eqref{eq_seq_linkages} is the finite periodic sequence, 
see Tables \ref{tab_partial Cartan_1}-\ref{tab_partial Cartan_3}. 
This sequence is said to be the {\it semi-Coxeter orbit}. Remember, that the linkage diagram $\gamma^{\vee}$
is said to be the $\alpha$-unicolored (resp. $\beta$-unicolored) linkage diagram if all $\beta$-labels (resp. $\alpha$-labels),  i.e., coordinates corresponding to all $\beta_i$ (resp. $\alpha_i$) of  $\gamma^{\vee}$ are zeros,
 \cite[p. 5]{St10.II}. Suppose, for some integer $m$, the element $({\SC}^*)^m\gamma^{\vee}$ in semi-Coxeter orbit is
 a certian unicolored linkage diagram. Let $({\SC}^*)^m\gamma^{\vee}$ be, for example, $\alpha$-unicolored. Then $(\SC^m\gamma, \beta_i) = 0$ for all $\beta$-labels. This means that $s_{{\SC}^m\gamma}$ commute with all $s_{\beta_i}$. By
 \eqref{semiCox_basic} and \eqref{eq_semiCox}
\begin{equation*}
     s_{\gamma}\SC = s_{{\SC}^m\gamma}\SC =
     s_{{\SC}^m\gamma}\prod\limits_{i=1}^k{s}_{\alpha_i} \prod\limits_{j=1}^h{s}_{\beta_j}
     \simeq
     \prod\limits_{i=1}^k{s}_{\alpha_i} (\prod\limits_{j=1}^h{s}_{\beta_j})s_{{\SC}^m\gamma}.
\end{equation*}
The latter product is the bicolored decomposition, since $(\prod\limits_{j=1}^h{s}_{\beta_j})s_{{\SC}^m\gamma}$
is involution. Thus, it suffices to prove that any semi-Coxeter orbit contains an unicolored linkage diagram.

\subsubsection{Unicolored linkage diagrams and exceptional orbits}  However, there are semi-Coxeter orbits containing no unicolored linkage diagrams. We call these orbits {\it exceptional semi-Coxeter orbits}. The total quantity of orbits for Carter diagrams from $\Gamma \in \mathsf{C4} \coprod \mathsf{DE4}$  (for $l \leq 7$) is $140$,
the number of exceptional orbits is $24$, see Table \ref{tan_numb_len_orbits}.
Instead $24$ orbits it suffices to consider only $10$ exceptional orbits,
 namely $(1a)$, $(2a)$, $(2c)$, $(3a)$, $(4a)$, $(4b)$, $(5a)$, $(6a)$, $(7a)$, $(7b)$,
 see Table \ref{tab_repres_exc},
 they are checked case-by-case in Section \ref{sec_exc_cases}.
\begin{table}[H]
  \centering
  \renewcommand{\arraystretch}{1.3}
  \begin{tabular} {||c|c|c|c|c|c||}
  \hline \hline
      The Carter &   \multicolumn{3}{c|}{Number of orbits}          & Lengths of & Number of  \cr
      \cline{2-4}
        diagram  &  All   & no unicolored    &  self-opposite &  orbits\footnotemark[1]    & linkages \cr
                 &        & linkages         &                &            &           \cr
    \hline \hline
       $D_4(a_1)$ &   $6$  & $2$   & $6$  & $6 \times 4$ & $24$ \\
    \hline
       $D_4$      &    $6$  & -   & $6$  & $3 \times 6 + 3 \times 2$ & $24$ \\
    \hline \hline
       $D_5(a_1)$ &   $6$  & -   & $2$ & $2 \times 12 + 3 \times 4 + 6$ & $42$  \\
    \hline
       $D_5$      &   $6$  & -   & $2$ & $5 \times 8 + 2$ & $42$  \\
    \hline \hline
       $E_6(a_1)$ &   $6$  & -  &  -  & $9 \times 9$ & $54$  \\
    \hline
       $E_6(a_2)$ &   $10$  & $4$  & - & $8 \times 6 + 2 \times 3$ & $54$  \\
    \hline
       $E_6$      &   $6$  &  -  &  -  & $4 \times 12 + 2 \times 3$ & $54$ \\
    \hline \hline
       $D_6(a_1)$ &   $10$  &  $4$  & 2  & $9 \times 8 + 4$ & $76$  \\
    \hline
       $D_6(a_2)$ &   $14$  &  $2$ & $14$ & $12 \times 6 + 2 \times 2$ & $76$  \\
    \hline
       $D_6$      &   $8$  &  - &  8 & $7 \times 10 + 3 \times 2$ & $76$\\
    \hline \hline
       $E_7(a_1)$ &   $4$  &   -  &  $4$ & $4 \times 14$ & $56$  \\
    \hline
       $E_7(a_2)$ &   $6$  &   $2$  & $2$ & $4 \times 12 + 6 + 2$ & $56$  \\
    \hline
       $E_7(a_3)$ &   $4$  &  -  & $4$  & $30 + 2 \times 10 + 6$ & $56$  \\
    \hline
       $E_7(a_4)$ &   $10$  &  $6$  & $10$  & $9 \times 6 + 2$ & $56$   \\
    \hline
       $E_7$      &   $4$  &  - & $4$  & $3 \times 18 + 2$ & $56$   \\
    \hline \hline
       $D_7(a_1)$ &   $10$  & -  & $2$  & $6 \times 20 + 2 \times 4 + 10 + 4$ & $142$  \\
    \hline
       $D_7(a_2)$ &   $10$  &  -  &  $2$  & $4 \times 24 + 4 \times 8 + 8 + 6$ & $142$   \\
    \hline
       $D_7$      &   $14$  &  4 &  $2$  & $10 \times 12 + 2 \times 4 + 12 + 2$ & $142$   \\
    \hline \hline
       $D_l(a_k)$, $l > 7$ &   $2$  &   -  & $2$ & $2 \times (k+1) + 2 \times (l - k - 1)$ & $2l$  \\
    \hline
       $D_l$, $l >  7$ &   $2$  &   -  & $2$ & $2 \times (l-1)$ & $2l$  \\
   \hline  \hline  
\end{tabular}
  \vspace{2mm}
  \caption{\hspace{3mm}Number and lengths of semi-Coxeter orbits }
  \label{tan_numb_len_orbits}
  \end{table}
\begin{remark}{\rm
  We observe that the number of unicolored linkage diagrams in every semi-Coxeter orbit is equal to $0$
  (exceptional orbit) or $2$, see Tables \ref{tab_D4} - \ref{E7a4_linkages} (where unicolored linkage diagrams
  are framed by a rectangle). Of course, this fact requires {\it a priori} reasoning.
}
\end{remark}
\footnotetext[1]{Explanation to the column. For example, expression $6 \times 20 + 2 \times 4 + 10 + 4$
  in the line $D_7(a_1)$ means that the total number of linkage diagrams in the linkage system for the Carter diagram  $D_7(a_1)$ is divided into the sum of $6$ orbits
  each of $20$ elements, $2$ orbits each of $4$ elements, one orbit containing $10$ elements and one orbit containing $4$ elements.}

\subsubsection{Semi-Coxeter orbits for infinite series $D_l(a_k)$ and $D_l$}

\begin{figure}[H]
\centering
\includegraphics[scale=1.5]{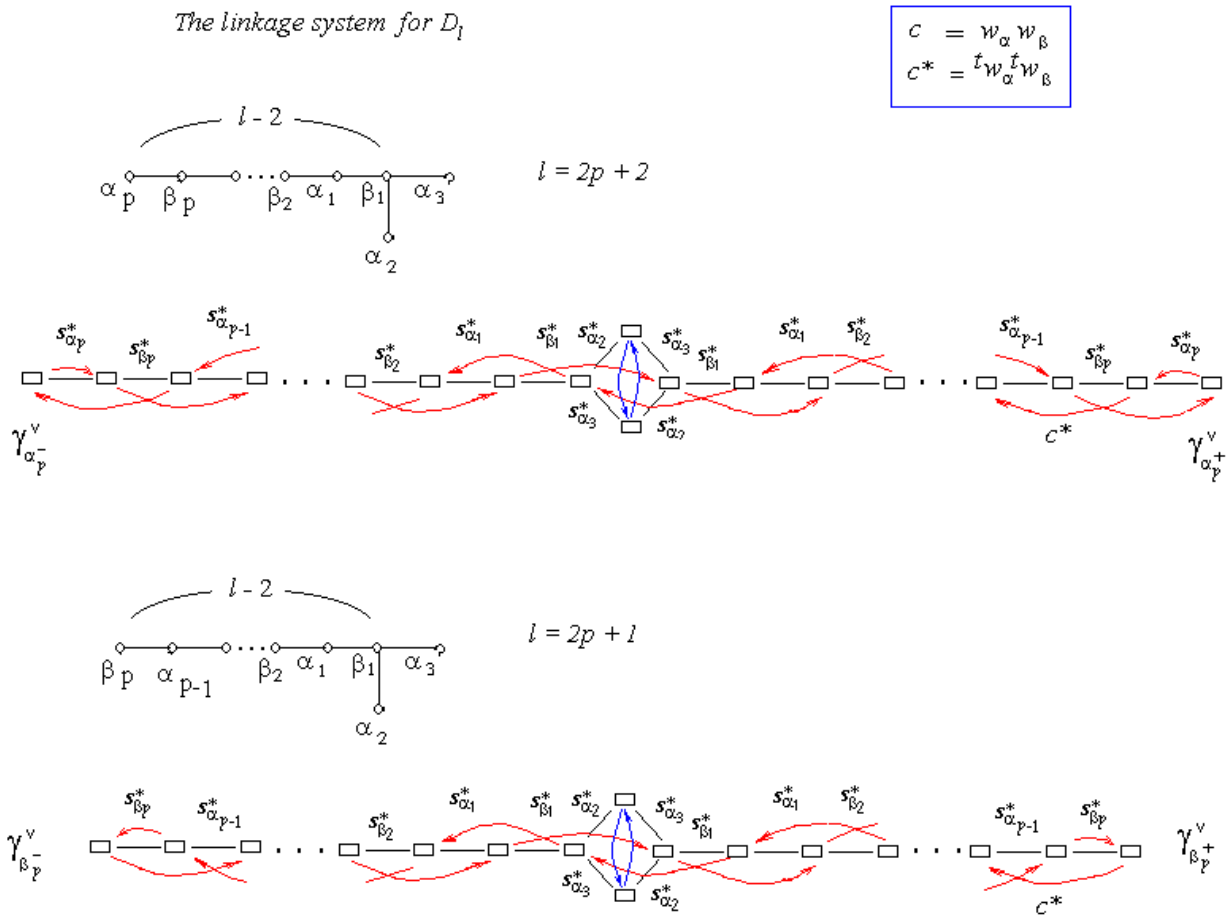}
\vspace{3mm}
\caption{\hspace{3mm}\footnotesize Two semi-Coxeter orbits of $D_l$, one of length $2(l-1)$ (= Coxeter number), one of length $2$}
\label{Dl_orbits_16c}
\end{figure}

\begin{figure}[H]
\centering
\includegraphics[scale=1.5]{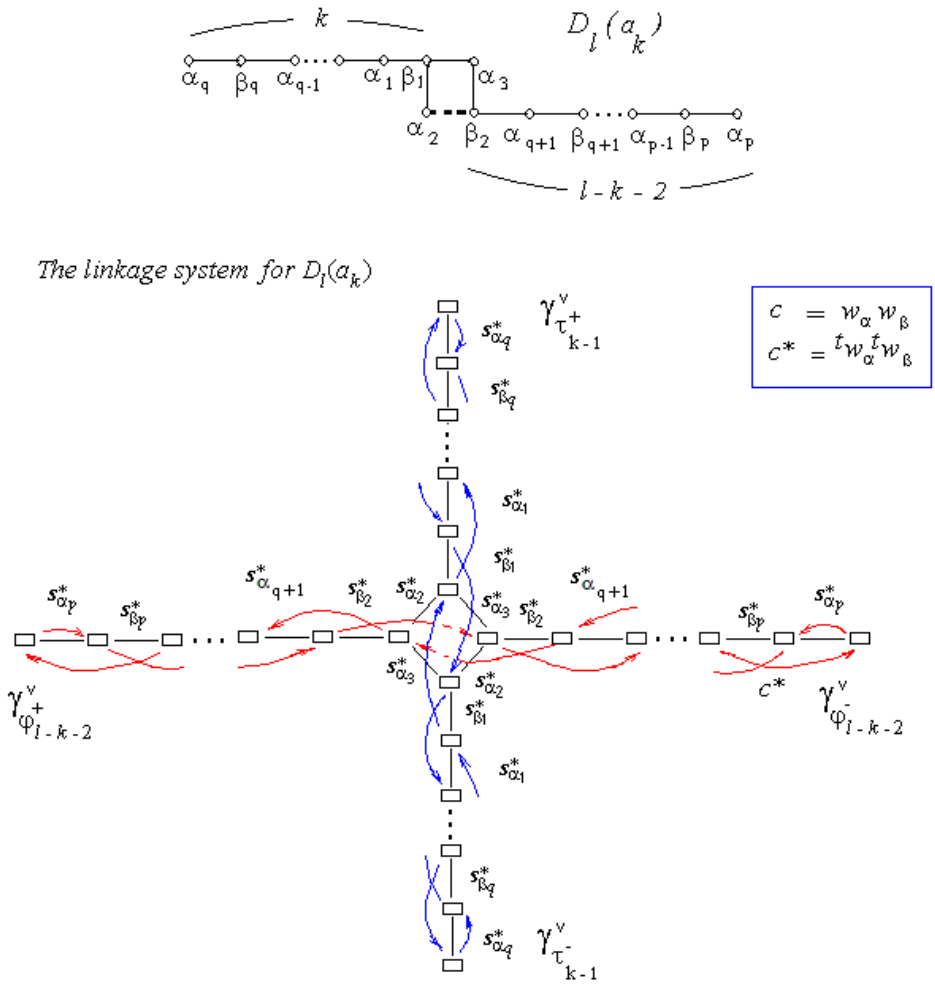}
\vspace{3mm}
\caption{\hspace{3mm}\footnotesize Two semi-Coxeter orbits of $D_l$:
one of length $2(k+1)$, one of length $2(l-k-1)$}
\label{Dlak_orbits_16c}
\end{figure}

 It is convenient to imagine a semi-Coxeter orbit for $D_l$ (resp. $D_l(a_k)$) as a cosine wave
 that runs in one direction and then returns in the opposite direction with a shift in the phase by half a period,
 see Fig. \ref{Dl_orbits_16c} and Fig. \ref{Dlak_orbits_16c}.

 \underline{Diagram  $D_l$.}
 We have one long orbit -- the red wave in the horizontal direction, and one $2$-element orbit -- the blue orbit in the vertical direction, see Fig. \ref{Dl_orbits_16c}. There are two cases: $l = 2p+2$ and $l = 2p+1$.
 For $l = 2p+2$, the linkage labels vector $\gamma^{\vee}_{\alpha^{+}_p}$ (resp. $\gamma^{\vee}_{\alpha^{-}_p}$)
 is the vector with the unit in the place $\alpha^{+}_p$ (resp. $\alpha^{-}_p$) and zeros in remaining places, see Fig. \ref{Dl_orbits_16c}. 
  They are two unicolored linkages for the long semi-Coxeter orbit (red wave).
 The blue orbit consists of following two unicolored linkages:
 $\gamma^{\vee}_4 = \{ 0, 1, -1, 0, \dots, 0 \}$, $\gamma^{\vee}_5 = - \gamma^{\vee}_4$
 with the only non-zero coordinates in coordinates $\alpha_2$ and $\alpha_3$. Notations of $\gamma^{\vee}_4$,  $\gamma^{\vee}_5$ are retained as in \cite[Fig. B.48]{St10.II}.  For $l = 2p+1$,
 the linkage labels vector $\gamma^{\vee}_{\beta^{+}_p}$ (resp. $\gamma^{\vee}_{\beta^{-}_p}$)
 is the vector with the unit in the place $\beta^{+}_p$ (resp. $\beta^{-}_p$) and zeros in remaining places, see Fig. \ref{Dl_orbits_16c}. The blue orbit is the same as in the case $l = 2p+2$.
  Linkages $\gamma^{\vee}_{\alpha^{+}_p}$ and $\gamma^{\vee}_{\alpha^{-}_p}$ (see Remark \ref{rem_abuse}) 
 for $l = 2p+2$, and linkages $\gamma^{\vee}_{\beta^{+}_p}$ and $\gamma^{\vee}_{\beta^{-}_p}$
 for $l = 2p+1$ are the same linkages as $\gamma^{\vee}_{\tau^{+}_{l-3}}$ and $\gamma^{\vee}_{\tau^{-}_{l-3}}$
 in \cite[Fig. B.48]{St10.II}.

 \underline{Diagram  $D_l(a_k)$.}
  Here, we have two long orbits: one red wave in the horizontal direction,
 and one blue wave in the vertical direction, see Fig. \ref{Dlak_orbits_16c}.
 The linkage labels vector  $\gamma^{\vee}_{\tau^+_{k-1}}$ (resp. $\gamma^{\vee}_{\tau^-_{k-1}}$,
 $\gamma^{\vee}_{\varphi^{+}_{l-k-2}}$, $\gamma^{\vee}_{\varphi^{-}_{l-k-2}}$)
 is the vector with the unit on the place
 $\tau^+_{k-1}$ (resp. $\tau^-_{k-1}$,  $\varphi^{+}_{l-k-2}$,  $\varphi^{-}_{l-k-2}$)
 and zeros on remaining places, see Fig. \ref{Dlak_orbits_16c}. As above, these vectors are unicolored,
 see \cite[Fig. B.46-B.47]{St10.II}.

 Two semi-Coxeter orbits of $D_l(a_k)$ are of lengths $2(k+1)$ and $2(l-k-1)$. 
 For the left (resp. right) branch of $D_l(a_k)$, 
 there are two options for endpoints: $\alpha_p$ or $\beta_p$ (resp. $\alpha_q$ or $\beta_q$).
 Thus, from the view of endpoints there are $4$ options for the Carter diagram $D_l(a_k)$:
 \begin{equation*}
    \{\alpha_p, \alpha_q \}, \quad \{\alpha_p, \beta_q \}, \quad
    \{\beta_p, \alpha_q \}, \quad  \{\beta_p, \beta_q \}.
 \end{equation*}
 In Fig. \ref{Dlak_orbits_16c} we depict only one from $4$ options for $D_l(a_k)$ and
 its linkage system $\mathscr{L}(D_l(a_k))$.

%% file: 3exceptional.tex
\newpage
~\\

\subsection{Exceptional semi-Coxeter orbits}
 \label{sec_exc_cases}
~\\

\begin{table}[H]
  \centering
  \scriptsize
  \renewcommand{\arraystretch}{1.5}  
  \begin{tabular} {||c|c|c|c|c|c|c|c|c||}
  \hline  \hline  
        & The Carter & Total  & \multicolumn{6}{c||}{Representatives of exceptional orbits} \\
        &  diagram $\Gamma$ & orbits & \multicolumn{6}{c||}{in the linkage system $\mathscr{L}(\Gamma)$} \cr
     \cline{2-7}
  \hline  \hline  
   1 & $\begin{array}{c} ~\\ \includegraphics[scale=0.5]{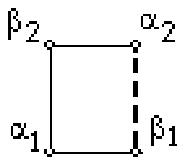} \\ {\bf D_4(a_1)} \end{array}$ & $6$ &
   \multicolumn{6}{c||}{
   $\begin{array}{c} (1a)  \quad \includegraphics[scale=0.6]{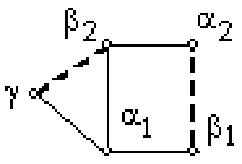} \\ $\{-1, 0, 0, 1 \}$ \end{array}$
    \qquad \qquad
   $\begin{array}{c} (1b)  \quad  \includegraphics[scale=0.6]{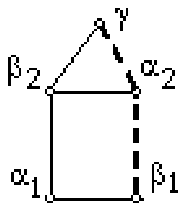} \\ $\{0, 1, 0, -1 \}$ \end{array}$ \quad
   } \\
  \hline  
   2 & $\begin{array}{c} ~\\ \includegraphics[scale=0.5]{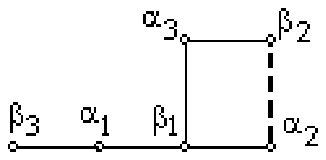} \\ {\bf D_6(a_1)} \end{array}$ & $10$ &
   \multicolumn{6}{c||}{
   $\begin{array}{c} (2a) ~\\ \includegraphics[scale=0.6]{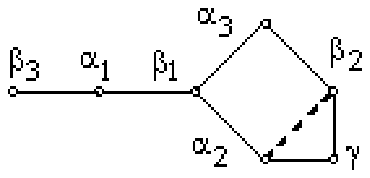} \\ $\{0, -1, 0, 0, -1, 0 \}$ \end{array}$ \quad
   $\begin{array}{c} (2b) ~\\ \includegraphics[scale=0.6]{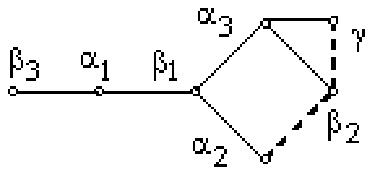} \\ $\{0, 0, -1, 0, 1, 0 \}$ \end{array}$ \quad
   $\begin{array}{c} (2c) ~\\ \includegraphics[scale=0.6]{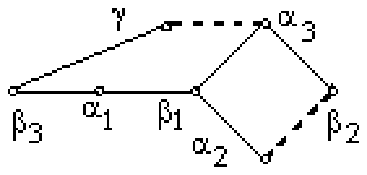} \\ $\{0, 0, 1, 0, 0, -1 \}$ \end{array}$ \quad
   $\begin{array}{c} (2d) ~\\ \includegraphics[scale=0.6]{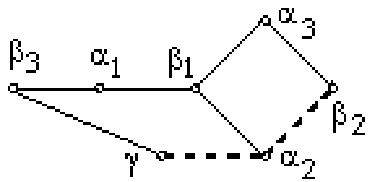} \\ $\{0, 1, 0, 0, 0, -1 \}$ \end{array}$ \quad
   } \\
  \hline  
   3 & $\begin{array}{c} ~\\ \includegraphics[scale=0.5]{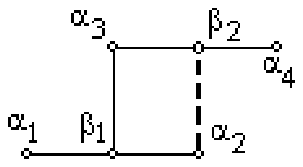} \\ {\bf D_6(a_2)} \end{array}$ & $14$ &
   \multicolumn{6}{c||}{
   $\begin{array}{c} (3a) ~\\ \includegraphics[scale=0.6]{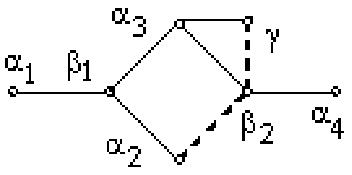} \\ $\{0, 0, -1, 0, 0, 1 \}$ \end{array}$ \quad
   $\begin{array}{c} (3b) ~\\ \includegraphics[scale=0.6]{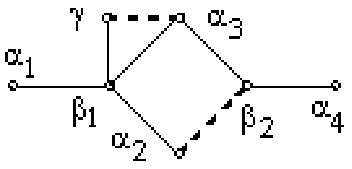} \\ $\{0, 0, 1, 0, -1, 0 \}$ \end{array}$ \quad
   $\begin{array}{c} (3c) ~\\ \includegraphics[scale=0.6]{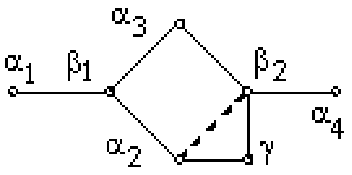} \\ $\{0, -1, 0, 0, 0, -1 \}$ \end{array}$ \quad
   $\begin{array}{c} (3d) ~\\ \includegraphics[scale=0.6]{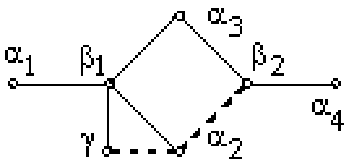} \\ $\{0, 1, 0, 0, -1, 0 \}$ \end{array}$ \quad
   } \\
  \hline  
   4 & $\begin{array}{c} ~\\ \includegraphics[scale=0.5]{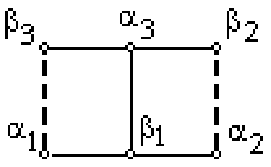} \\ {\bf E_6(a_2)} \end{array}$ & $10$ &
   \multicolumn{6}{c||}{
   $\begin{array}{c} (4a) \quad \includegraphics[scale=0.6]{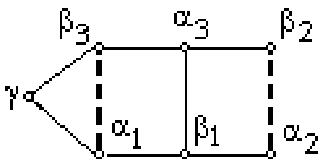} \\ $\{-1, 0, 0, 0, 0, -1 \}$ \end{array}$ \quad
   $\begin{array}{c} (4b) \quad \includegraphics[scale=0.6]{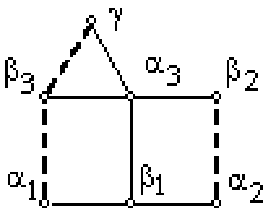} \\ $\{0, 0, -1, 0, 0, 1 \}$ \end{array}$ \quad
   } \\
  \hline  
   5 & $\begin{array}{c} ~\\ \includegraphics[scale=0.5]{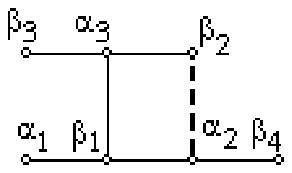} \\ {\bf E_7(a_2)} \end{array}$ & $6$ &
   \multicolumn{6}{c||}{
   $\begin{array}{c} (5a) \quad \includegraphics[scale=0.6]{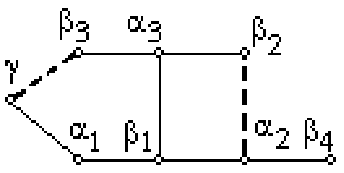} \\ $\{-1, 0, 0, 0, 0, 1, 0 \}$ \end{array}$ } \\
  \hline  
   6 & $\begin{array}{c} ~\\ \includegraphics[scale=0.5]{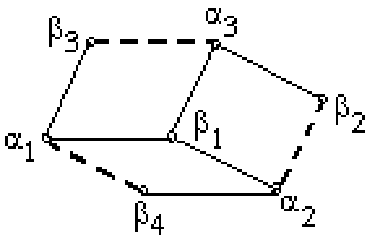} \\ {\bf E_7(a_4)} \end{array}$ & $10$ &
   \multicolumn{6}{c||}{
   $\begin{array}{c} (6a) ~\\ \includegraphics[scale=0.5]{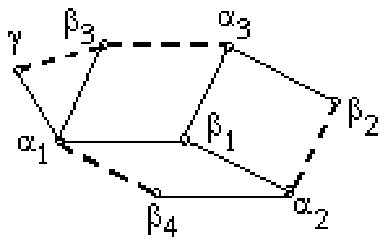} \\
      $\{-1, 0, 0, 0, 0, 1, 0 \}$ \end{array}$ \quad
   $\begin{array}{c} (6b) ~\\ \includegraphics[scale=0.5]{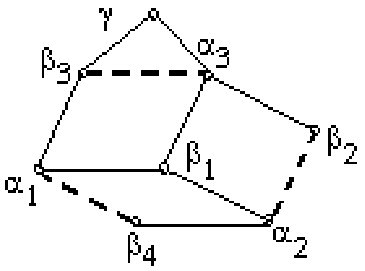} \\
      $\{0, 0, -1, 0, 0, -1, 0 \}$ \end{array}$ \quad
   $\begin{array}{c} (6c) ~\\ \includegraphics[scale=0.5]{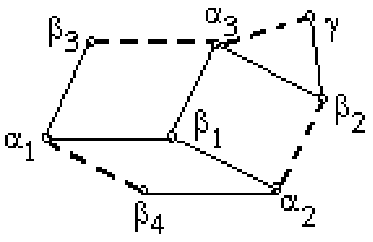} \\
       $\{0, 0, 1, 0, -1, 0, 0 \}$ \end{array}$ } \cr
     & &  &
   \multicolumn{6}{c||}{
   $\begin{array}{c} (6d) ~\\ \includegraphics[scale=0.5]{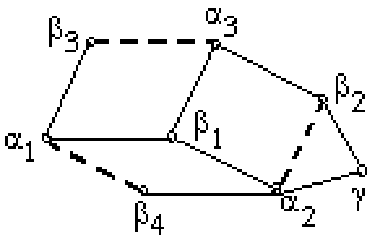} \\
      $\{0, -1, 0, 0, -1, 0, 0 \}$ \end{array}$ \quad
   $\begin{array}{c} (6e) ~\\ \includegraphics[scale=0.5]{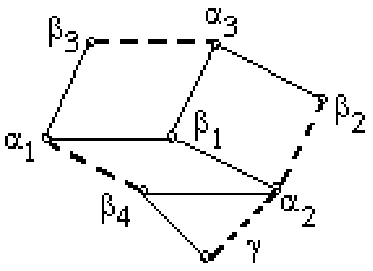} \\
      $\{0, -1, 0, 0, 0, 0, 1 \}$ \end{array}$ \quad
   $\begin{array}{c} (6f) ~\\ \includegraphics[scale=0.5]{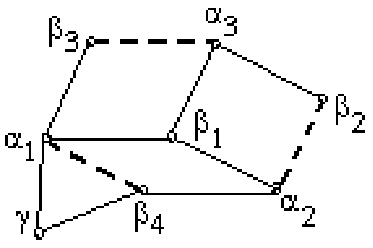} \\
      $\{1, 0, 0, 0, 0, 0, 1 \}$ \end{array}$ } \\
  \hline 
  7  & $\begin{array}{c} ~\\ \includegraphics[scale=0.5]{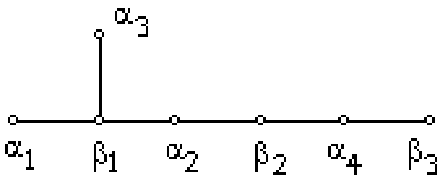} \\ {\bf D_7} \end{array}$ & $14$ &
   \multicolumn{6}{c||}{
   $\begin{array}{c} (7a) \quad \includegraphics[scale=0.6]{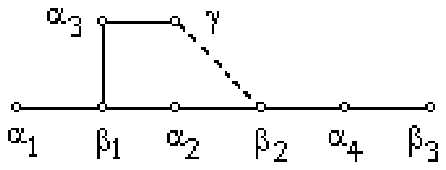} \\
      $\{0, 0, -1, 0, 0, 1, 0 \}$ \end{array}$ \quad
   $\begin{array}{c} (7b) \quad \includegraphics[scale=1.1]{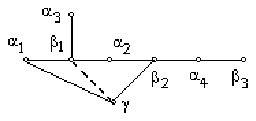} \\
      $\{-1, 0, 0, 0, 1, -1, 0 \}$ \end{array}$ }  \\
   \hline
\end{tabular}
  \vspace{2mm}
  \caption{\hspace{3mm}\footnotesize Exceptional orbits}
  \label{tab_repres_exc}
\end{table}

\newpage
~\\

 \subsubsection{Diagram $D_4(a_1)$. Case $(1a)$} For the Carter diagram $D_4(a_1)$, there are $2$  exceptional
  semi-Coxeter orbits with representatives $(1a)$ and $(1b)$. These cases are similar, see Table \ref{tab_repres_exc}.  We consider only $(1a)$.

  \underline{Case $(1a)$.}

  \begin{equation*}
    \begin{split}
     w = & s_{\alpha_1}s_{\alpha_2}s_{\beta_1}s_{\beta_2}s_{\gamma} =
        s_{\alpha_2}s_{\beta_2}s_{\beta_1}s_{\beta_1 + \beta_2 + \alpha_1}s_{\gamma}
          \stackrel{s_{\alpha_2}}{\simeq}
         (s_{\beta_2}s_{\beta_1})(s_{\beta_1 + \beta_2 + \alpha_1}s_{\gamma}s_{\alpha_2}).
    \end{split}
  \end{equation*}

\begin{figure}[h]
\centering
\includegraphics[scale=0.8]{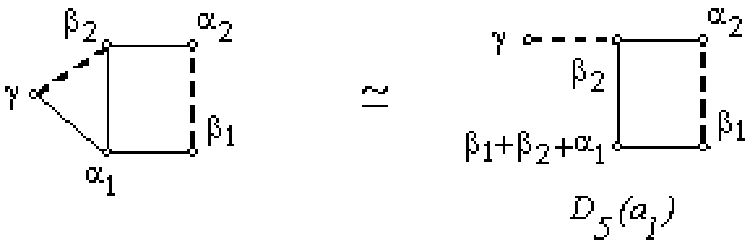}
\caption{\hspace{3mm} \footnotesize The linkage diagram $(1a)$ for the Carter diagram {\bf $D_4(a_1)$}}
\label{D4a1_exception_a}
\end{figure}
  Thus, the linkage diagram $(1a)$ from Table \ref{tab_repres_exc} is equivalent to the Carter diagram $D_5(a_1)$. This case was also considered in \cite[Lemma 1.8]{St10.I}.
~\\

  \subsubsection{Diagram $D_6(a_1)$. Cases $(2a)$ and $(2c)$}
   For the diagram $D_6(a_1)$, there are $4$ exceptional semi-Coxeter orbits $(2a)$, $(2b)$, $(2c)$, $(2d)$.
   Cases $(2a)$ and $(2b)$ are similar; cases $(2c)$ and $(2d)$ are also similar.
   We consider only $(2a)$ and $(2c)$.

  \underline{Case $(2a)$}. Here, we have
  \begin{equation*}
    \begin{split}
     w = & s_{\alpha_1}s_{\alpha_2}s_{\alpha_3}s_{\beta_1}s_{\beta_2}s_{\beta_3}s_{\gamma} =
        s_{\alpha_1}(s_{\alpha_2}s_{\alpha_3}s_{\beta_2})s_{\beta_1}s_{\beta_3}s_{\gamma} = \\
         &  s_{\alpha_1}s_{\beta_2 + \alpha_3 - \alpha_2}s_{\alpha_2}s_{\alpha_3}s_{\beta_1}s_{\beta_3}s_{\gamma}
       \stackrel{s_{\beta_2 + \alpha_3 - \alpha_2}}{\simeq}
       (s_{\alpha_1}s_{\alpha_2}s_{\alpha_3})
        (s_{\beta_1}s_{\beta_3}s_{\gamma}s_{\beta_2 + \alpha_3 - \alpha_2}).
    \end{split}
  \end{equation*}

\begin{figure}[h]
\centering
\includegraphics[scale=0.8]{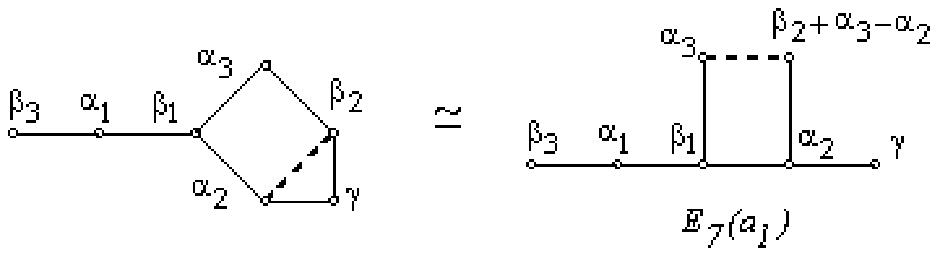}
\caption{\hspace{3mm} \footnotesize The linkage diagram $(2a)$ for the Carter diagram {\bf $D_6(a_1)$}}
\label{D6a1_exception_a}
\end{figure}
 Hence, the linkage diagram $(2a)$ from Table \ref{tab_repres_exc} is equivalent to the Carter diagram $E_7(a_1)$.
 ~\\

 \underline{Case $(2c)$}. This case is reduced to the exception case $(4a)$ in the exceptional orbit for $E_6(a_2)$:

 \begin{equation*}
    \begin{split}
     w = & s_{\alpha_1}s_{\alpha_2}s_{\alpha_3}s_{\beta_1}s_{\beta_2}s_{\beta_3}s_{\gamma} =
        (s_{\alpha_1}s_{\alpha_2}s_{\alpha_3})(s_{\beta_1}s_{\beta_2}s_{\beta_3 + \gamma}s_{\beta_3}) =
         (s_{\alpha_1}s_{\alpha_2}s_{\alpha_3})(s_{\beta_1}s_{\beta_2}s_{-(\beta_3 + \gamma)}s_{\beta_3}) = \\
         & (s_{\alpha_1}s_{\alpha_2}s_{\alpha_3})s_{-(\beta_3 + \gamma)}(s_{\beta_1}s_{\beta_2}s_{\beta_3})
          ~ \simeq ~
          s_{-(\beta_3 + \gamma)}(s_{\beta_1}s_{\beta_2}s_{\beta_3})(s_{\alpha_1}s_{\alpha_2}s_{\alpha_3}).
    \end{split}
  \end{equation*}

  \begin{figure}[h]
\centering
\includegraphics[scale=0.8]{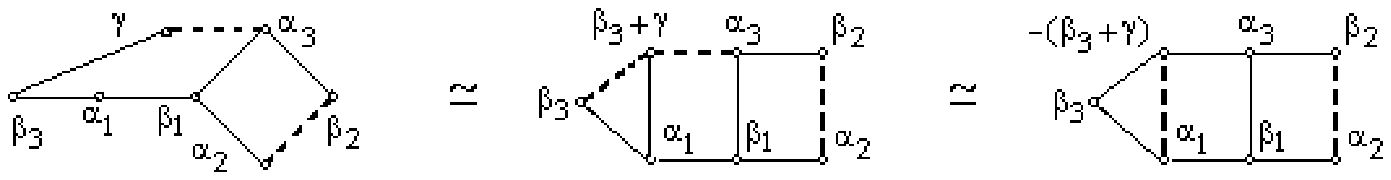}
\caption{\hspace{3mm} \footnotesize The linkage diagram $(2c)$  for the Carter diagram  {\bf $D_6(a_1)$}}
\label{D6a1_exception_c}
\end{figure}
~\\

  \subsubsection{Diagram $D_6(a_2)$. Case $(3a)$}
  For the Carter diagram $D_6(a_2)$, there are $4$ exceptional
  semi-Coxeter orbits with representatives $(3a)$, $(3b)$, $(3c)$, $(3d)$.
  These cases are similar to each other, see Table \ref{tab_repres_exc}. We consider only $(3a)$.

  \underline{Case $(3a)$}.  Here, we have

  \begin{equation*}
    \begin{split}
     w = & s_{\alpha_1}s_{\alpha_2}s_{\alpha_3}s_{\alpha_4}s_{\beta_1}s_{\beta_2}s_{\gamma} =
        s_{\alpha_1}s_{\beta_2 + \alpha_3 - \alpha_2 + \alpha_4}
        s_{\alpha_2}s_{\alpha_3}s_{\alpha_4}s_{\beta_1} s_{\gamma} = \\
        & \\
        & s_{\beta_2 + \alpha_3 - \alpha_2 + \alpha_4}
        (s_{\alpha_1}s_{\alpha_2}s_{\alpha_3}s_{\alpha_4})s_{\beta_1} s_{\gamma}
        \stackrel{s_{\beta_2 + \alpha_3 - \alpha_2 + \alpha_4}}{\simeq}
         (s_{\alpha_1}s_{\alpha_2}s_{\alpha_3}s_{\alpha_4})
         (s_{\beta_1} s_{\gamma}s_{\beta_2 + \alpha_3 - \alpha_2 + \alpha_4}) = \\
        & \\
        &   (s_{\alpha_1}s_{-\alpha_2}s_{\alpha_3}s_{\alpha_4})
         (s_{\beta_1} s_{\gamma}s_{-(\beta_2 + \alpha_3 - \alpha_2 + \alpha_4)}).
    \end{split}
  \end{equation*}
  Therefore, the linkage diagram $(3a)$ from Table \ref{tab_repres_exc} is equivalent to the Carter diagram $E_7(a_2)$.
\begin{figure}[h]
\centering
\includegraphics[scale=0.8]{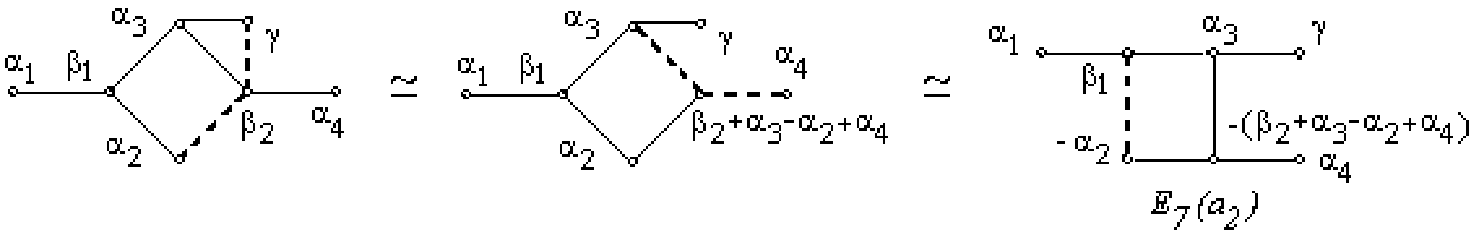}
\caption{\hspace{3mm} \footnotesize The linkage diagram $(3a)$  for the Carter diagram {\bf $D_6(a_2)$}}
\label{D6a2_exception_a}
\end{figure}

\begin{remark}[on $s$-permutation] {\rm
  Let $s_\alpha$ and $s_\beta$ be two adjacent reflections in the decomposition of $w$:
  \begin{equation}
    \label{perm_1}
     w = \dots s_\alpha s_\beta \dots.
  \end{equation}
  If decomposition \eqref{perm_1} is written down in one of the equivalent forms
\begin{equation}
    \label{perm_2}
   \begin{split}
     & w = \dots s_\alpha s_\beta \dots = \dots s_{s_{\alpha}(\beta)}s_\alpha  \dots, \text{ or } \\
     & w = \dots s_\alpha s_\beta \dots = \dots s_\beta s_{s_{\beta}(\alpha)}  \dots,
   \end{split}
  \end{equation}
  we say that elements $s_\alpha$ and $s_\beta$ are {\it $s$-permuted}.
  The linkage diagram related to \eqref{perm_1} is respectively changed.
   The corresponding transformation of the word $w$ and related linkage diagram we call
   the {\it $s$-permutation}. In \cite[\S1.4.1]{St10.I}, we considered $s$-permutation in the framework of equivalent transformations of connection diagrams. \qed
  }
\end{remark}

\subsubsection{Diagram $E_6(a_2)$. Cases $(4a)$, $(4b)$}
 These two cases are different.

\begin{figure}[h]
\centering
\includegraphics[scale=0.8]{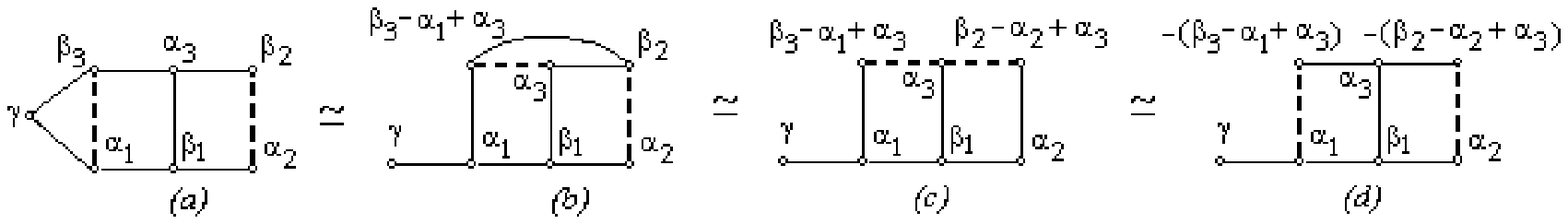}
\caption{\hspace{3mm} \footnotesize The linkage diagram $(4a)$ for the Carter diagram {\bf $E_6(a_2)$}}
\label{E6a2_exception_a}
\end{figure}

   \underline{Case $(4a)$}. First, reflections $s_{\beta_3}$ and $s_{\alpha_1}s_{\alpha_3}$
 are $s$-permuted. The new connection between $\beta_2$ and $\beta_3  - \alpha_1 + \alpha_3$
  appears, see Fig. \ref{E6a2_exception_a},(b).
  After that,  $s_{\beta_2}$ and $s_{\alpha_2}s_{\alpha_3}$ are $s$-permuted.
  The the new connection disappears, see Fig. \ref{E6a2_exception_a},(c).
 \begin{equation*}
     w = s_{\alpha_1}s_{\alpha_2}s_{\alpha_3}s_{\beta_1}s_{\beta_2}s_{\beta_3}s_{\gamma} =
        s_{\alpha_2}s_{\beta_3  - \alpha_1 + \alpha_3}
        s_{\alpha_1}s_{\alpha_3}s_{\beta_1}s_{\beta_2}s_{\gamma} =
         s_{\beta_3  - \alpha_1 + \alpha_3}
        s_{\alpha_1}s_{\alpha_2}s_{\alpha_3}s_{\beta_1}s_{\beta_2}s_{\gamma},
  \end{equation*}
      since $s_{\alpha_2}$ and $s_{\beta_3  - \alpha_1 + \alpha_3}$ commute, see Fig. \ref{E6a2_exception_a},(b).
 \begin{equation*}
     w = s_{\beta_3  - \alpha_1 + \alpha_3}s_{\alpha_1}s_{\beta_2  - \alpha_2 + \alpha_3}
         s_{\alpha_2}s_{\alpha_3}s_{\beta_1}s_{\gamma} =
         s_{\beta_3  - \alpha_1 + \alpha_3}s_{\beta_2  - \alpha_2 + \alpha_3}
         s_{\alpha_1}s_{\alpha_2}s_{\alpha_3}s_{\beta_1}s_{\gamma},
  \end{equation*}
    since $s_{\alpha_1}$ and $s_{\beta_2  - \alpha_2 + \alpha_3}$ commute, see Fig. \ref{E6a2_exception_a},(c).
    Further,
 \begin{equation*}
     w = s_{\beta_3  - \alpha_1 + \alpha_3}s_{\beta_2  - \alpha_2 + \alpha_3}
         s_{\alpha_1}s_{\alpha_2}s_{\alpha_3}s_{\beta_1}s_{\gamma} =
         s_{-(\beta_3  - \alpha_1 + \alpha_3)}s_{-(\beta_2  - \alpha_2 + \alpha_3)}
         s_{\alpha_1}s_{\alpha_2}s_{\alpha_3}s_{\beta_1}s_{\gamma},
  \end{equation*}
    see Fig. \ref{E6a2_exception_a},(d). Hence, the linkage diagram $(4a)$ from Table \ref{tab_repres_exc} is equivalent to the Carter diagram $E_7(a_3)$.

\begin{figure}[h]
\centering
\includegraphics[scale=0.8]{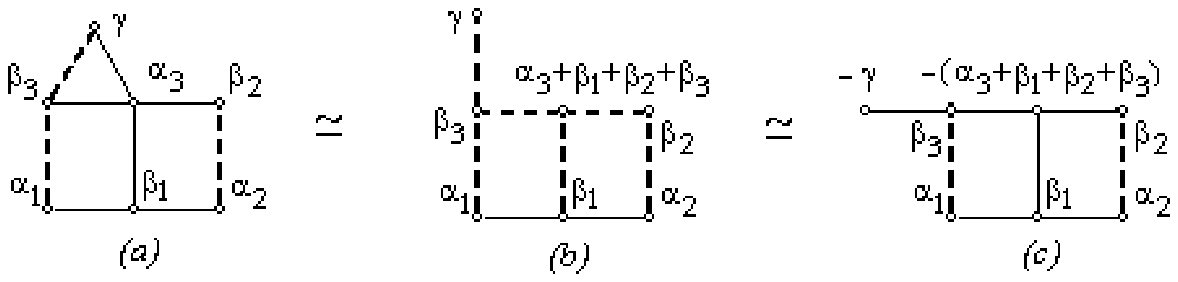}
\caption{\hspace{3mm} \footnotesize The linkage diagram $(4b)$  for the Carter diagram {\bf $E_6(a_2)$}}
\label{E6a2_exception_b}
\end{figure}

    \underline{Case $(4b)$}.

 \begin{equation*}
     w = s_{\alpha_1}s_{\alpha_2}(s_{\alpha_3}s_{\beta_1}s_{\beta_2}s_{\beta_3})s_{\gamma} =
         s_{\alpha_1}s_{\alpha_2}s_{\beta_1}s_{\beta_2}s_{\beta_3}s_{\alpha_3  + \beta_1 + \beta_2 + \beta_3}
         s_{\gamma},
  \end{equation*}
where
 \begin{equation*}
  \begin{split}
 & (\alpha_3  + \beta_1 + \beta_2 + \beta_3, \alpha_1) =
      (\beta_3, \alpha_1) + (\beta_1, \alpha_1) = 0,  \\
 & (\alpha_3  + \beta_1 + \beta_2 + \beta_3, \alpha_2) =
      (\beta_2, \alpha_2) + (\beta_1, \alpha_2) =  0.
  \end{split}
  \end{equation*}
  Thus, $w$ is described by the diagram Fig. \ref{E6a2_exception_b},(b). Further,
  \begin{equation*}
     w = s_{\alpha_1}s_{\alpha_2}s_{\beta_1}s_{\beta_2}s_{\beta_3}s_{\alpha_3  + \beta_1 + \beta_2 + \beta_3}
         s_{\gamma} =
         s_{\alpha_1}s_{\alpha_2}s_{\beta_1}s_{\beta_2}s_{\beta_3}s_{-(\alpha_3  + \beta_1 + \beta_2 + \beta_3)}
         s_{\gamma},
  \end{equation*}
  see Fig. \ref{E6a2_exception_b},(c). Hence, the linkage diagram $(4b)$ from Table \ref{tab_repres_exc} is equivalent to the Carter diagram $E_7(a_3)$.

\subsubsection{Diagram $E_7(a_2)$. Case $(5a)$}
~\\

\begin{figure}[h]
\centering
\includegraphics[scale=0.8]{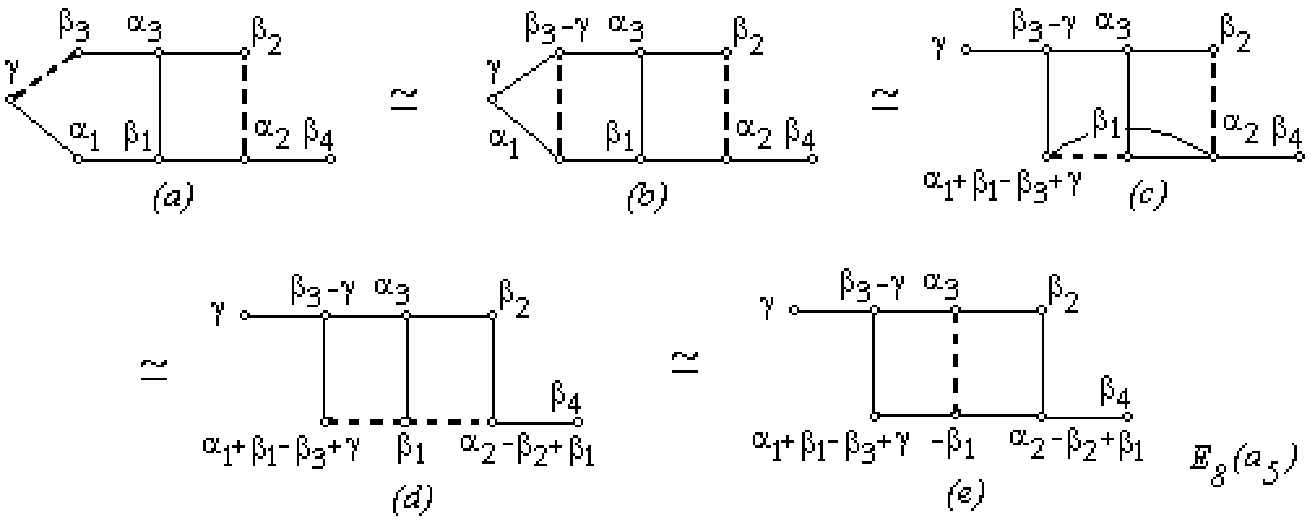}
\caption{\hspace{3mm} \footnotesize The linkage diagram $(5a)$  for the Carter diagram {\bf $E_7(a_2)$}}
\label{E7a2_exception_a}
\end{figure}

\underline{Case $(5a)$}.  First, we reduce the $5$-cycle to the usual contour of the $4$-cycle with the adjoined triangle, see Fig. \ref{E7a2_exception_a},(b).
  \begin{equation*}
     w = s_{\alpha_1}s_{\alpha_2}s_{\alpha_3}s_{\beta_1}s_{\beta_2}s_{\beta_4}(s_{\beta_3}s_{\gamma}) =
         s_{\alpha_1}s_{\alpha_2}s_{\alpha_3}s_{\beta_1}s_{\beta_2}s_{\beta_4}s_{\gamma}s_{\beta_3 - \gamma}.
  \end{equation*}
 Secondly, reflections $s_{\alpha_1}$  and $s_{\beta_1}s_{\beta_3 - \gamma}$ are $s$-permuted. The new connection
 between $\alpha_1 + \beta_1 - \beta_3  + \gamma$ and $\alpha_2$ appears:
  \begin{equation*}
    \begin{split}
     w = & s_{\alpha_1}s_{\alpha_2}s_{\alpha_3}s_{\beta_2}s_{\beta_4}s_{\gamma}s_{\beta_1}s_{\beta_3 - \gamma}
         \stackrel{s_{\beta_1}s_{\beta_3 - \gamma}}{\simeq}
         (s_{\beta_1}s_{\beta_3 - \gamma}s_{\alpha_1})s_{\alpha_2}s_{\alpha_3}s_{\beta_2}s_{\beta_4}s_{\gamma} = \\
         & s_{\alpha_1 + \beta_1 - \beta_3  + \gamma}s_{\beta_1}s_{\beta_3 - \gamma}
           s_{\alpha_2}s_{\alpha_3}s_{\beta_2}s_{\beta_4}s_{\gamma},
    \end{split}
  \end{equation*}
  where
 \begin{equation*}
  \begin{split}
   & (\alpha_1 + \beta_1 - \beta_3  + \gamma, \gamma) =
      (\gamma, \gamma) - (\beta_3, \gamma) + (\alpha_1, \gamma) = 1 - \frac{1}{2} - \frac{1}{2} = 0,  \\
   & (\alpha_1 + \beta_1 - \beta_3  + \gamma, \alpha_3) = (\beta_1 - \beta_3, \alpha_3) = 0, \\
   & (\alpha_1 + \beta_1 - \beta_3  + \gamma, \alpha_2) = (\beta_1, \alpha_2) = - \frac{1}{2},
  \end{split}
  \end{equation*}
  see Fig. \ref{E7a2_exception_a},(c).

~\\
   Further, reflections $\alpha_2$ and $s_{\beta_1}s_{\beta_2}$ are $s$-permuted.
   Then the new connection disappears:
 \begin{equation*}
    \begin{split}
     w = & s_{\alpha_1 + \beta_1 - \beta_3  + \gamma}s_{\beta_1}s_{\beta_3 - \gamma}
           s_{\alpha_2}s_{\alpha_3}s_{\beta_2}s_{\beta_4}s_{\gamma} \stackrel{s_{\beta_2}}{\simeq}
           s_{\alpha_1 + \beta_1 - \beta_3  + \gamma}s_{\beta_3 - \gamma}(s_{\beta_2}s_{\beta_1}s_{\alpha_2})
               s_{\alpha_3}s_{\beta_4}s_{\gamma} = \\
         & s_{\alpha_1 + \beta_1 - \beta_3  + \gamma}s_{\beta_3 - \gamma}s_{\alpha_2 - \beta_2 + \beta_1}
           s_{\beta_2}s_{\beta_1}s_{\alpha_3}s_{\beta_4}s_{\gamma},
    \end{split}
  \end{equation*}
    where
 \begin{equation*}
    (\alpha_1 + \beta_1 - \beta_3  + \gamma, \alpha_2 - \beta_2 + \beta_1) =
    (\alpha_2, \beta_1) + (\beta_1, \beta_1) + (\beta_1, \alpha_1) = -\frac{1}{2} + 1 - \frac{1}{2} = 0.  \\
  \end{equation*}

    Further,
 \begin{equation*}
    \begin{split}
     w = & s_{\alpha_1 + \beta_1 - \beta_3  + \gamma}s_{\beta_3 - \gamma}s_{\alpha_2 - \beta_2 + \beta_1}
           s_{\beta_2}s_{\beta_1}s_{\alpha_3}s_{\beta_4}s_{\gamma} =
           s_{\alpha_1 + \beta_1 - \beta_3  + \gamma}s_{\alpha_2 - \beta_2 + \beta_1}
           s_{\beta_3 - \gamma}s_{\beta_2}s_{\beta_1}s_{\alpha_3}s_{\beta_4}s_{\gamma}
           \stackrel{s_{\alpha_3}}{\simeq} \\
         &  (s_{\alpha_3}s_{\alpha_1 + \beta_1 - \beta_3  + \gamma}s_{\alpha_2 - \beta_2 + \beta_1})
            (s_{\beta_3 - \gamma}s_{\beta_2}s_{\beta_1}s_{\beta_4}s_{\gamma}).
    \end{split}
  \end{equation*}
    see Fig. \ref{E7a2_exception_a},(d).
    In the last step, $s_{\beta_1}$ is replaced by $s_{-\beta_1}$, the corresponding diagram is depicted in
    Fig. \ref{E7a2_exception_a},(e). Thus, the linkage diagram $(5a)$ from Table \ref{tab_repres_exc} is equivalent to the Carter diagram $E_8(a_5)$.
~\\

 \subsubsection{Diagram $E_7(a_4)$. Case $(6a)$}
  For the Carter diagram $E_7(a_4)$, there are $6$ exceptional
  semi-Coxeter orbits with representatives
  $(6a)$, $(6b)$, $(6c)$, $(6d)$, $(6e)$, $(6f)$, that are similar to each other, see Table \ref{tab_repres_exc}.
  Let us consider $(6a)$.

\begin{figure}[h]
\centering
\includegraphics[scale=0.8]{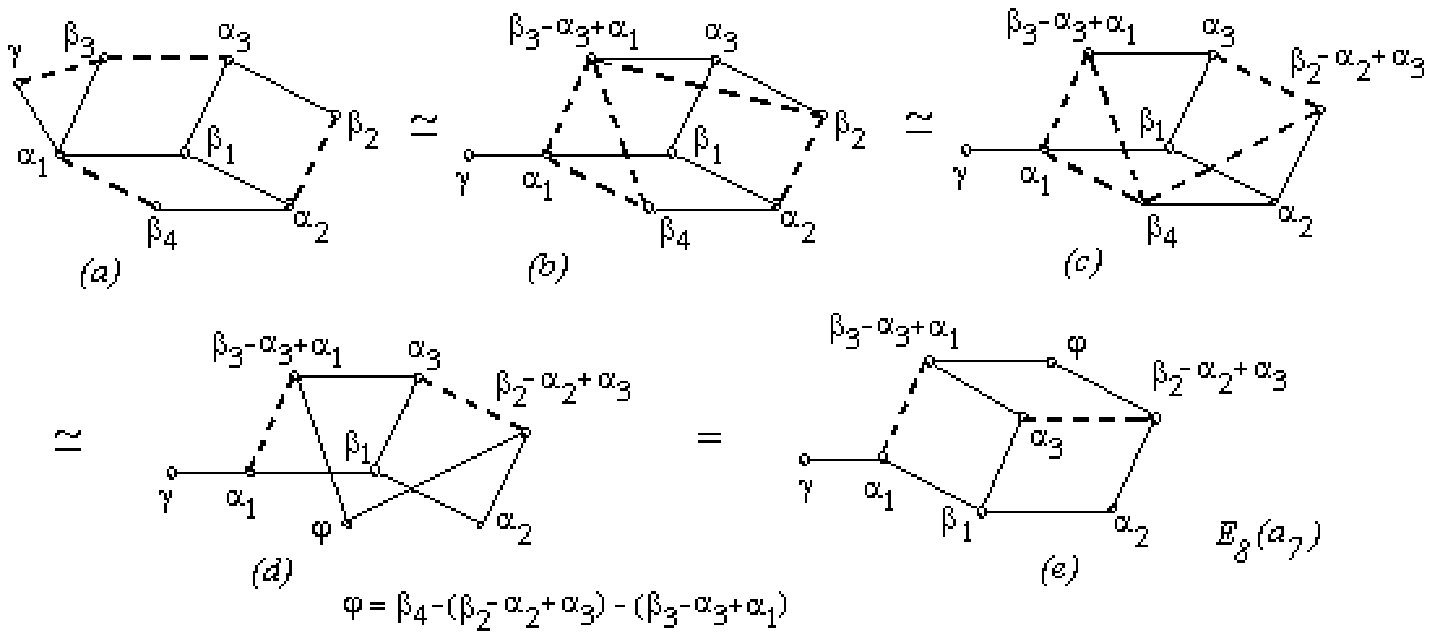}
\caption{\hspace{3mm} \footnotesize The linkage diagram $(6a)$  for the Carter diagram {\bf $E_7(a_4)$}}
\label{E7a4_exception_a}
\end{figure}

  \underline{Case $(6a)$}. First, reflections $s_{\beta_3}$  and $s_{\alpha_1}s_{\alpha_3}$ are $s$-permuted. Two new connections  $\{ \beta_3 - \alpha_3  + \alpha_1, \beta_2 \}$  and
  $\{ \beta_3 - \alpha_3  + \alpha_1, \beta_4 \}$  appear:
  \begin{equation*}
    \begin{split}
    w = & s_{\alpha_2}(s_{\alpha_1}s_{\alpha_3}s_{\beta_3})s_{\beta_1}s_{\beta_2}s_{\beta_4}s_{\gamma} =
         s_{\alpha_2}s_{\beta_3 - \alpha_3 + \alpha_1}s_{\alpha_1}s_{\alpha_3}
         s_{\beta_1}s_{\beta_2}s_{\beta_4}s_{\gamma} = \\
        &  s_{\beta_3 - \alpha_3 + \alpha_1}s_{\alpha_1}s_{\alpha_2}s_{\alpha_3}
         s_{\beta_1}s_{\beta_2}s_{\beta_4}s_{\gamma},
      \end{split}
  \end{equation*}
 since $\alpha_2$ and $\beta_3 - \alpha_3 + \alpha_1$ commute, see Fig. \ref{E7a4_exception_a},(b).
 After that, reflections $s_{\beta_2}$  and $s_{\alpha_2}s_{\alpha_3}$ are $s$-permuted. Then
 the connection $\{ \beta_3 - \alpha_3  + \alpha_1, \beta_2 \}$ disappears, and the connection
 $\{ \beta_2 - \alpha_2  + \alpha_3, \beta_4 \}$ appears:
    \begin{equation*}
    \begin{split}
    w =  &  s_{\beta_3 - \alpha_3 + \alpha_1}s_{\alpha_1}(s_{\alpha_2}s_{\alpha_3}
         s_{\beta_2})s_{\beta_1}s_{\beta_4}s_{\gamma} =
         s_{\beta_3 - \alpha_3 + \alpha_1}s_{\alpha_1}s_{\beta_2 - \alpha_2 + \alpha_3}
         s_{\alpha_2}s_{\alpha_3}s_{\beta_1}s_{\beta_4}s_{\gamma} = \\
        &  s_{\beta_3 - \alpha_3 + \alpha_1}s_{\beta_2 - \alpha_2 + \alpha_3}
         s_{\alpha_1}s_{\alpha_2}s_{\alpha_3}s_{\beta_1}s_{\beta_4}s_{\gamma},
      \end{split}
  \end{equation*}
  since $s_{\alpha_1}$ and $s_{\beta_2 - \alpha_2 + \alpha_3}$ commute,
  see Fig. \ref{E7a4_exception_a},(c). We have
  \begin{equation*}
     (\beta_3 - \alpha_3 + \alpha_1, \beta_2 - \alpha_2 + \alpha_3)  =
     (\beta_3, \alpha_3) - (\alpha_3,\alpha_3) - (\alpha_3, \beta_2) =
     \frac{1}{2} - 1 + \frac{1}{2} = 0.
  \end{equation*}

~\\
  Further, reflections $s_{\beta_3 - \alpha_3 + \alpha_1}s_{\beta_2 - \alpha_2 + \alpha_3}$
  and $s_{\beta_4}$ are $s$-permuted:
   \begin{equation*}
    \begin{split}
    w =  &  s_{\beta_3 - \alpha_3 + \alpha_1}s_{\beta_2 - \alpha_2 + \alpha_3}
         s_{\alpha_1}s_{\alpha_2}s_{\alpha_3}s_{\beta_1}s_{\beta_4}s_{\gamma}
         \stackrel{s_{\beta_4}}{\simeq}  \\
        & (s_{\beta_4}s_{\beta_3 - \alpha_3 + \alpha_1}s_{\beta_2 - \alpha_2 + \alpha_3})
          s_{\alpha_1}s_{\alpha_2}s_{\alpha_3}s_{\beta_1}s_{\gamma} =
          s_{\beta_3 - \alpha_3 + \alpha_1}s_{\beta_2 - \alpha_2 + \alpha_3}
          s_{\varphi}s_{\alpha_1}s_{\alpha_2}s_{\alpha_3}s_{\beta_1}s_{\gamma}, \\
        & \text{ where }
          \varphi = \beta_4 - (\beta_3 - \alpha_3 + \alpha_1)  - (\beta_2 - \alpha_2 + \alpha_3) =
              \beta_4 - \beta_3 - \beta_2 - \alpha_1 + \alpha_2.
      \end{split}
  \end{equation*}
  Then
  \begin{equation*}
   \begin{split}
     & (\varphi, \alpha_1)  =
         (\beta_4 - \beta_3 - \alpha_1, \alpha_1) = \frac{1}{2} + \frac{1}{2} - 1  = 0, \\
     & (\varphi, \alpha_2)  =
         (\beta_4 - \beta_2 + \alpha_2, \alpha_2)  = -\frac{1}{2} - \frac{1}{2} + 1 = 0,
    \end{split}
  \end{equation*}
  i.e.,  connections $\{ \varphi, \alpha_1 \}$ and $\{ \varphi, \alpha_2 \}$ disappear,
  see Fig. \ref{E7a4_exception_a},(d). Hence, the linkage diagram $(6a)$ from Table \ref{tab_repres_exc} is equivalent to the Carter diagram $E_8(a_7)$, see Fig. \ref{E7a4_exception_a},(e).  \qed

    \subsubsection{Diagram $D_7$. Cases $(7a)$, $(7b)$}
     These $2$ cases are different.
~\\

    \underline{Case $(7a)$}.  We have
    \begin{equation*}
     w =  s_{\gamma}s_{\alpha_1}s_{\alpha_2}s_{\alpha_3}s_{\alpha_4}s_{\beta_1}s_{\beta_2}s_{\beta_3} =
       s_{\alpha_3}(s_{\gamma + \alpha_3}s_{\alpha_1}s_{\alpha_2}s_{\alpha_4}s_{\beta_1})
         s_{\beta_2}s_{\beta_3}.
   \end{equation*}
   Let $s$-permute $s_{\gamma + \alpha_3}s_{\alpha_1}s_{\alpha_2}s_{\alpha_4}$ and $s_{\beta_1}$:
    \begin{equation*}
      \begin{split}
     w =  & s_{\alpha_3}s_{\beta_1 + \alpha_1 + \alpha_2 + \gamma + \alpha_3}
      (s_{\gamma + \alpha_3}s_{\alpha_1}s_{\alpha_2}s_{\alpha_4})s_{\beta_2}s_{\beta_3}
        \stackrel{s_{\beta_2}s_{\beta_3}}{\simeq}  \\
          & (s_{\beta_2}s_{\beta_3}s_{\alpha_3}s_{\beta_1 + \alpha_1 + \alpha_2 + \gamma + \alpha_3})
            (s_{\gamma + \alpha_3}s_{\alpha_1}s_{\alpha_2}s_{\alpha_4}).
     \end{split}
   \end{equation*}

\begin{figure}[h]
\centering
\includegraphics[scale=1.4]{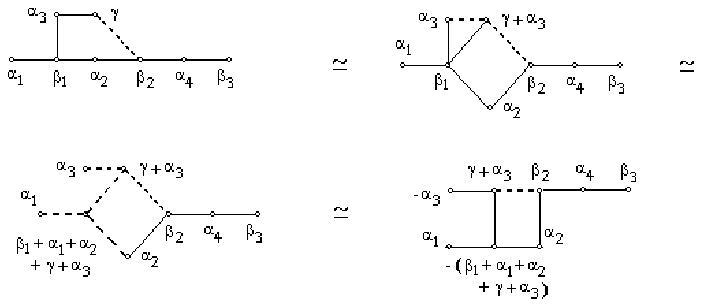}
\caption{\hspace{3mm} \footnotesize The linkage diagram $(7a)$  for the Carter diagram {\bf $D_7$}}
\label{D7_exception_a}
\end{figure}

    Since
    \begin{equation*}
      \begin{split}
        & (\beta_1 + \alpha_1 + \alpha_2 + \gamma + \alpha_3, \alpha_3) =
          (\beta_1 + \gamma + \alpha_3, \alpha_3) = 1 - \frac{1}{2} - \frac{1}{2} = 0, \\
        & (\beta_1 + \alpha_1 + \alpha_2 + \gamma + \alpha_3, \beta_2) =  (\gamma + \alpha_2, \beta_2)  =
             \frac{1}{2} - \frac{1}{2} = 0,
     \end{split}
   \end{equation*}
    we get the bicolored decomposition, see Fig. \ref{D7_exception_a}.
    Since  $s_{\beta_1 + \alpha_1 + \alpha_2 + \gamma + \alpha_3} = s_{-(\beta_1 + \alpha_1 + \alpha_2 + \gamma + \alpha_3)}$ and $s_{\alpha_3} = s_{-\alpha_3}$ we can change $\beta_1 + \alpha_1 + \alpha_2 + \gamma + \alpha_3$
to the opposite vector $-(\beta_1 + \alpha_1 + \alpha_2 + \gamma + \alpha_3)$ and $\alpha_3$ to $-\alpha_3$.
Thus, we obtain the last diagram in Fig. \ref{D7_exception_a}. The dotted edge $\{\gamma + \alpha_3, \beta_2 \}$
(i.e., the property \lq\lq{be dotted edge\rq\rq)
can be moved to any edge of the square. Hence, we get the Carter diagram $E_8(a_2)$.

    \underline{Case $(7b)$}. Let $s$-permute $s_\gamma{s}_{\alpha_2}$ and $s_{\beta_1}$ as follows:

\begin{equation*}
   \begin{split}
     w =  & s_{\gamma}s_{\alpha_1}s_{\alpha_2}s_{\alpha_3}s_{\alpha_4}s_{\beta_1}s_{\beta_2}s_{\beta_3}
     \stackrel{s_{\gamma}s_{\alpha_2}}{\simeq}
       (s_{\alpha_1}s_{\alpha_3}s_{\alpha_4})s_{\beta_2}s_{\beta_3}
       (s_{\beta_1}s_{\gamma}s_{\alpha_2}) = \\
        & s_{\alpha_1}s_{\alpha_3}s_{\alpha_4}s_{\beta_2}s_{\beta_3}s_{\gamma}s_{\alpha_2}
          s_{\beta_1 - \gamma + \alpha_2}
     \stackrel{s_{\beta_1 - \gamma + \alpha_2}}{\simeq}
       s_{\alpha_4}(s_{\beta_1 - \gamma + \alpha_2}s_{\alpha_1}s_{\beta_2}s_{\beta_3})
        (s_{\alpha_3}s_{\gamma}s_{\alpha_2}) \stackrel{s_{\alpha_4}}{\simeq} \\
        & (s_{\beta_1 - \gamma + \alpha_2}s_{\alpha_1}s_{\beta_2}s_{\beta_3})
          (s_{\alpha_3}s_{\gamma}s_{\alpha_2}s_{\alpha_4}),
   \end{split}
\end{equation*}
    where
    \begin{equation*}
      \begin{split}
        & (\beta_1 - \gamma + \alpha_2, \alpha_1) =
          (\beta_1 - \gamma, \alpha_1) = -\frac{1}{2} + \frac{1}{2} = 0, \\
        & (\beta_1 - \gamma + \alpha_2, \beta_2) =  (-\gamma + \alpha_2, \beta_2)  =
             \frac{1}{2} - \frac{1}{2} = 0.
     \end{split}
   \end{equation*}
\begin{figure}[h]
\centering
\includegraphics[scale=1.4]{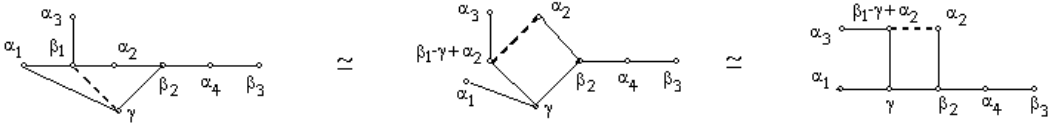}
\caption{\hspace{3mm} \footnotesize The linkage diagram $(7a)$ for the Carter diagram {\bf $D_7$}}
\label{D7_exception_b}
\end{figure}

  Thus, as in the Case (7a), $w$ has the bicolored decomposition corresponding to the Carter diagram $E_8(a_2)$,
  see Fig. \ref{D7_exception_b}.

  \qed

%% file: 4finding.tex
\subsection{Finding semi-Coxeter orbits}

\subsubsection{How to find semi-Coxeter orbits?}
We use two ways to find semi-Coxeter orbits. The first one is the matrix approach:

 1) Calculation of powers of dual semi-Coxeter elements $\SC^*$, see Tables \ref{tab_partial Cartan_1} - \ref{tab_partial Cartan_3}.

 2) Applying $(\SC^*)^k$ to any unicolored linkage diagram $\gamma^{\vee}$ until finding the period of $\SC^*$ on
    this linkage.

 3) Search any new linkage from the corresponding linkage system, preferably unicolored and back to step 2).
   The linkage systems for all Carter diagrams from $\mathsf{DE4}$ and $\mathsf{C4}$ are in \cite{St10.II}.
\\

   The second way is the diagram approach. We find all semi-Coxeter orbits as a closed cycles in
   the linkage system.  We call such a semi-Coxeter orbit the {\it $\SC^*$-cycle}.
   The link connecting $\gamma^{\vee}$ and $\SC\gamma^{\vee}$ we call the {\it $\SC^*$-transition}.
   Every $\SC^*$-cycle consists of $\SC^*$-transitions}
   $\gamma^{\vee} \longrightarrow \SC\gamma^{\vee}$. Each $\SC^*$-transition
   consists of $2$ passages, one after the other:  $\gamma^{\vee} \longrightarrow {}^t{w}_{\beta}\gamma^{\vee}$
   and ${}^t{w}_{\beta}\gamma^{\vee} \longrightarrow {}^t{w}_{\alpha}{}^t{w}_{\beta}\gamma^{\vee}$,
   where
   \begin{equation*}
       {}^t{w}_{\alpha} = \prod\limits_{i=1}^k{s}^*_{\alpha_i}, \quad
       {}^t{w}_{\beta} = \prod\limits_{j=1}^h{s}^*_{\beta_j},
   \end{equation*}
   Reflections ${s}^*_{\beta_j}$ and ${s}^*_{\alpha_i}$ act on the linkage diagrams in the linkage system $\mathscr{L}(\Gamma)$.
   The order of actions of ${s}^*_{\beta_j}$ within ${}^t{w}_{\beta}$ 
   (resp. ${s}^*_{\alpha_i}$ within ${}^t{w}_{\beta}$) does not matter
   since all ${s}^*_{\beta_j}$ (resp.  ${s}^*_{\alpha_i}$) mutually commute. Examples of semi-Coxeter orbits
   are presented in Appendix \ref{sec_orbits}, where the orbits are differed by colors or bold and dotted
   lines.  Let $\mathscr{L}(\Gamma)$ be the linkage system for the Carter diagram of $\Gamma$.
   Note that for any linkage $\gamma^{\vee} \in \mathscr{L}(\Gamma)$, we have
   $-\gamma^{\vee} \in \mathscr{L}(\Gamma)$,
   since $\mathscr{B}^{\vee}_L(\gamma^{\vee}) = \mathscr{B}^{\vee}_L(-\gamma^{\vee})$ and
  \begin{equation*}
      \mathscr{B}^{\vee}_L(\gamma^{\vee}) < 2 \Longleftrightarrow \gamma^{\vee} \in \mathscr{L}(\Gamma),
  \end{equation*}
   see \cite[Theorem 2.14]{St10.II}. Two orbits are said to be the {\it opposite orbits} if for every linkage
   $\gamma^{\vee}$ in one of the orbits there exists the linkage $-\gamma^{\vee}$ in another one. 
   There are some orbits which are opposite to themselves, 
   such an orbit is said to be the {\it self-opposite orbit}.

   For Carter diagrams $D_4$, $D_4(a_1)$, $D_5(a_1)$, $D_5$, $E_6(a_1)$, $E_6(a_2)$, $E_6$,
   the figures of linkage systems with semi-Coxeter orbits are depicted 
   in Fig. \ref{D4_orbits_16c}-\ref{E6pure_orbits_16c}.

%% file: A01matrSemiCox.tex
\newpage
\section{\sc\bf The dual semi-Coxeter element for the Carter diagrams}
 \begin{table}[H]
\tiny
  \centering
  \renewcommand{\arraystretch}{1.5}
  \begin{tabular} {|c|c|c|c|}
  \hline
   The Carter  & The transpose semi-Coxeter   &  The dual semi-Coxeter      & Order  \cr
   diagram     & element $\DualSC$            &  element $\InvDualSC = \DualSC^{-1}$         & of $\DualSC$ \\
  \hline  
       &  & &   \\
     $\begin{array}{c} \includegraphics[scale=0.6]{D4a1.eps} \\
        {\bf D_4(a_1)} \end{array}$  &
     $\begin{array}{c}
       \left [
   \begin{array}{cccc}
     1 & 0 &   1 & 1  \\
     0 & 1 &  -1 & 1 \\
     -1 & 1 &  -1 & 0  \\
     -1 & -1 &  0 & -1  \\
  \end{array}
  \right ]
     \end{array}$
     &
      $\begin{array}{c}
       \left [
   \begin{array}{cccc}
     -1  &  0  & -1  & -1 \\
      0  & -1  &  1  & -1 \\
      1  & -1  &  1  &  0 \\
      1  &  1  &  0  &  1 \\
  \end{array}
  \right ]
     \end{array}$
     & $4$ \\
   & & & \\
    \hline
    &  & &   \\
     $\begin{array}{c} \includegraphics[scale=0.6]{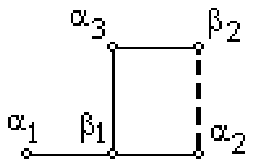} \\
        {\bf D_5(a_1) = D_5(a_2)} \end{array}$ &
     $\begin{array}{c}
       \left [
   \begin{array}{ccccc}
     0 & 1 &  1 & 1 & 0  \\
     1 & 1 &  0 & 1 & -1 \\
     1 & 0 &  1 & 1 & 1  \\
     -1 & -1 &  -1 & -1 & 0  \\
      0 &  1 &   -1  & 0 & -1 \\
  \end{array}
  \right ]
     \end{array}$
     &
      $\begin{array}{c}
       \left [
   \begin{array}{ccccc}
     -1   &  0  &   0  &  -1  &   0 \\
      0   & -1  &   0  &  -1  &   1 \\
      0   &  0  &  -1  &  -1  &  -1 \\
      1   &  1  &   1  &   2  &   0 \\
      0   & -1  &   1  &   0  &   1 \\
  \end{array}
  \right ]
     \end{array}$
     & $12$ \\
   & & & \\
    \hline
    &  & &   \\
     $\begin{array}{c} \includegraphics[scale=0.6]{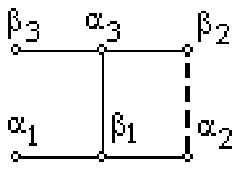} \\
        {\bf E_6(a_1)} \end{array}$ &
     $\begin{array}{c}
      \left [
   \begin{array}{cccccc}
     0 & 1 &  1 & 1 & 0  & 0 \\
     1 & 1 &  0 & 1 & -1 & 0 \\
     1 & 0 &  2 & 1 & 1 &  1 \\
    -1 & -1 & -1 & -1 & 0 & 0 \\
      0 &  1 &  -1  & 0 & -1 & 0 \\
      0 &  0 &  -1   & 0 & 0 & -1 \\
  \end{array}
  \right ]
     \end{array}$
     &
     $\begin{array}{c}
      \left [
   \begin{array}{cccccc}
     -1  &  0  &  0  & -1  &  0  &  0 \\
      0  & -1  &  0  & -1  &  1  &  0 \\
      0  &  0  & -1  & -1  & -1  & -1 \\
      1  &  1  &  1  &  2  &  0  &  1 \\
      0  & -1  &  1  &  0  &  1  &  1 \\
      0  &  0  &  1  &  1  &  1  &  0 \\
  \end{array}
  \right ]
     \end{array}$
     & $9$ \\
    &  & &   \\
     \hline
    &  & &   \\
     $\begin{array}{c} \includegraphics[scale=0.6]{E6a2.eps} \\
        {\bf E_6(a_2)} \end{array}$ &
     $\begin{array}{c}
      \left [
   \begin{array}{cccccc}
      1 & 1 &    0 &  1 & 0  & -1 \\
      1 & 1 &    0 &  1 & -1 &  0 \\
      0 & 0 &    2 &  1 & 1  &  1 \\
     -1 & -1 &  -1 & -1 & 0  &  0 \\
      0 &  1 &  -1 &  0 & -1 &  0 \\
      1 &  0 &  -1 &  0 & 0  & -1 \\
  \end{array}
  \right ]
     \end{array}$
     &
     $\begin{array}{c}
      \left [
   \begin{array}{cccccc}
     -1  &  0  &  0  & -1  &  0  &  0 \\
      0  & -1  &  0  & -1  &  1  &  0 \\
      0  &  0  & -1  & -1  & -1  & -1 \\
      1  &  1  &  1  &  2  &  0  &  1 \\
      0  & -1  &  1  &  0  &  1  &  1 \\
      0  &  0  &  1  &  1  &  1  &  0 \\
  \end{array}
  \right ]
     \end{array}$
     & $6$ \\
    &  & &   \\
     \hline
    &  & &   \\
     $\begin{array}{c} \includegraphics[scale=0.6]{D6a1.eps} \\
        {\bf D_6(a_1) = D_6(a_3)} \end{array}$  &
     $\begin{array}{c}
       \left [
   \begin{array}{cccccc}
      1  &  1  &   1  &   1  &   0  &  1 \\
      1  &  1  &   0  &   1  &  -1  &  0 \\
      1  &  0  &   1  &   1  &   1  &  0 \\
     -1  &  -1  &  -1  &  -1  &   0  & 0 \\
      0  &   1  &  -1  &   0  &  -1  & 0 \\
     -1  &   0  &   0  &   0  &   0  & -1 \\
  \end{array}
  \right ]
     \end{array}$
     &
     $\begin{array}{c}
       \left [
   \begin{array}{cccccc}
     -1  &  0  &  0  & -1  &  0  & -1 \\
      0  & -1  &  0  & -1  &  1  &  0 \\
      0  &  0  & -1  & -1  & -1  &  0 \\
      1  &  1  &  1  &  2  &  0  &  1 \\
      0  & -1  &  1  &  0  &  1  &  0 \\
      1  &  0  &  0  &  1  &  0  &  0 \\
  \end{array}
  \right ]
     \end{array}$
     &  $8$ \\
    &  & &   \\
   \hline
    &  & &   \\
  $\begin{array}{c} \includegraphics[scale=0.6]{D6a2.eps} \\
        {\bf D_6(a_2)} \end{array}$  &
     $\begin{array}{c}
       \left [
   \begin{array}{cccccc}
       0  &  1  &  1  &  0  &  1  &  0 \\
       1  &  1  &  0  &  -1 &  1  & -1 \\
       1  &  0  &  1  &  1  &  1  &  1 \\
       0  & -1  &  1  &  0  &  0  &  1 \\
      -1  & -1  & -1  &  0  & -1  &  0 \\
       0  &  1  & -1  & -1  &  0  & -1 \\
  \end{array}
  \right ]
     \end{array}$
     &
     $\begin{array}{c}
       \left [
   \begin{array}{cccccc}
      -1  &  0  &  0  &  0  & -1  &  0 \\
       0  & -1  &  0  &  0  & -1  &  1 \\
       0  &  0  & -1  &  0  & -1  & -1 \\
       0  &  0  &  0  & -1  &  0  & -1 \\
       1  &  1  &  1  &  0  &  2  &  0 \\
       0  & -1  &  1  &  1  &  0  &  2 \\
  \end{array}
  \right ]
     \end{array}$
     & $6$ \\
     &  & & \\
     \hline
   \end{tabular}
  \vspace{2mm}
  \caption{\small\hspace{3mm} The dual semi-Coxeter element $\InvDualSC$ for $l < 7$}
  \label{tab_partial Cartan_1}
\end{table}

 \begin{table}[H]
 \tiny
  \centering
  \renewcommand{\arraystretch}{1.5}
  \begin{tabular} {|c|c|c|c|}
  \hline
   The Carter  & The transpose semi-Coxeter  &  The dual semi-Coxeter                & Order \cr
   diagram     & element $\DualSC$           &  element $\InvDualSC = \DualSC^{-1}$  & of $\DualSC$  \\
  \hline  
     &  & &   \\
     $\begin{array}{c} \includegraphics[scale=0.6]{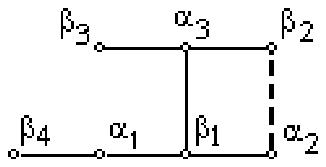} \\
        {\bf E_7(a_1)} \end{array}$  &
     $\begin{array}{c}
       \left [
   \begin{array}{ccccccc}
      1  &  1  &  1 &  1  &  0  &  0  &  1 \\
      1  &  1  &  0 &  1  & -1  &  0  &  0 \\
      1  &  0  &  2 &  1  &  1  &  1  &  0 \\
     -1  & -1  & -1 & -1  &  0  &  0  &  0 \\
      0  &  1  & -1 &  0  & -1  &  0  &  0 \\
      0  &  0  & -1 &  0  &  0  & -1  &  0 \\
     -1  &  0  &  0 &  0  &  0  &  0  & -1 \\
  \end{array}
  \right ]
     \end{array}$
     &
      $\begin{array}{c}
       \left [
   \begin{array}{ccccccc}
     -1  &  0  &  0  & -1 &  0  &  0  & -1 \\
      0  & -1  &  0  & -1 &  1  &  0  &  0 \\
      0  &  0  & -1  & -1 & -1  & -1  &  0 \\
      1  &  1  &  1  &  2 &  0  &  1  &  1 \\
      0  & -1  &  1  &  0 &  1  &  1  &  0 \\
      0  &  0  &  1  &  1 &  1  &  0  &  0 \\
      1  &  0  &  0  &  1 &  0  &  0  &  0 \\
  \end{array}
  \right ]
     \end{array}$
     & $14$ \\
     &  & &   \\
     \hline
   &  & &   \\
     $\begin{array}{c} \includegraphics[scale=0.6]{E7a2.eps} \\
        {\bf E_7(a_2)} \end{array}$  &
      $\begin{array}{c}
       \left [
   \begin{array}{ccccccc}
     0  &  1  &  1  &  1  &  0  &  0  &  0 \\
     1  &  2  &  0  &  1  & -1  &  0  &  1 \\
     1  &  0  &  2  &  1  &  1  &  1  &  0 \\
    -1  & -1  & -1  & -1  &  0  &  0  &  0 \\
     0  &  1  & -1  &  0  & -1  &  0  &  0 \\
     0  &  0  & -1  &  0  &  0  & -1  &  0 \\
     0  & -1  &  0  &  0  &  0  &  0  & -1 \\
  \end{array}
  \right ]
     \end{array}$
     &
      $\begin{array}{c}
       \left [
   \begin{array}{ccccccc}
    -1  &  0  &  0 & -1  &  0  &  0  &  0 \\
     0  & -1  &  0 & -1  &  1  &  0  & -1 \\
     0  &  0  & -1 & -1  & -1  & -1  &  0 \\
     1  &  1  &  1 &  2  &  0  &  1  &  1 \\
     0  & -1  &  1 &  0  &  1  &  1  & -1 \\
     0  &  0  &  1 &  1  &  1  &  0  &  0 \\
     0  &  1  &  0 &  1  & -1  &  0  &  0 \\
   \end{array}
  \right ]
     \end{array}$
     & $12$ \\
     &  & &   \\
     \hline
   &  & &   \\
     $\begin{array}{c} \includegraphics[scale=0.6]{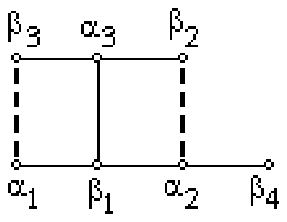} \\
        {\bf E_7(a_3)} \end{array}$  &
      $\begin{array}{c}
       \left [
   \begin{array}{ccccccc}
      1  &  1  &  0  &  1  &  0  & -1  &  0  \\
      1  &  2  &  0  &  1  & -1  &  0  &  1  \\
      0  &  0  &  2  &  1  &  1  &  1  &  0  \\
     -1  & -1  & -1  & -1  &  0  &  0  &  0  \\
      0  &  1  & -1  &  0  & -1  &  0  &  0  \\
      1  &  0  & -1  &  0  &  0  & -1  &  0  \\
      0  & -1  &  0  &  0  &  0  &  0  & -1  \\
  \end{array}
  \right ]
     \end{array}$
     &
      $\begin{array}{c}
       \left [
   \begin{array}{ccccccc}
    -1  &  0  &  0  & -1  &  0  &  1  &  0 \\
     0  & -1  &  0  & -1  &  1  &  0  & -1 \\
     0  &  0  & -1  & -1  & -1  & -1  &  0 \\
     1  &  1  &  1  &  2  &  0  &  0  &  1 \\
     0  & -1  &  1  &  0  &  1  &  1  & -1 \\
    -1  &  0  &  1  &  0  &  1  &  1  &  0 \\
     0  &  1  &  0  &  1  & -1  &  0  &  0 \\
  \end{array}
  \right ]
     \end{array}$
     & $30$ \\
     &  & &   \\
     \hline
   &  & &   \\
     $\begin{array}{c} \includegraphics[scale=0.6]{E7a4_upd_28feb2011.eps} \\
        {\bf E_7(a_4)} \end{array}$  &
      $\begin{array}{c}
       \left [
   \begin{array}{ccccccc}
      2  &  0  &  0  &  1  &  0  &  1  & -1 \\
      0  &  2  &  0  &  1  & -1  &  0  &  1 \\
      0  &  0  &  2  &  1  &  1  & -1  &  0 \\
     -1  & -1  & -1  & -1  &  0  &  0  &  0 \\
      0  &  1  & -1  &  0  & -1  &  0  &  0 \\
     -1  &  0  &  1  &  0  &  0  & -1  &  0 \\
      1  & -1  &  0  &  0  &  0  &  0  & -1 \\
  \end{array}
  \right ]
     \end{array}$
     &
      $\begin{array}{c}
       \left [
   \begin{array}{ccccccc}
     -1  &  0  &  0  & -1  &  0  & -1  & 1 \\
      0  & -1  &  0  & -1  &  1  &  0  & -1 \\
      0  &  0  & -1  & -1  & -1  &  1  &  0 \\
      1  &  1  &  1  &  2  &  0  &  0  &  0 \\
      0  & -1  &  1  &  0  &  1  & -1  & -1 \\
      1  &  0  & -1  &  0  & -1  &  1  & -1 \\
     -1  &  1  &  0  &  0  & -1  & -1  &  1 \\
  \end{array}
  \right ]
     \end{array}$
     & $6$ \\
     &  & &   \\
     \hline
     &  & &   \\
     $\begin{array}{c} \includegraphics[scale=0.6]{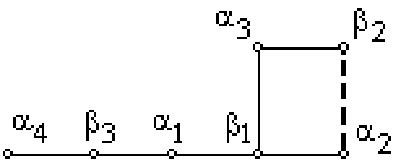} \\
        {\bf D_7(a_1) = D_7(a_4)} \end{array}$  &
     $\begin{array}{c}
       \left [
   \begin{array}{ccccccc}
     1 & 1 &    1  & 1 & 1  & 0  & 1 \\
     1 & 1 &    0  & 0 & 1  & -1 & 0 \\
     1 & 0 &    1  & 0 & 1  &  1 & 0 \\
     1 & 0 &    0  & 0 & 0  &  0 & 1 \\
    -1 & -1 &  -1  & 0 & -1 &  0 & 0 \\
     0 &  1 &  -1  & 0 &  0 & -1 & 0 \\
    -1 &  0 &   0  & -1 & 0 &  0 & -1 \\
  \end{array}
  \right ]
     \end{array}$
     &
     $\begin{array}{c}
       \left [
   \begin{array}{ccccccc}
    -1  &  0  &  0 &  0 & -1  &  0 &  -1 \\
     0  & -1  &  0 &  0 & -1  &  1 &   0 \\
     0  &  0  & -1 &  0 & -1  & -1 &   0 \\
     0  &  0  &  0 & -1 &  0  &  0 &  -1 \\
     1  &  1  &  1 &  0 &  2  &  0 &   1 \\
     0  & -1  &  1 &  0 &  0  &  1 &   0 \\
     1  &  0  &  0 &  1 &  1  &  0 &   1 \\
  \end{array}
  \right ]
     \end{array}$
     &  $20$ \\
    &  & &   \\
     \hline
    &  & &   \\
     $\begin{array}{c} \includegraphics[scale=0.6]{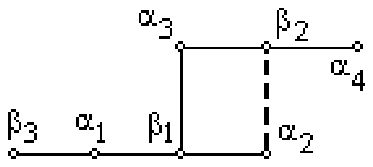} \\
        {\bf D_7(a_2) = D_7(a_3)} \end{array}$  &
     $\begin{array}{c}
       \left [
   \begin{array}{ccccccc}
     1 & 1 &    1  & 0 & 1  & 0  & 1 \\
     1 & 1 &    0  & -1 & 1  & -1 & 0 \\
     1 & 0 &    1  & 1 & 1  &  1 & 0 \\
     0 & -1 &   1  & 0 & 0  &  1 & 0 \\
    -1 & -1 &  -1  & 0 & -1 &  0 & 0 \\
     0 &  1 &  -1  & -1 &  0 & -1 & 0 \\
    -1 &  0 &   0  & 0 & 0 &  0 & -1 \\
  \end{array}
  \right ]
     \end{array}$
     &
     $\begin{array}{c}
       \left [
   \begin{array}{ccccccc}
     -1  &  0  &  0  &  0  & -1 &  0  & -1 \\
      0  & -1  &  0  &  0  & -1 &  1  &  0 \\
      0  &  0  & -1  &  0  & -1 &  -1 &  0 \\
      0  &  0  &  0  & -1  &  0 &  -1 &  0 \\
      1  &  1  &  1  &  0  &  2 &  0  &  1 \\
      0  & -1  &  1  &  1  &  0 &  2  &  0 \\
      1  &  0  &  0  &  0  &  1 &  0  &  0 \\
  \end{array}
  \right ]
     \end{array}$
     & $24$ \\
    &  & &   \\
     \hline
\end{tabular}
  \vspace{2mm}
  \caption{\small\hspace{3mm} (cont.) The dual semi-Coxeter element $\InvDualSC$ for
  $l = 7$}
  \label{tab_partial Cartan_2}
\end{table}

\newpage
~\\

 \begin{table}[H]
 \tiny
  \centering
  \renewcommand{\arraystretch}{1.5}
  \begin{tabular} {|c|c|c|c|}
  \hline
   The Carter  & The transpose semi-Coxeter  &  The dual semi-Coxeter                & Order \cr
   diagram     & element $\DualSC$           &  element $\InvDualSC = \DualSC^{-1}$  & of $\DualSC$  \\
   \hline
   &  & &   \\
     $\begin{array}{c} \includegraphics[scale=0.6]{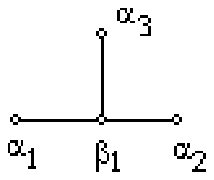} \\
        {\bf D_4} \end{array}$  &
      $\begin{array}{c}
       \left [
   \begin{array}{cccc}
     0 & 1 &   1 & 1  \\
     1 & 0 &   1 & 1 \\
     1 & 1 &   0 & 1  \\
    -1 & -1 &  -1 & -1  \\
  \end{array}
  \right ]
     \end{array}$
     &
     $\begin{array}{c}
       \left [
   \begin{array}{cccc}
     -1  &  0  &  0  & -1 \\
      0  & -1  &  0  & -1 \\
      0  &  0  & -1  & -1 \\
      1  &  1  &  1  &  2 \\
  \end{array}
  \right ]
     \end{array}$
     & $6$ \\
    &  & &   \\
   \hline
   &  & &   \\
     $\begin{array}{c} \includegraphics[scale=0.6]{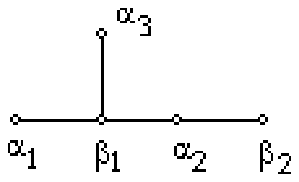} \\
        {\bf D_5} \end{array}$ &
     $\begin{array}{c}
       \left [
   \begin{array}{ccccc}
     0 & 1 &  1 & 1 & 0  \\
     1 & 1 &  1 & 1 & 1 \\
     1 & 1 &  0 & 1 & 0  \\
     -1 & -1 &  -1 & -1 & 0  \\
      0 &  -1 &   0  & 0 & -1 \\
  \end{array}
  \right ]
     \end{array}$
     &
     $\begin{array}{c}
       \left [
   \begin{array}{ccccc}
     -1  &  0  &  0  & -1  & 0  \\
      0  & -1  &  0  & -1  & -1 \\
      0  &  0  & -1  & -1  &  0 \\
      1  &  1  &  1  &  2  &  1 \\
      0  &  1  &  0  &  1  &  0 \\
  \end{array}
  \right ]
     \end{array}$
     & $8$ \\
    &  & &   \\
      \hline
    &  & & \\
     $\begin{array}{c} \includegraphics[scale=0.6]{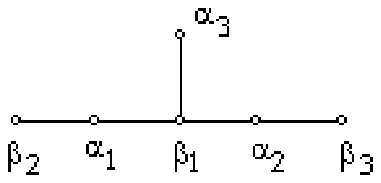} \\
        {\bf E_6} \end{array}$  &
     $\begin{array}{c}
       \left [
   \begin{array}{cccccc}
      1  &  1  &   1  &   1  &   1  &   0 \\
      1  &  1  &   1  &   1  &   0  &   1 \\
      1  &  1  &   0  &  1  &   0  &   0 \\
     -1  &  -1  &  -1 &  -1  &   0  &   0 \\
     -1  &   0  &   0 &   0  &  -1  &   0 \\
      0  &  -1  &   0 &    0 &   0  &  -1  \\
  \end{array}
  \right ]
     \end{array}$
     &
   $\begin{array}{c}
       \left [
   \begin{array}{cccccc}
      -1  &  0  &  0 &  -1  & -1  &  0 \\
       0  & -1  &  0 &  -1  &  0  & -1 \\
       0  &  0  & -1 &  -1  &  0  &  0 \\
       1  &  1  &  1 &   2  &  1  &  1 \\
       1  &  0  &  0 &   1  &  0  &  0 \\
       0  &  1  &  0 &   1  &  0  &  0  \\
  \end{array}
  \right ]
     \end{array}$
     &  $12$ \\
   &  & &   \\
     \hline
     $\begin{array}{c} \includegraphics[scale=0.6]{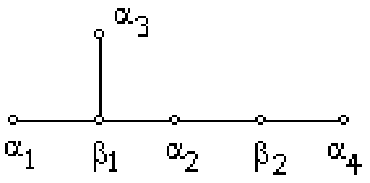} \\
        {\bf D_6} \end{array}$  &
     $\begin{array}{c}
       \left [
   \begin{array}{cccccc}
       0  &  1  &  1  &  0  &  1  &  0 \\
       1  &  1  &  1  &  1  &  1  &  1 \\
       1  &  1  &  0  &  0  &  1  &  0 \\
       0  &  1  &  0  &  0  &  0  &  1 \\
      -1  & -1  & -1  &  0  & -1  &  0 \\
       0  & -1  &  0  & -1  &  0  & -1 \\
  \end{array}
  \right ]
     \end{array}$
     &
     $\begin{array}{c}
       \left [
   \begin{array}{cccccc}
       -1  &  0  &  0  &  0  & -1  &  0 \\
        0  & -1  &  0  &  0  & -1  & -1 \\
        0  &  0  & -1  &  0  & -1  &  0 \\
        0  &  0  &  0  & -1  &  0  & -1 \\
        1  &  1  &  1  &  0  &  2  &  1 \\
        0  &  1  &  0  &  1  &  1  &  1 \\
  \end{array}
  \right ]
     \end{array}$
     & $10$ \\
     &  & &   \\
     \hline
      &  & &   \\
    $\begin{array}{c} \includegraphics[scale=0.6]{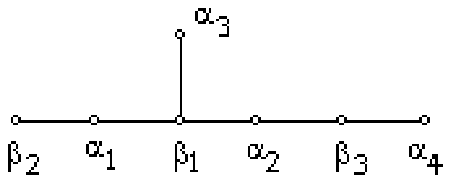} \\
        {\bf E_7} \end{array}$  &
     $\begin{array}{c}
       \left [
   \begin{array}{ccccccc}
      1  &  1  &  1  &  0  &  1  &  1  &  0 \\
      1  &  1  &  1  &  1  &  1  &  0  &  1 \\
      1  &  1  &  0  &  0  &  1  &  0  &  0 \\
      0  &  1  &  0  &  0  &  0  &  0  &  1 \\
     -1  & -1  & -1  &  0  & -1  &  0  &  0 \\
     -1  &  0  &  0  &  0  &  0  & -1  &  0 \\
      0  & -1  &  0  & -1  &  0  &  0  & -1 \\
  \end{array}
  \right ]
     \end{array}$
     &
     $\begin{array}{c}
       \left [
   \begin{array}{ccccccc}
      -1  &  0  &  0  &  0  & -1  & -1  &  0 \\
       0  & -1  &  0  &  0  & -1  &  0  & -1 \\
       0  &  0  & -1  &  0  & -1  &  0  &  0 \\
       0  &  0  &  0  & -1  &  0  &  0  & -1 \\
       1  &  1  &  1  &  0  &  2  &  1  &  1 \\
       1  &  0  &  0  &  0  &  1  &  0  &  0 \\
       0  &  1  &  0  &  1  &  1  &  0  &  1 \\
  \end{array}
  \right ]
     \end{array}$
     & $18$ \\
     &  & &   \\
     \hline
      &  & &   \\
    $\begin{array}{c} \includegraphics[scale=0.6]{D7_pure.eps} \\
        {\bf D_7} \end{array}$  &
     $\begin{array}{c}
       \left [
   \begin{array}{ccccccc}
      0  &  1  &  1  &  0  &  1  &  0  &  0 \\
      1  &  1  &  1  &  1  &  1  &  1  &  0 \\
      1  &  1  &  0  &  0  &  1  &  0  &  0 \\
      0  &  1  &  0  &  1  &  0  &  1  &  1 \\
     -1  & -1  & -1  &  0  & -1  &  0  &  0 \\
      0  & -1  &  0  & -1  &  0  & -1  &  0 \\
      0  &  0  &  0  & -1  &  0  &  0  & -1 \\
  \end{array}
  \right ]
     \end{array}$
     &
      $\begin{array}{c}
       \left [
   \begin{array}{ccccccc}
      -1  &  0  &  0  &  0  & -1  &  0  &  0 \\
       0  & -1  &  0  &  0  & -1  & -1  &  0 \\
       0  &  0  & -1  &  0  & -1  &  0  &  0 \\
       0  &  0  &  0  & -1  &  0  & -1  & -1 \\
       1  &  1  &  1  &  0  &  2  &  1  &  0 \\
       0  &  1  &  0  &  1  &  1  &  1  &  1 \\
       0  &  0  &  0  &  1  &  0  &  1  &  0 \\
  \end{array}
  \right ]
     \end{array}$
     & $12$ \\
     &  & &   \\
     \hline
 \end{tabular}
  \vspace{2mm}
  \caption{\small\hspace{3mm} (cont.) The dual semi-Coxeter element $\InvDualSC$}
  \label{tab_partial Cartan_3}
\end{table}

%% file: A1orbits.tex
\newpage
\section{\sc\bf Semi-Coxeter orbits}
   \label{sec_orbits}

Recall, that orbits of dual semi-Coxeter element acting on the linkage diagrams
are said to be {\it semi-Coxeter orbits}, see Section \ref{sec_semi_Coxeter}.

\subsection{Semi-Coxeter orbits for $D_l(a_i)$, $E_l(a_i)$, $D_l$, $E_l$, where $l < 7$}
  \label{sec_DEless7}
~\\


\begin{table}[H]
  \centering
  \scriptsize
  \renewcommand{\arraystretch}{1.7}  
  \begin{tabular} {|c|c|c|c|}
  \hline
  $\begin{array}{c} ~\\ \includegraphics[scale=0.5]{D4_pure.eps} \\ {\bf D_4} \end{array}$
      & $\begin{array}{c} \text{Orbit $1$} \\  \text{(red\footnotemark[1], $L_{12}$)} \end{array}$
      & $\begin{array}{c} \text{Orbit $2$} \\  \text{(green, $L_{13}$)} \end{array}$
      & $\begin{array}{c} \text{Orbit $3$} \\  \text{(blue, $L_{23}$)} \end{array}$ \cr
  \hline
   $\begin{array}{c}
     \gamma^{\vee} \\
     \InvDualSC\gamma^{\vee} \\
     (\InvDualSC)^2\gamma^{\vee} \\
     (\InvDualSC)^3\gamma^{\vee}\\
     (\InvDualSC)^4\gamma^{\vee}\\
     (\InvDualSC)^5\gamma^{\vee}\\
    \end{array}$
    & 
   $\begin{array}{c}
    \fbox{0, 0, -1, 0}  \\
      0, 0, 1, -1 \\
      1, 1, 0, -1  \\
    \fbox{0, 0, 1, 0}  \\
      0, 0, -1, 1 \\
      -1, -1, 0, 1 \\
    \end{array}$
    & 
   $\begin{array}{c}
   \fbox{0, -1, 0, 0}  \\
       0, 1, 0, -1  \\
       1, 0, 1, -1  \\
   \fbox{0, 1, 0, 0} \\
       0, -1, 0, 1  \\
       -1, 0, -1, 1  \\
    \end{array}$
    &
    $\begin{array}{c}
    \fbox{-1, 0, 0, 0}  \\
       1, 0, 0, -1  \\
       0, 1, 1, -1  \\
    \fbox{1,  0,  0, 0}  \\
       -1, 0, 0, 1  \\
       0, -1, -1, 1  \\
    \end{array}$ \\
   \hline
\end{tabular}
  \vspace{2mm}
\begin{tabular} {|c|c|c|c|}
  \hline   
       & $\begin{array}{c} \text{Orbit $4$} \\ \text{(blue, $L_{12}$)}  \end{array}$
       & $\begin{array}{c} \text{Orbit $5$} \\ \text{(brown, $L_{13}$)}  \end{array}$
       & $\begin{array}{c} \text{Orbit $6$} \\ \text{(green, $L_{23}$)}  \end{array}$ \cr
  \hline
   $\begin{array}{c}
     \gamma^{\vee} \\
     (\InvDualSC)\gamma^{\vee} \\
    \end{array}$
    &
   $\begin{array}{c}
    \fbox{-1, 1, 0, 0}  \\
    \fbox{1, -1, 0, 0}  \\
    \end{array}$
    &
  $\begin{array}{c}
    \fbox{1, 0, -1, 0}  \\
    \fbox{-1, 0, 1, 0}  \\
    \end{array}$
    &
  $\begin{array}{c}
    \fbox{0, 1, -1, 0}  \\
    \fbox{0, -1, 1, 0}  \\
    \end{array}$ \\
  \hline
\end{tabular}
  \caption{\hspace{3mm}\footnotesize ${\bf D_4}$, there exist $6$ semi-Coxeter orbits. All orbits are self-opposite}
  \label{tab_D4}
\end{table}

\footnotetext[1]{Here and below in all tables unicolored linkage labels vectors are framed by a rectangle.}


\begin{figure}[H]
\centering
\includegraphics[scale=1.5]{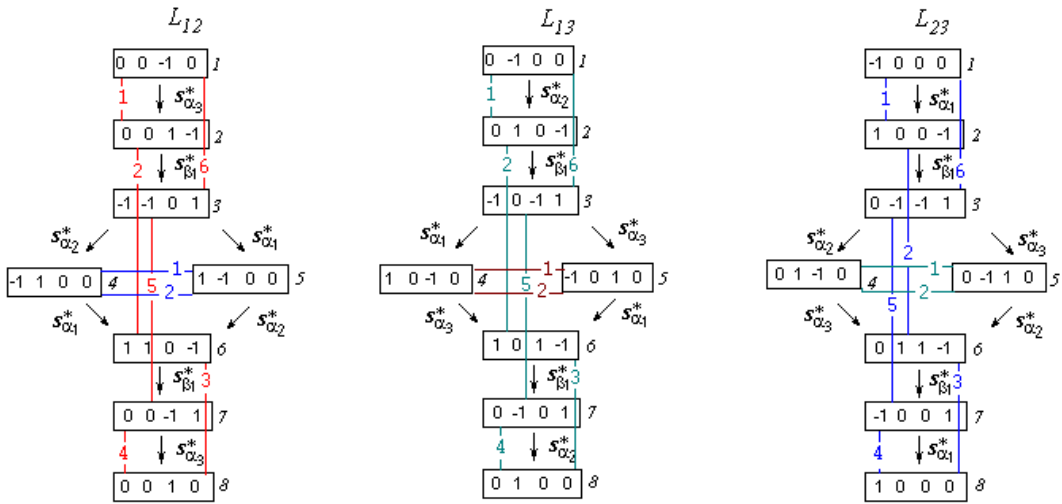}
\vspace{3mm}
\caption{\hspace{3mm}\footnotesize $D_4$, $6$ semi-Coxeter orbits, three of length $6$, three of length $2$}
\label{D4_orbits_16c}
\end{figure}


\begin{table}[H]
  \centering
  \scriptsize
  \renewcommand{\arraystretch}{1.3}  
  \begin{tabular} {|c|c|c|c|}
  \hline
  $\begin{array}{c}  ~\\ \includegraphics[scale=0.5]{D4a1.eps} \end{array}$
  $\begin{array}{c}  ~\\ {\bf D_4(a_1)} \end{array}$
      & $\begin{array}{c} \text{Orbit $1$} \\ \text{(green, II)} \\   \end{array}$
      & $\begin{array}{c} \text{Orbit $2$} \\ \text{(red, II)} \\     \end{array}$
      & $\begin{array}{c} \text{Orbit $3$} \\ \text{(brown, III)} \\  \end{array}$ \cr
  \hline
   $\begin{array}{c}
     \gamma^{\vee} \\
     \InvDualSC\gamma^{\vee} \\
     (\InvDualSC)^2\gamma^{\vee} \\
     (\InvDualSC)^3\gamma^{\vee}\\
    \end{array}$
    & 
   $\begin{array}{c}
    \fbox{1, 0, 0, 0}  \\
      1, 0, -1, -1  \\
    \fbox{-1, 0, 0, 0}  \\
      -1, 0, 1, 1 \\
    \end{array}$
    & 
   $\begin{array}{c}
       0, 1, 1, -1  \\
   \fbox{0, 1, 0, 0}  \\
       0, -1, -1, 1  \\
   \fbox{0, -1, 0, 0} \\
    \end{array}$
    &
    $\begin{array}{c}
    \fbox{0, 0, 0, 1}  \\
      -1, -1, 0, 1  \\
    \fbox{0,  0,  0, -1}  \\
       1,  1, 0, -1 \\
    \end{array}$ \\
   \hline
\end{tabular}
  \vspace{2mm}
\begin{tabular} {|c|c|c|c|c|c|c|}
  \hline   
       & $\begin{array}{c} \text{Orbit $4$} \\ \text{(green, III)}  \end{array}$
       & $\begin{array}{c} \text{Orbit $5$} \\  \text{(red, I)}     \end{array}$
       & $\begin{array}{c} \text{Orbit $6$} \\  \text{(blue, I)} \end{array}$ \cr
       &                       &  (no unicolored diagrams)      & (no unicolored diagrams)  \\
  \hline
   $\begin{array}{c}
     \gamma^{\vee} \\
     (\InvDualSC)\gamma^{\vee} \\
     (\InvDualSC)^2\gamma^{\vee} \\
     (\InvDualSC)^3\gamma^{\vee}\\
    \end{array}$
    &
   $\begin{array}{c}
       1, -1, -1, 0  \\
    \fbox{0, 0, 1, 0}  \\
       -1,  1,  1, 0 \\
    \fbox{0,  0,  -1, 0}  \\
    \end{array}$
    &
  $\begin{array}{c}
       0, 1, 1, 0 \\
       -1, 0, 0, 1 \\
       0, -1, -1, 0 \\
       1, 0, 0, -1 \\
    \end{array}$
    &
  $\begin{array}{c}
       1, 0, -1, 0 \\
       0, -1, 0, 1 \\
       -1, 0, 1, 0 \\
       0, 1, 0, -1 \\
    \end{array}$ \\
  \hline
\end{tabular}
  \caption{\hspace{3mm}\footnotesize ${\bf D_4(a_1)}$, there exist $6$ semi-Coxeter orbits. All orbits are self-opposite. Orbits $1-4$ contain unicolored linkage diagrams.}
  \label{tab_D4a1}
\end{table}

\begin{figure}[H]
\centering
\includegraphics[scale=0.95]{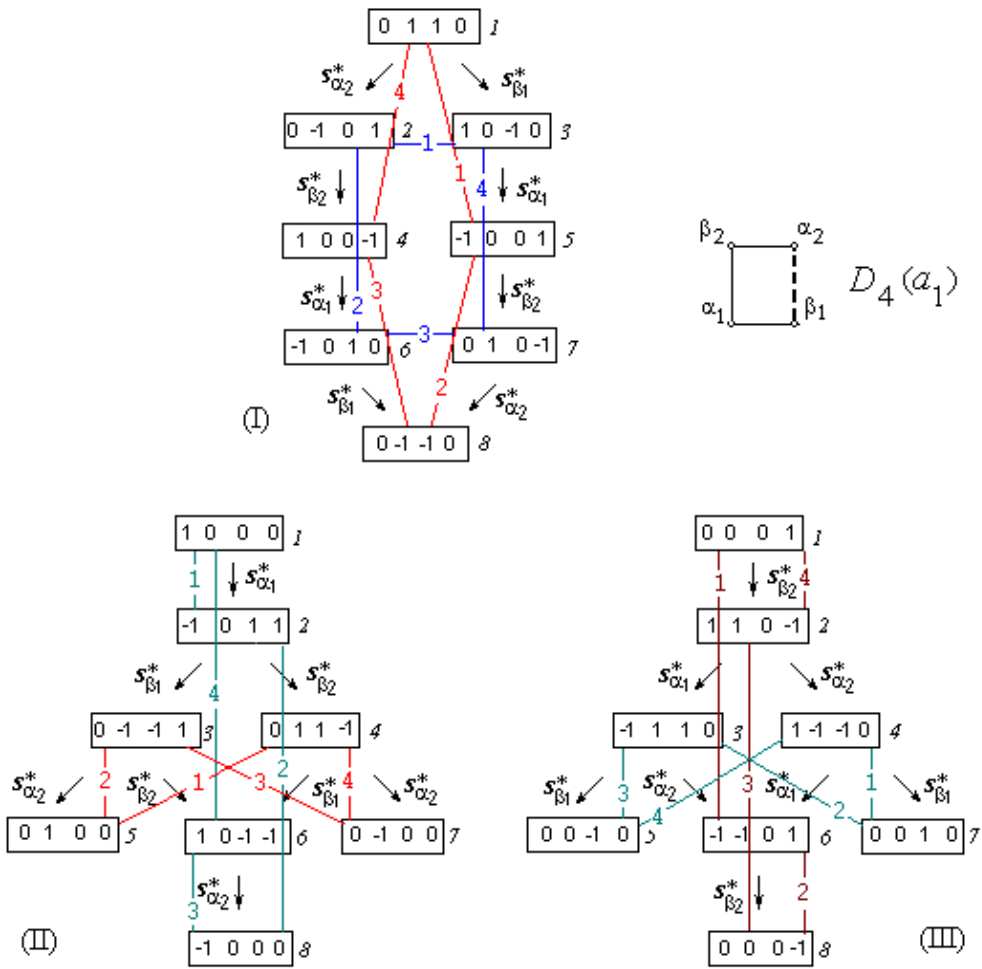}
\vspace{3mm}
\caption{\hspace{3mm}\footnotesize The linkage system of $D_4(a_1)$, three components, $24$ linkage diagrams, $3$ components}
\label{D4a1_orbits_16c}
\end{figure}


\newpage
~\\

\begin{table}[H]
  \centering
  \scriptsize
  \renewcommand{\arraystretch}{1.1}  
  \begin{tabular} {|c|c|c|c|c|}
  \hline   
        $\begin{array}{c} ~\\ \includegraphics[scale=0.5]{D5a1.eps} \\ {\bf D_5(a_1) = D_5(a_2)} \end{array}$
      & $\begin{array}{c} \text{Orbit $1$} \\  \text{(red)} \\ ~ \end{array}$
      & $\begin{array}{c} \text{Orbit $2$} \\  \text{(green)}\\ ~ \end{array}$
      & $\begin{array}{c} \text{Orbit $3$} \\  \text{(brown)} \\ \text{(self-opposite)} \end{array}$
      & $\begin{array}{c} \text{Orbit $4$} \\  \text{(blue)} \\ \text{(self-opposite)} \end{array}$ \cr
  \hline   
   $\begin{array}{c}
     \gamma^{\vee} \\
     (\InvDualSC)\gamma^{\vee} \\
     (\InvDualSC)^2\gamma^{\vee} \\
     (\InvDualSC)^3\gamma^{\vee}\\
     (\InvDualSC)^4\gamma^{\vee}\\
     (\InvDualSC)^5\gamma^{\vee}\\
     (\InvDualSC)^6\gamma^{\vee}\\
     (\InvDualSC)^7\gamma^{\vee}\\
     (\InvDualSC)^8\gamma^{\vee}\\
     (\InvDualSC)^9\gamma^{\vee}\\
     (\InvDualSC)^{10}\gamma^{\vee}\\
     (\InvDualSC)^{11}\gamma^{\vee}\\
    \end{array}$
    & 
   $\begin{array}{c}
    \fbox{0, 0, -1, 0, 0}  \\
      0, 0, 1, -1, -1  \\
      1, 0,  1, -1, 0 \\
      0, 1, 0,  0, 1  \\ 
      0, 0, -1, 1,  0 \\
     -1, -1, 0,  1, -1 \\  
    \fbox{0, -1,  0,  0, 0} \\
     0, 1, 0,  -1, 1 \\ 
     1, 1, 0,  -1, 0 \\
     0, 0,  1,  0, -1 \\ 
     0, -1, 0,  1, 0 \\
     -1, 0, -1,  1, 1 \\ 
    \end{array}$
    &
   $\begin{array}{c}
       \fbox{-1, 1, 0, 0, 0}  \\
       1, -1, 0, 0, -1  \\
       \fbox{-1, 0,  1, 0, 0} \\
       1, 0, -1,  0, 1  \\ 
      ~ ~ ~ ~ \\
      ~ ~ ~ ~ \\
      ~ ~ ~ ~ \\
      ~ ~ ~ ~ \\ 
      ~ ~ ~ ~ \\
      ~ ~ ~ ~ \\
      ~ ~ ~ ~ \\
      ~ ~ ~ ~ \\  
    \end{array}$
    &
   $\begin{array}{c}
    \fbox{-1, 0, 0, 0, 0} \hspace{0.2mm}  \\
      1, 0, 0, -1, 0   \\
      0, 1,  1, -1, 0  \\
    \fbox{1, 0, 0,  0, 0}   \\ 
      -1, 0,  0, 1,  0  \\
     0, -1, -1,  1, 0  \\  
       ~ ~ ~ ~ \\
      ~ ~ ~ ~ \\
      ~ ~ ~ ~ \\
      ~ ~ ~ ~ \\
      ~ ~ ~ ~ \\
      ~ ~ ~ ~ \\
    \end{array}$
    &
   $\begin{array}{c}
       \fbox{0, 0, 0, 0, -1}  \\
       0, -1,  1, 0, -1  \\
       \fbox{0, 0, 0, 0, 1}  \\
       0, 1,  -1, 0, 1 \\
       ~ ~ ~ ~ \\
       ~ ~ ~ ~ \\
       ~ ~ ~ ~ \\
      ~ ~ ~ ~ \\
      ~ ~ ~ ~ \\
      ~ ~ ~ ~ \\
      ~ ~ ~ ~ \\
      ~ ~ ~ ~ \\
    \end{array}$ \\
  \hline
  \end{tabular}
  \vspace{2mm}
  \caption{\hspace{3mm}   \footnotesize
  ${\bf D_5(a_1) = D_5(a_2)}$, $6$ semi-Coxeter orbits.
   All orbits contain unicolored linkage diagrams.
   Orbits $1$, $2$ have opposite orbits (in bottom component).
   Orbits $3$, $4$ are self-opposite}
  \label{tab_D5a1}
\end{table}

\begin{figure}[H]
\centering
\includegraphics[scale=1.3]{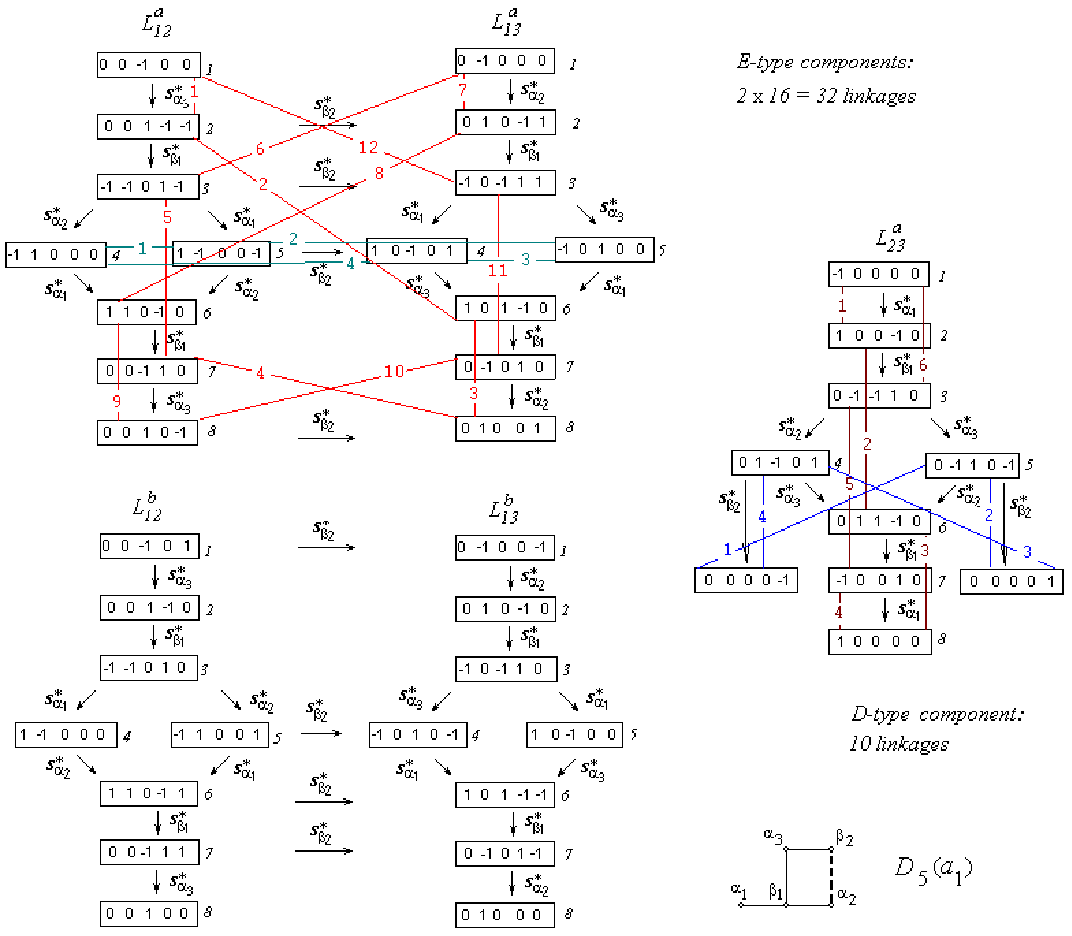}
\vspace{3mm}
\caption{\hspace{3mm} \footnotesize The linkage system of ${\bf D_5(a_1)}$, $42$ linkage diagrams, $3$ components}
\label{D5a1_orbits_16c}
\end{figure}


\newpage
~\\

\begin{table}[H]
  \centering
  \scriptsize
  \renewcommand{\arraystretch}{1.1}  
  \begin{tabular} {|c||c|c|c|c|}
  \hline   
    $\begin{array}{c} ~\\ \includegraphics[scale=0.5]{D5_pure.eps} \\ {\bf D_5}  \end{array}$
          & $\begin{array}{c} \text{Orbit $1$} \\  \text{(red, $E$-type)} \\ ~ \end{array}$
          & $\begin{array}{c} \text{Orbit $2$} \\  \text{(blue, $E$-type)} \\ ~ \end{array}$
          & $\begin{array}{c} \text{Orbit $3$} \\  \text{(red, $D$-type)} \\ \text{(self-opposite)} \end{array}$
          & $\begin{array}{c} \text{Orbit $4$} \\  \text{(blue, $D$-type)} \\ \text{(self-opposite)} \end{array}$ \\
  \hline  
   $\begin{array}{c}
     \gamma^{\vee} \\
      (\InvDualSC)\gamma^{\vee} \\
      (\InvDualSC)^2\gamma^{\vee} \\
      (\InvDualSC)^3\gamma^{\vee}\\
      (\InvDualSC)^4\gamma^{\vee}\\
      (\InvDualSC)^5\gamma^{\vee}\\
      (\InvDualSC)^6\gamma^{\vee}\\
      (\InvDualSC)^7\gamma^{\vee}\\
    \end{array}$
    & 
   $\begin{array}{c}
    \fbox{0, 0, -1, 0, 0}  \\
      0, 0, 1, -1, 0  \\
      1, 1,  0, -1, -1 \\
      0, 1, 1,  -1, 0  \\ 
    \fbox{1, 0, 0, 0, 0} \\
     -1, 0, 0,  1, 0 \\  
     0, -1, -1,  1, 1 \\
     -1, -1, 0,  1, 0 \\ 
    \end{array}$
    &
   $\begin{array}{c}
      0, 0, 1, 0, -1  \\
    \fbox{0, 1, -1, 0, 0}  \\
      0, -1, 1,  0, 1 \\
      0, 0, -1,  1, -1  \\ 
     -1, 0, 0, 0, 1 \\
    \fbox{1, -1, 0,  0, 0} \\  
     -1, 1,  0, 0, -1 \\
      1, 0, 0,  -1, 1 \\ 
    \end{array}$
    &
   $\begin{array}{c}
    \fbox{0, 0, 0, 0, -1} \hspace{1.2mm} \\
      0, 1, 0, -1, 0   \\
      1, 0,  1, -1, 0  \\
      0, 1, 0,  0, -1 \\ 
    \fbox{0, 0,  0, 0, 1}  \\
      0, -1, 0,  1, 0  \\  
      -1, 0, -1, 1, 0  \\
      0, -1,  0, 0, 1  \\ 
    \end{array}$
    &
   $\begin{array}{c}
      \fbox{1, 0, -1, 0, 0}   \\
      \fbox{-1, 0,  1, 0, 0}  \\
       ~ ~ ~ ~ \\
       ~ ~ ~ ~ \\
       ~ ~ ~ ~ \\
       ~ ~ ~ ~ \\
       ~ ~ ~ ~ \\
       ~ ~ ~ ~ \\
    \end{array}$ \\
  \hline
  \end{tabular}
  \vspace{2mm}
  \caption{\hspace{3mm}   \footnotesize
  ${\bf D_5}$, there exist $6$ semi-Coxeter orbits.
   All orbits contain unicolored linkage diagrams.
   Orbits $1$, $2$ have opposite orbits (in $2$nd $E$-type component).
   Orbits $3$, $4$ are self-opposite}
  \label{tab_D5pu}
\end{table}

\begin{figure}[H]
\centering
\includegraphics[scale=1.25]{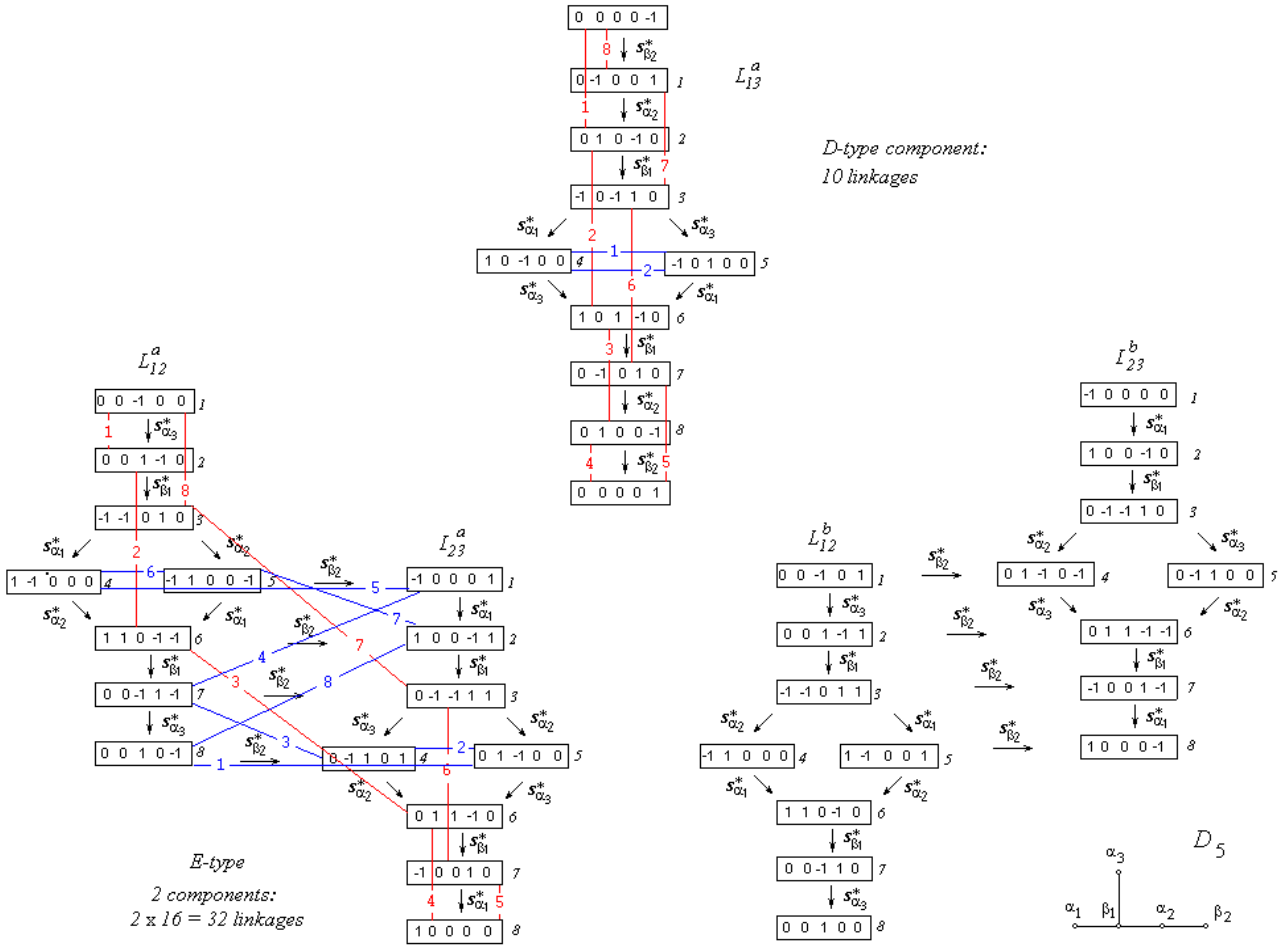}
\vspace{3mm}
\caption{\hspace{3mm} \footnotesize The linkage system of $D_5$, $42$ linkage diagrams, $3$ components}
\label{D5pu_orbits_16c}
\end{figure}


\newpage
~\\

\begin{table}[H]
  \centering
  \scriptsize
  \renewcommand{\arraystretch}{1.1}  
  \begin{tabular} {|c|c|c|c|c|c|c|}
  \hline
    $\begin{array}{c}  \includegraphics[scale=0.5]{E6a1.eps} \end{array}$
    $\begin{array}{c}  ~ \\  {\bf E_6(a_1)} \\ ~  \end{array}$
       & $\begin{array}{c} \text{Orbit $1$} \\  \text{(red)}   \end{array}$
       & $\begin{array}{c} \text{Orbit $2$} \\  \text{(green)} \end{array}$
       & $\begin{array}{c} \text{Orbit $2$} \\  \text{(blue)}  \end{array}$ \cr
  \hline
   $\begin{array}{c}
     \gamma^{\vee} \\
     \InvDualSC\gamma^{\vee} \\
     (\InvDualSC)^2\gamma^{\vee} \\
     (\InvDualSC)^3\gamma^{\vee}\\
     (\InvDualSC)^4\gamma^{\vee}\\
     (\InvDualSC)^5\gamma^{\vee}\\
     (\InvDualSC)^6\gamma^{\vee}\\
     (\InvDualSC)^7\gamma^{\vee}\\
     (\InvDualSC)^8\gamma^{\vee}\\
    \end{array}$
    & 
   $\begin{array}{c}
    \fbox{0, 0, 0, 0, 0, -1}  \\
     {0, 0, 1, -1, -1, 0}  \\
     {1, 0, 1, -1,  0, -1} \\
     {0, 1, 1, -1,  0, 0}  \\
    \fbox{1, 0, 0, 0,  0, 0} \\
     -1, 0, 0,  1, 0, 0 \\
     0, -1, -1,  1, 0, 1 \\
     -1, 0, -1,  1, 1, 0 \\
     0, 0, -1,  0, 0, 1 \\
    \end{array}$
    &
   $\begin{array}{c}
      \fbox{0, 0, 0, 0, 1, 0} \\
      {0, 1, -1, 0, 1,  1} \\
      0, 0, -1, 1, 0, 0 \\
      -1, -1, 0, 1, -1, 0 \\
      \fbox{0, -1, 0, 0, 0, 0} \\
      0, 1, 0, -1, 1, 0 \\
      1, 1, 0, -1, 0, 0 \\
      0, 0, 1, 0, -1, -1 \\
      0, -1, 1, 0, -1, 0 \\
    \end{array}$
   &
    $\begin{array}{c}
     \fbox{0, 0, 0, 0, -1, 1} \\
     0, -1, 0, 1, 0, -1 \\
     -1, 0, 0, 0, 0, 1 \\
     1, 0, -1, 0, 1, 0  \\
     \fbox{-1, 1, 0, 0, 0, 0} \\
     1, -1, 0, 0, -1, 0 \\
     -1, 0, 1, 0, 0, -1 \\
     1, 0, 0, -1, 0, 1 \\
     0, 1, 0, 0, 1, -1 \\
    \end{array}$ \\
  \hline
\end{tabular}
  \vspace{2mm}
  \caption{\hspace{3mm}   \footnotesize
  ${\bf E_6(a_1)}$, there exist $6$ semi-Coxeter orbits, each of length $9$.
   Every orbit contains the $\beta$-unicolored linkage diagram $\gamma^{\vee}$.
   Orbits $1$, $2$, $3$ have opposite orbits (starting from $-\gamma^{\vee}$)}
  \label{tab_E6a1}
\end{table}

\begin{figure}[H]
\centering
\includegraphics[scale=1.00]{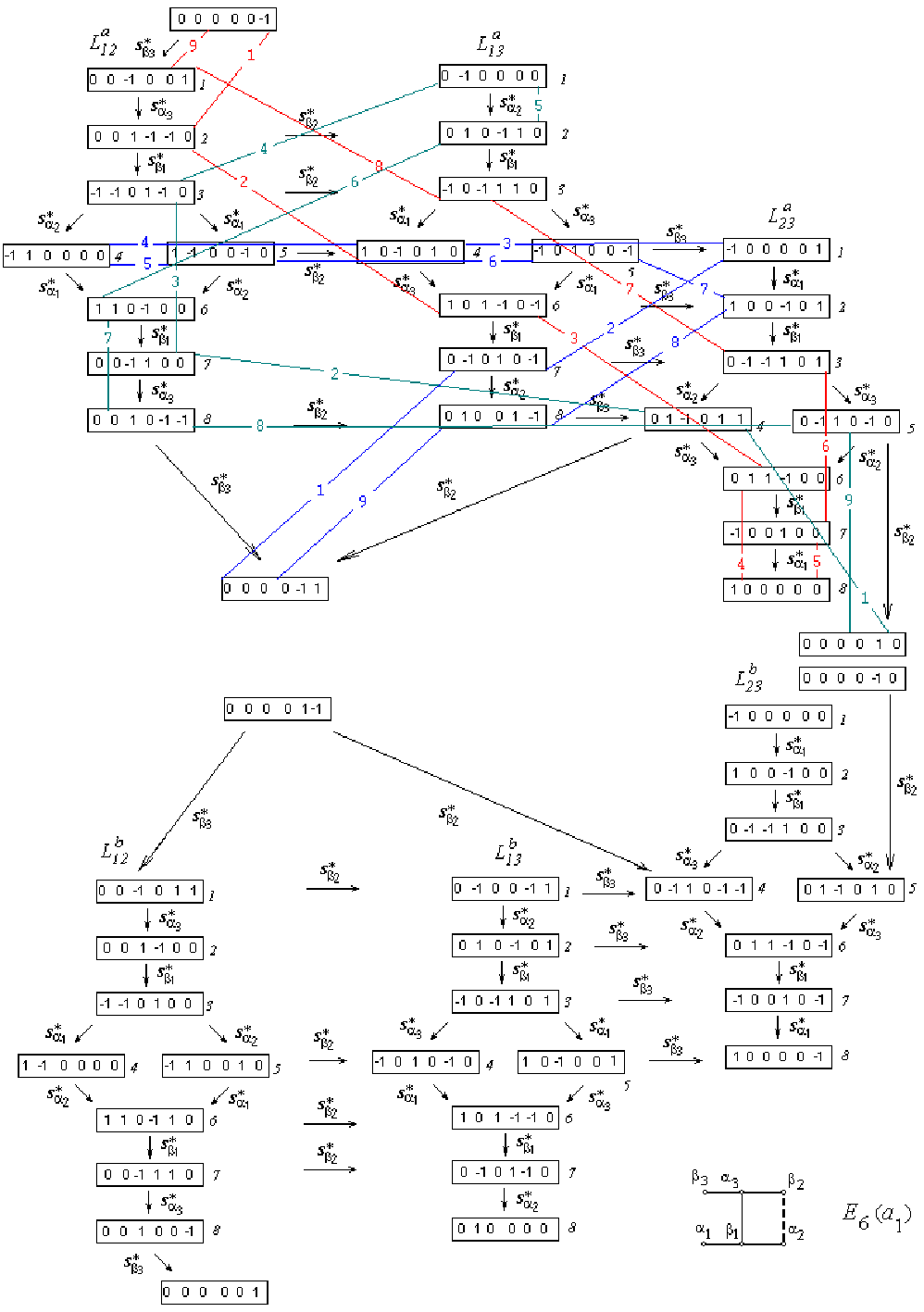}
\vspace{3mm}
\caption{\hspace{3mm} \footnotesize The linkage system of $E_6(a_1)$, two components, $54$ linkage diagrams, $6$ loctets}
\label{E6a1_orbits_16c}
\end{figure}


\begin{table}[t]
  \centering
  \scriptsize
  \renewcommand{\arraystretch}{1.1}  
  \begin{tabular} {|c|c|c|c|c|}
  \hline
    $\begin{array}{c}  \includegraphics[scale=0.5]{E6a2.eps}  \end{array}$
    $\begin{array}{c}  ~\\ {\bf E_6(a_2)} \\ ~ \end{array}$
       & $\begin{array}{c} \text{Orbit $1$ (blue)} \\  ~ \end{array}$
       & $\begin{array}{c} \text{Orbit $2$ (green)} \\  ~ \end{array}$ \cr
  \hline
   $\begin{array}{c}
     \gamma^{\vee} \\
     \InvDualSC\gamma^{\vee} \\
     (\InvDualSC)^2\gamma^{\vee} \\
     (\InvDualSC)^3\gamma^{\vee}\\
     (\InvDualSC)^4\gamma^{\vee}\\
     (\InvDualSC)^5\gamma^{\vee}\\
    \end{array}$
    & 
   $\begin{array}{c}
    \fbox{0, -1, 0, 0, 0, 0}  \\ 
      0, 1, 0,  -1, 1, 0 \\ 
      1, 1, 0, -1,  0, 1 \\ 
    \fbox{1, 0, 0, 0,  0, 0}  \\ 
     -1, 0, 0, 1,  0, -1 \\
     -1, -1, 0, 1, -1, 0  \\
    \end{array}$
    &
   $\begin{array}{c}
       \fbox{0, 0, 0, 0, -1, 0} \\
       0, -1, 1, 0, -1, -1 \\
       -1, 0, 1, 0, 0, -1 \\
       \fbox{0, 0, 0, 0, 0, 1} \\
       1, 0, -1, 0, 1, 1 \\
       0, 1, -1, 0, 1, 0 \\
    \end{array}$ \\
   \hline
\end{tabular}
  \vspace{2mm}
\begin{tabular} {|c|c|c|c|c|c|c|}
  \hline
  & $\begin{array}{c} \text{Orbit $3$ (red)}   \\  \text{(no unicolored diagrams)} \end{array}$
  & $\begin{array}{c} \text{Orbit $4$ (brown)} \\  \text{(no unicolored diagrams)} \end{array}$
  & $\begin{array}{c} \text{Orbit $5$ (turquoise)} \\ ~ \end{array}$ \cr
  \hline
   $\begin{array}{c}
     \gamma^{\vee} \\
     \InvDualSC\gamma^{\vee} \\
     (\InvDualSC)^2\gamma^{\vee} \\
     (\InvDualSC)^3\gamma^{\vee}\\
     (\InvDualSC)^4\gamma^{\vee}\\
     (\InvDualSC)^5\gamma^{\vee}\\
    \end{array}$
    &
   $\begin{array}{c}
     -1, 0, 0, 0, 0, -1  \\
      0, 0, 1, -1, -1, 0  \\
      1, 0, 1, -1,  0, 0 \\
      0, 1, 0,  0,  1, 0 \\
      0, 0, -1, 1,  0, 1  \\
      0, -1, -1, 1, 0, 0 \\
    \end{array}$
    &
  $\begin{array}{c}
       0, 0, -1, 0, 0, 1 \\
       1, 0, 0, -1, 0, 0 \\
       0, 1, 1, -1, 0, -1 \\
       0, 0, 1, 0, -1, 0 \\
       0, -1, 0, 1, 0, 0 \\
      -1, 0, -1, 1, 1, 0 \\
    \end{array}$
   &
    $\begin{array}{c}
     \fbox{0, 0, 0, 0, 1, -1} \\
     \fbox{-1, 1, 0, 0, 0, 0} \\
     1, -1, 0, 0, -1, 1 \\
     ~ ~ ~ \\
     ~ ~ ~ \\
     ~ ~ ~ \\
    \end{array}$ \\
  \hline
\end{tabular}
  \caption{\hspace{3mm} \footnotesize ${\bf E_6(a_2)}$, there exist $10$ semi-Coxeter orbits. Orbits $1-5$ have opposite orbits. Only orbits $1$, $2$, $5$ contain unicolored linkage diagrams}
  \label{E6a2_linkages}
\end{table}

\begin{figure}[H]
\centering
\includegraphics[scale=1.2]{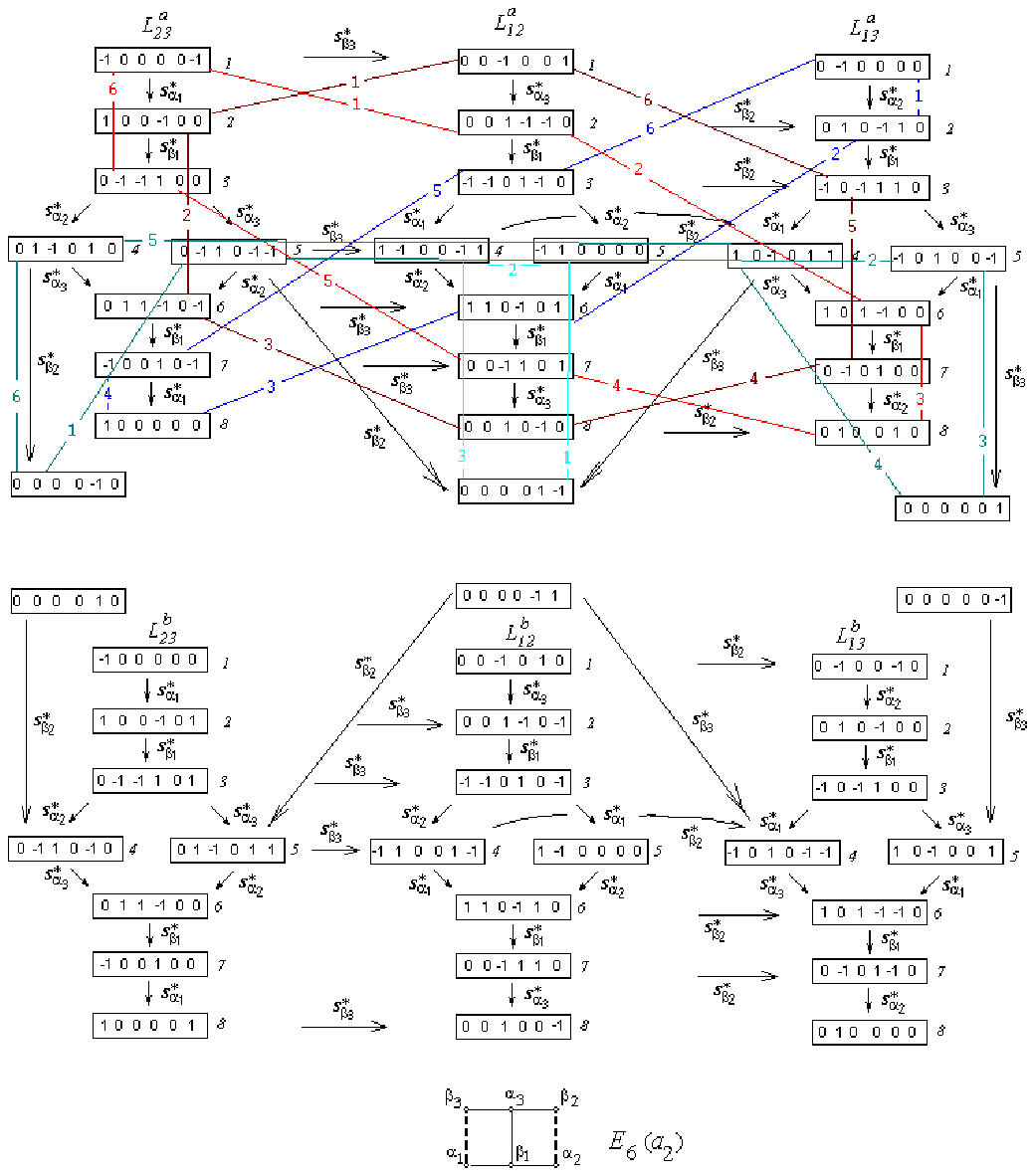}
\vspace{3mm}
\caption{\hspace{3mm} \footnotesize The linkage system of $E_6(a_2)$, two components, $54$ linkage diagrams, $6$ loctets}
\label{E6a2_orbits_16c}
\end{figure}


\newpage
~\\

\begin{table}[H]
  \centering
  \scriptsize
  \renewcommand{\arraystretch}{1.3}  
  \begin{tabular} {|c|c|c|c|}
  \hline
    $\begin{array}{c}  \includegraphics[scale=0.3]{E6_pure.eps}  \end{array}$
    $\begin{array}{c}  {\bf E_6} \end{array}$
     & $\begin{array}{c} \text{Orbit $1$}   \\  \text{(red)}  \end{array}$
     & $\begin{array}{c} \text{Orbit $2$}   \\  \text{(blue)} \end{array}$
     & $\begin{array}{c} \text{Orbit $3$}   \\  \text{(green)} \end{array}$ \cr
  \hline
   $\begin{array}{c}
     \gamma^{\vee} \\
     \InvDualSC\gamma^{\vee} \\
     (\InvDualSC)^2\gamma^{\vee} \\
     (\InvDualSC)^3\gamma^{\vee}\\
     (\InvDualSC)^4\gamma^{\vee}\\
     (\InvDualSC)^5\gamma^{\vee}\\
     (\InvDualSC)^6\gamma^{\vee}\\
     (\InvDualSC)^7\gamma^{\vee}\\
     (\InvDualSC)^8\gamma^{\vee}\\
     (\InvDualSC)^9\gamma^{\vee}\\
     (\InvDualSC)^{10}\gamma^{\vee}\\
     (\InvDualSC)^{11}\gamma^{\vee}\\
    \end{array}$
    &
   $\begin{array}{c}
    \fbox{0, 0, 0, 0, 0, -1}  \\
      0, 1, 0, -1, 0, 0  \\
      1, 0, 1, -1, -1, 0 \\
      1, 1, 0, -1,  0, -1  \\
      0, 1, 1, -1, 0, 0 \\ 
      1, 0, 0,  0, -1, 0 \\
    \fbox{0, 0, 0,  0, 1, 0} \\
     -1, 0,  0, 1, 0, 0 \\ 
     0, -1, -1,  1, 0, 1 \\
     -1, -1, 0, 1, 1, 0 \\ 
     -1, 0, -1,  1, 0, 0 \\
     0, -1, 0,  0, 0, 1 \\ 
    \end{array}$
    &
   $\begin{array}{c}
      0, 0, -1, 0, 1, 0 \\
      \fbox{-1, 0, 1, 0, 0, 0} \\ 
      1, 0, -1, 0, -1, 0 \\
      0, 0, 1, -1,  1, 0 \\ 
      0, 1, 0, 0, -1, -1 \\ 
      1, 0, 0, -1, 0, 1 \\ 
      0, 0, 1,  0, 0, -1 \\
      \fbox{0, 1, -1, 0, 0, 0} \\  
      0, -1, 1, 0, 0, 1 \\
      0, 0, -1, 1, 0, -1 \\ 
      -1, 0, 0, 0, 1,  1 \\
      0, -1, 0, 1, -1, 0 \\ 
    \end{array}$
   &
    $\begin{array}{c}
     \fbox{1, -1, 0, 0, 0, 0} \\
     -1, 1, 0, 0, 1, -1 \\
     \fbox{0, 0, 0, 0, -1, 1} \\
     ~ ~ ~\\
     ~ ~ ~\\
     ~ ~ ~\\
     ~ ~ ~\\
     ~ ~ ~\\
     ~ ~ ~\\
     ~ ~ ~\\
     ~ ~ ~\\
     ~ ~ ~\\
    \end{array}$ \\
  \hline
\end{tabular}
  \vspace{2mm}
  \caption{\hspace{3mm} \footnotesize ${\bf E_6}$, there exist $6$ semi-Coxeter orbits: four of length $12$,
  and two of length $3$.
   Orbits $1$, $2$, $3$ have opposite orbits lying in the second component}
  \label{E6pure_orbits}
\end{table}

\begin{figure}[H]
\centering
\includegraphics[scale=1.1]{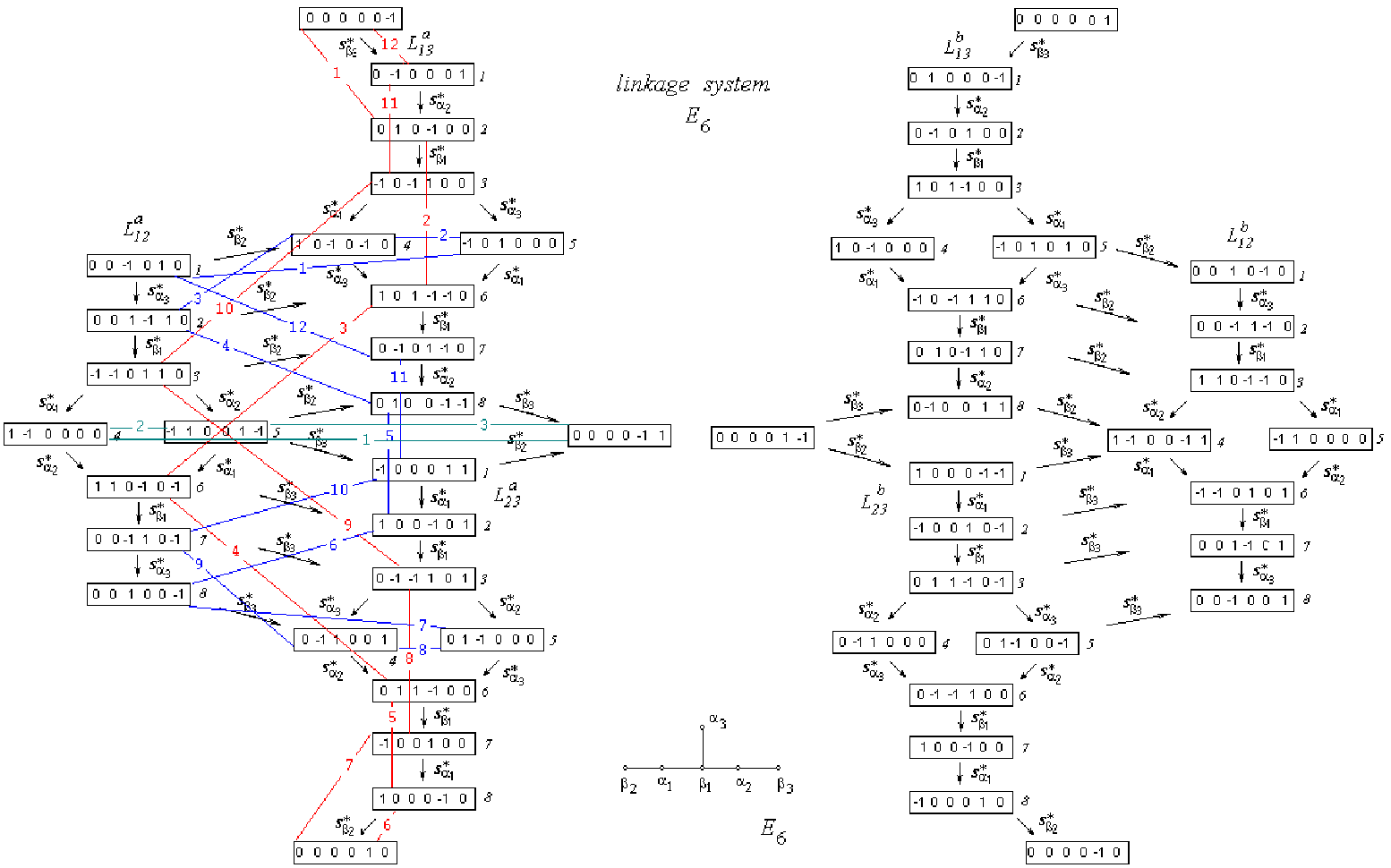}
\vspace{3mm}
\caption{\hspace{3mm} \footnotesize The linkage system of $E_6$, two components, $54$ linkage diagrams, $6$ loctets}
\label{E6pure_orbits_16c}
\end{figure}


\newpage
~\\
\begin{table}[t]
  \centering
  \scriptsize
  \renewcommand{\arraystretch}{1.8}  
  \begin{tabular} {|c|c|c|c|c|}
  \hline
    $\begin{array}{c} ~\\ \includegraphics[scale=0.5]{D6a1.eps} \\ {\bf D_6(a_1)} \end{array}$
    & $\begin{array}{c} \text{Orbit $1$}   \\  \text{(no unicolored)}  \end{array}$
    & $\begin{array}{c} \text{Orbit $2$}   \\  \text{(no unicolored)}  \end{array}$
    & $\begin{array}{c} \text{Orbit $3$}   \\  \text{(no unicolored)}  \end{array}$
    & $\begin{array}{c} \text{Orbit $4$}   \\  \text{(no unicolored)} \end{array}$ \cr
  \hline
   $\begin{array}{c}
     \gamma^{\vee} \\
     \InvDualSC\gamma^{\vee} \\
     (\InvDualSC)^2\gamma^{\vee} \\
     (\InvDualSC)^3\gamma^{\vee}\\
     (\InvDualSC)^4\gamma^{\vee}\\
     (\InvDualSC)^5\gamma^{\vee}\\
     (\InvDualSC)^6\gamma^{\vee}\\
     (\InvDualSC)^7\gamma^{\vee}\\
    \end{array}$
    & 
   $\begin{array}{c}
      0, -1, 0, 0, -1, 0  \\ 
      0, 0, 1,  -1, 0, 0  \\ 
      1, 1, 0,  -1, 1, -1 \\ 
      1, 1, 0,  -1, 0, 0  \\ 
      0, 0, 1,   0, -1, 0 \\ 
      0, -1, 0, 1,  0, 0  \\ 
     -1, 0, -1, 1,  1, 1  \\ 
     -1, 0, -1, 1, 0, 0   \\ 
   \end{array}$
    &
   $\begin{array}{c}
       0, 0, -1, 0, 1, 0 \\  
       0, 1, 0, -1, 0, 0 \\  
       1, 0, 1, -1, -1, -1 \\ 
       1, 0, 1, -1, 0, 0 \\  
       0, 1, 0, 0, 1, 0  \\  
       0, 0, -1, 1, 0, 0 \\  
      -1, -1, 0, 1, -1, 1 \\ 
      -1, -1, 0, 1, 0, 0 \\  
    \end{array}$
   &
     $\begin{array}{c}
      1, 0, -1, 0, 0, -1  \\ 
      0, 0,  1, -1, -1, 1  \\
      0, 0,  1, 0,  0, -1 \\
      1, 0, -1, 0,  1, 0  \\ 
      -1, 1, 0, 0, 0, 1 \\
      0, -1, 0,  1, -1, -1 \\  
      0, -1, 0,  0, 0, 1 \\
      -1, 1, 0,  0, 1, 0 \\  
   \end{array}$
    &
   $\begin{array}{c}
       1, -1, 0, 0, 0, -1 \\  
       0, 1, 0, -1, 1, 1 \\
       0, 1, 0, 0, 0, -1 \\
       1, -1, 0, 0, -1, 0 \\  
       -1, 0, 1, 0, 0, 1 \\
       0, 0, -1, 1, 1, -1 \\
       0, 0, -1, 0, 0, 1 \\
       -1, 0, 1, 0, -1, 0 \\  
    \end{array}$ \\
   \hline
\end{tabular}


  \begin{tabular} {|c|c|c|c|c|}
   \hline
    & $\begin{array}{c} \text{Orbit $5$}   \\   ~  \end{array}$
    & $\begin{array}{c} \text{Orbit $6$}   \\   ~  \end{array}$
    & $\begin{array}{c} \text{Orbit $7$}   \\  ~ \end{array}$
    & $\begin{array}{c} \text{Orbit $8$}   \\  ~ \end{array}$ \cr
  \hline
   $\begin{array}{c}
     \gamma^{\vee} \\
     \InvDualSC\gamma^{\vee} \\
     (\InvDualSC)^2\gamma^{\vee} \\
     (\InvDualSC)^3\gamma^{\vee}\\
     (\InvDualSC)^4\gamma^{\vee}\\
     (\InvDualSC)^5\gamma^{\vee}\\
     (\InvDualSC)^6\gamma^{\vee}\\
     (\InvDualSC)^7\gamma^{\vee}\\
    \end{array}$
    &
    $\begin{array}{c}
      \fbox{0, 0,  0, 0, 0, -1}  \\ 
      1, 0,  0, -1, 0, 0  \\
      0, 1,  1, -1,  0, 0 \\
      1, 0, 0, 0,  0, -1  \\ 
      \fbox{0, 0, 0, 0, 0, 1} \\
      -1, 0, 0,  1, 0, 0 \\  
      0, -1, -1,  1, 0, 0 \\
      -1, 0, 0,  0, 0, 1 \\  
   \end{array}$
    &
   $\begin{array}{c}
      \fbox{0, 0, -1, 0, 0, 0}  \\ 
       0, 0, 1,  -1, -1, 0 \\  
      1, 0, 1,  -1, 0, -1 \\
      1, 1, 0,  -1, 1, 0 \\  
      \fbox{0, 1, 0, 0, 0, 0} \\
      0, -1,  0, 1, -1, 0  \\
      -1, -1,  0, 1,  0, 1 \\ 
      -1, 0, -1, 1,  1, 0  \\
   \end{array}$
    & 
   $\begin{array}{c}
       \fbox{0, -1, 0, 0, 0, 0} \\  
       0, 1, 0, -1, 1, 0 \\
       1, 1, 0, -1, 0, -1 \\
       1, 0, 1, -1, -1, 0 \\  
       \fbox{0, 0, 1, 0, 0, 0} \\
       0, 0, -1, 1, 1, 0 \\
       -1, 0, -1, 1, 0, 1 \\
       -1, -1, 0, 1, -1, 0 \\  
    \end{array}$
    &
    $\begin{array}{c}
      \fbox{-1, 1,  0, 0, 0, 0}  \\ 
      1, -1,  0, 0, -1, -1  \\
      0, 0,  1, -1,  0, 1 \\
      0, 1, 0, 0,  1, -1  \\ 
      \fbox{1, 0, -1, 0, 0, 0} \\
      -1, 0, 1,  0, -1, 1 \\  
      0, -1, 0,  1, 0, -1 \\ 
      0, 0, -1,  0, 1, 1 \\  
   \end{array}$ \\
   \hline
\end{tabular}


  \begin{tabular} {|c|c|c|c|c|c|c|}
   \hline  
    & $\begin{array}{c} \text{Orbit $9$ (brown, right)}   \\  ~ \end{array}$
    & $\begin{array}{c} \text{Orbit $10$ (blue, dotted)}   \\  ~ \end{array}$ \cr
  \hline
   $\begin{array}{c}
     \gamma^{\vee} \\
     \InvDualSC\gamma^{\vee} \\
     (\InvDualSC)^2\gamma^{\vee} \\
     (\InvDualSC)^3\gamma^{\vee}\\
     (\InvDualSC)^4\gamma^{\vee}\\
     (\InvDualSC)^5\gamma^{\vee}\\
     (\InvDualSC)^6\gamma^{\vee}\\
     (\InvDualSC)^7\gamma^{\vee}\\
    \end{array}$
    &
   $\begin{array}{c}
      \fbox{-1, 0, 1, 0, 0, 0}  \\ 
       1, 0, -1,  0, 1, -1 \\  
      0, 1, 0,  -1, 0, 1 \\
      0, 0, 1,  0, -1, -1 \\  
      \fbox{1, -1, 0, 0, 0, 0} \\
      -1, 1,  0, 0, 1, 1  \\
      0, 0, -1, 1,  0, -1 \\ 
      0, -1, 0, 0, -1, 1  \\
   \end{array}$
       &
   $\begin{array}{c}
       0, 1, -1, 0, 1, 0 \\  
       \fbox{0, 0, 0, 0, -1, 0 }\\
       0, -1, 1, 0, -1, 0 \\
       \fbox{0, 0, 0, 0, 1, 0} \\  
       ~ ~ \\
       ~ ~ \\
       ~ ~ \\
       ~ ~ \\
   \end{array}$ \\
   \hline
\end{tabular}
  \vspace{2mm}
  \caption{\hspace{3mm} \footnotesize ${\bf D_6(a_1)}$, there exist $10$ semi-Coxeter orbits: nine of length $8$,
  and one of length $4$. Pairs of orbits $\{1, 2\}$, $\{3, 4\}$, $\{6, 7\}$, $\{8, 9\}$  are pairs of opposite orbits.
   Orbits $5$ and $10$ are self-opposite}
  \label{D6a1_orbits}
\end{table}



\newpage
~\\
\begin{table}[t]
  \centering
  \scriptsize
  \renewcommand{\arraystretch}{1.8}  
  \begin{tabular} {|c|c|c|c|c|c|}
  \hline
    $\begin{array}{c} ~\\ \includegraphics[scale=0.5]{D6a2.eps} \\  {\bf D_6(a_2)} \end{array}$
    & $\begin{array}{c} \text{Orbit $1$}   \\ ~ \end{array}$
    & $\begin{array}{c} \text{Orbit $2$}   \\  \text{(no unicolored)}  \end{array}$
    & $\begin{array}{c} \text{Orbit $3$}   \\ ~ \end{array}$
    & $\begin{array}{c} \text{Orbit $4$}   \\  \text{(no unicolored)} \end{array}$
    & $\begin{array}{c} \text{Orbit $5$}   \\ ~ \end{array}$ \cr
  \hline
   $\begin{array}{c}
     \gamma^{\vee} \\
     \InvDualSC\gamma^{\vee} \\
     (\InvDualSC)^2\gamma^{\vee} \\
     (\InvDualSC)^3\gamma^{\vee}\\
     (\InvDualSC)^4\gamma^{\vee}\\
     (\InvDualSC)^5\gamma^{\vee}\\
    \end{array}$
    & 
   $\begin{array}{c}
      \fbox{0, 0, -1, 0, 0, 0}  \\ 
      0, 0,  1, 0, -1, -1  \\
      1, 0,  1, 1,  -1, -1 \\
      \fbox{0, 0, 1, 0,  0, 0} \\ 
      0, 0, -1, 0,  1, 1  \\
      -1, 0, -1, -1, 1, 1 \\  
   \end{array}$
    & 
   $\begin{array}{c}
       0, 0, -1, 0, 0, 1 \\  
       0, 1, 0, -1, -1, 1 \\
       1, 1, 0, 0, -1, 0 \\
       0, 0,  1, 0, 0, -1 \\  
       0, -1, 0, 1, 1, -1 \\
       -1, -1, 0, 0, 1, 0 \\
    \end{array}$
   &  
     $\begin{array}{c}
     \fbox{0, -1, 0, -1, 0, 0}  \\ 
      0, 1,  0, 1, -1, 0  \\
      1, 0,  1, -1,  -1, 0 \\
     \fbox{0, 1, 0, 1, 0, 0}  \\ 
      0, -1, 0, -1, 1, 0 \\
      -1, 0, -1,  1, 1, 0 \\  
   \end{array}$
    & 
   $\begin{array}{c}
       0, -1, 0, 1, 0, -1 \\  
       0, 0, 1, 0, -1, 0 \\
       1, 1, 0, 0, -1, 1 \\
       0, 1, 0, -1, 0, 1 \\  
       0, 0, -1, 0, 1, 0 \\
       -1, -1, 0, 0, 1, -1 \\
    \end{array}$
    &
    $\begin{array}{c}
      1, -1,  0, 0, 0, -1  \\ 
      -1, 0,  1, 1, 0, -1  \\
     \fbox{1, -1,  0, 0,  0, 0} \\
      -1, 1, 0, 0,  0, 1  \\ 
      1, 0, -1, -1, 0, 1 \\
     \fbox{-1, 1, 0,  0, 0, 0} \\  
   \end{array}$ \\
   \hline
\end{tabular}


  \begin{tabular} {|c|c|c|c|c|c|}
   \hline
    & $\begin{array}{c} \text{Orbit $6$}   \\  ~ \end{array}$
    & $\begin{array}{c} \text{Orbit $7$}   \\  \text{(no unicolored)} \end{array}$
    & $\begin{array}{c} \text{Orbit $8$}   \\  ~ \end{array}$
    & $\begin{array}{c} \text{Orbit $9$}   \\ \text{(no unicolored)} \end{array}$
    & $\begin{array}{c} \text{Orbit $10$}  \\  ~ \end{array}$ \cr
  \hline
   $\begin{array}{c}
     \gamma^{\vee} \\
     \InvDualSC\gamma^{\vee} \\
     (\InvDualSC)^2\gamma^{\vee} \\
     (\InvDualSC)^3\gamma^{\vee}\\
     (\InvDualSC)^4\gamma^{\vee}\\
     (\InvDualSC)^5\gamma^{\vee}\\
    \end{array}$
    &  
   $\begin{array}{c}
      \fbox{0, -1, 0, 0, 0, 0}  \\ 
      0, 1,  0, 0, -1, 1 \\  
      1, 1, 0,  -1, -1, 1 \\
      \fbox{0, 1, 0,  0, 0, 0} \\  
      0, -1, 0, 0, 1, -1 \\
      -1, -1,  0, 1, 1, -1  \\
   \end{array}$
    & 
   $\begin{array}{c}
       0, -1, 0, 0, 0, -1 \\  
       0, 0, 1, 1, -1, -1 \\
       1, 0, 1, 0, -1, 0 \\
       0, 1, 0, 0, 0, 1 \\  
       0, 0, -1, -1, 1, 1 \\
       -1, 0, -1, 0, 1, 0 \\
    \end{array}$
    &  
    $\begin{array}{c}
      \fbox{0, 0,  -1, 1, 0, 0}  \\ 
      0, 0, 1, -1,  -1, 0  \\
      1, 1,  0, 1,  -1, 0 \\
      \fbox{0, 0, 1, -1,  0, 0}  \\ 
       0, 0, -1, 1, 1, 0 \\
      -1, -1, 0,  -1, 1, 0 \\  
   \end{array}$
    &  
   $\begin{array}{c}
       0, 0, -1, -1, 0, 1  \\ 
       0, 1, 0,  0, -1, 0 \\  
       1, 0,  1,  0, -1, -1 \\
       0, 0, 1,  1, 0, -1 \\  
       0, -1, 0, 0, 1, 0 \\
      -1, 0,  -1, 0, 1, 1  \\
   \end{array}$
    &  
   $\begin{array}{c}
      1, 0,  -1, 0, 0, 1  \\ 
      -1, 1,  0, -1, 0, 1  \\
     \fbox{1, 0,  -1, 0,  0, 0} \\
      -1, 0, 1, 0,  0, -1  \\ 
       1, -1, 0, 1, 0, -1 \\
     \fbox{-1, 0, 1,  0, 0, 0} \\  
    \end{array}$ \\
   \hline
\end{tabular}


  \begin{tabular} {|c|c|c|c|c|c|c|}
   \hline  
    & $\begin{array}{c} \text{Orbit $11$}  \end{array}$
    & $\begin{array}{c} \text{Orbit $12$}  \end{array}$
    & $\begin{array}{c} \text{Orbit $13$}  \end{array}$
    & $\begin{array}{c} \text{Orbit $14$}  \end{array}$ \cr
  \hline
   $\begin{array}{c}
     \gamma^{\vee} \\
     \InvDualSC\gamma^{\vee} \\
     (\InvDualSC)^2\gamma^{\vee} \\
     (\InvDualSC)^3\gamma^{\vee}\\
     (\InvDualSC)^4\gamma^{\vee}\\
     (\InvDualSC)^5\gamma^{\vee}\\
    \end{array}$
    &  
    $\begin{array}{c}
      \fbox{-1, 0,  0, 0, 0, 0}  \\ 
      1, 0,  0, 0, -1, 0  \\
      0, 1,  1, 0,  -1, 0 \\
      \fbox{1, 0, 0, 0, 0, 0} \\
      -1, 0, 0,  0, 1, 0 \\  
       0, -1, -1, 0,  1, 0  \\ 
   \end{array}$
    &  
   $\begin{array}{c}
       0, 1, -1, 0, 0, 1 \\
      \fbox{0, 0, 0, -1, 0, 0}  \\ 
       0, 0, 0,  1, 0, -1 \\  
      0, -1, 1,  0, 0, -1 \\
      \fbox{0, 0, 0, 1, 0, 0} \\
       0, 0, 0,  -1, 0, 1 \\  
   \end{array}$
    &
   $\begin{array}{c}
       \fbox{-1, 0, 1, -1, 0, 0} \\  
       \fbox{1, 0, -1, 1, 0, 0} \\
       ~ \\
       ~ \\
       ~ \\
       ~ \\
    \end{array}$
    &
   $\begin{array}{c}
       \fbox{-1, 1, 0, 1, 0, 0}\\  
       \fbox{1, -1, 0, -1, 0, 0 }\\
       ~ \\
       ~ \\
       ~ \\
       ~ \\
   \end{array}$ \\
   \hline
\end{tabular}
  \vspace{2mm}
  \caption{\hspace{3mm} \footnotesize ${\bf D_6(a_2)}$, there exist $14$ semi-Coxeter orbits: orbits $1-12$
  are of length $6$, and orbits $13,14$ are of length $2$. All orbits are self-opposite}
  \label{D6a2_orbits}
\end{table}



\newpage
~\\
\begin{table}[t]
  \centering
  \scriptsize
  \renewcommand{\arraystretch}{1.9}  
  \begin{tabular} {|c|c|c|c|c|c|c|}
  \hline
    $\begin{array}{c} ~\\ \includegraphics[scale=0.5]{D6_pure.eps} \\ {\bf D_6} \end{array}$
     & $\begin{array}{c} \text{Orbit $1$}   \\  ~ \end{array}$
     & $\begin{array}{c} \text{Orbit $2$}   \\  ~ \end{array}$
     & $\begin{array}{c} \text{Orbit $3$}   \\  ~ \end{array}$
     & $\begin{array}{c} \text{Orbit $4$}   \\  ~ \end{array}$ \cr
  \hline
   $\begin{array}{c}
     \gamma^{\vee} \\
     \InvDualSC\gamma^{\vee} \\
     (\InvDualSC)^2\gamma^{\vee} \\
     (\InvDualSC)^3\gamma^{\vee}\\
     (\InvDualSC)^4\gamma^{\vee}\\
     (\InvDualSC)^5\gamma^{\vee}\\
     (\InvDualSC)^6\gamma^{\vee} \\
     (\InvDualSC)^7\gamma^{\vee}\\
     (\InvDualSC)^8\gamma^{\vee}\\
     (\InvDualSC)^9\gamma^{\vee}\\
    \end{array}$
    & 
   $\begin{array}{c}
      \fbox{0, 0, -1, 0, 0, 0}  \\ 
      0, 0,  1, 0, -1, 0  \\  
      1, 1,  0, 0,  -1, -1 \\  
      0, 1, 1, 1,  -1, -1 \\ 
      1, 1, 0, 0,  -1, 0  \\ 
      \fbox{0, 0, 1, 0, 0, 0} \\  
      0, 0, -1, 0, 1, 0 \\ 
      -1, -1,  0, 0,  1, 1 \\ 
      0, -1, -1, -1,  1, 1  \\ 
      -1, -1, 0, 0, 1, 0 \\  
   \end{array}$
    & 
   $\begin{array}{c}
       \fbox{-1, 0, 0, 1, 0, 0} \\  
       1, 0, 0, -1, -1, 1 \\ 
       0, 0, 1, 0, 0, -1 \\ 
       0, 1,  -1, 1, 0, -1 \\ 
       0, 0, 1, 0, -1, 1 \\ 
       \fbox{1, 0, 0, -1, 0, 0} \\ 
       -1, 0, 0, 1, 1, -1 \\  
       0, 0,  -1, 0,  0,  1 \\ 
       0, -1,  1, -1,  0, 1 \\ 
       0, 0, -1, 0,  1, -1  \\ 
    \end{array}$
   &  
     $\begin{array}{c}
     \fbox{1, -1, 0, 0, 0, 0}  \\ 
      -1, 1, 0, 0, 0, -1 \\ 
      1, 0,  0, 1, -1, 0 \\ 
      0, 1, 1, -1, -1, 0  \\ 
      1,  0, 0, 1, 0, -1 \\ 
     \fbox{-1, 1, 0,  0, 0, 0} \\  
      1, -1,  0, 0, 0, 1  \\ 
      -1, 0,  0, -1,  1, 0 \\  
      0, -1, -1, 1, 1, 0  \\ 
      -1, 0, 0, -1, 0, 1 \\  
   \end{array}$
    & 
   $\begin{array}{c}
    \fbox{-1, 0, 0, 0, 0, 0} \\  
       1, 0, 0, 0, -1, 0 \\ 
       0, 1, 1, 0, -1, -1 \\ 
       1, 1, 0, 1, -1, -1 \\  
       0, 1, 1, 0, -1, 0 \\ 
    \fbox{1, 0, 0, 0, 0, 0} \\ 
       -1, 0, 0, 0, 1, 0 \\ 
       0, -1, -1, 0, 1, 1 \\  
       -1, -1, 0, -1, 1, 1 \\ 
       0, -1, -1, 0, 1, 0 \\ 
    \end{array}$ \\
   \hline
\end{tabular}


  \begin{tabular} {|c|c|c|c|c|c|c|}
   \hline  
     & $\begin{array}{c} \text{Orbit $5$} \end{array}$
     & $\begin{array}{c} \text{Orbit $6$} \end{array}$
     & $\begin{array}{c} \text{Orbit $7$} \end{array}$
     & $\begin{array}{c} \text{Orbit $8$} \end{array}$ \cr
  \hline  
   $\begin{array}{c}
     \gamma^{\vee} \\
     \InvDualSC\gamma^{\vee} \\
     (\InvDualSC)^2\gamma^{\vee} \\
     (\InvDualSC)^3\gamma^{\vee}\\
     (\InvDualSC)^4\gamma^{\vee}\\
     (\InvDualSC)^5\gamma^{\vee}\\
     (\InvDualSC)^6\gamma^{\vee} \\
     (\InvDualSC)^7\gamma^{\vee}\\
     (\InvDualSC)^8\gamma^{\vee}\\
     (\InvDualSC)^9\gamma^{\vee}\\
    \end{array}$
    & 
    $\begin{array}{c}
    \fbox{0, 0, -1, 1, 0, 0} \\  
       0, 0,  1, -1, -1,  1  \\  
       1, 0,  0, 0,  0, -1 \\  
      -1, 1, 0, 1,  0, -1  \\  
       1, 0, 0, 0, -1, 1 \\  
     \fbox{0, 0, 1, -1, 0, 0} \\  
       0, 0, -1, 1, 1, -1 \\  
      -1, 0,  0, 0, 0, 1 \\  
       1, -1, 0, -1, 0, 1 \\  
      -1, 0,  0, 0, 1, -1 \\ 
   \end{array}$
    & 
   $\begin{array}{c}
      \fbox{0, -1, 1, 0, 0, 0}  \\ 
      0, 1,  -1, 0, 0, -1 \\  
      0, 0, 1,  1, -1, 0 \\  
      1, 1, 0,  -1, -1, 0 \\  
      0, 0, 1, 1, 0, -1 \\  
    \fbox{0, 1, -1, 0, 0, 0}  \\ 
      0, -1, 1, 0, 0, 1 \\  
      0, 0, -1,  -1, 1, 0 \\ 
      -1, -1,  0, 1, 1, 0 \\ 
      0, 0, -1,  -1, 0, 1 \\ 
   \end{array}$
    & 
   $\begin{array}{c}
     \fbox{0, 0, 0, -1, 0, 0} \\  
       0, 0, 0, 1, 0, -1 \\  
       0, 1, 0, 0, -1, 0 \\  
       1, 0, 1, 0, -1, 0 \\  
       0, 1, 0, 0, 0, -1 \\  
     \fbox{0, 0, 0, 1, 0, 0} \\ 
       0, 0, 0, -1, 0, 1 \\ 
       0, -1, 0, 0, 1, 0 \\ 
       -1, 0, -1, 0, 1, 0 \\  
       0, -1, 0, 0, 0,  1 \\ 
    \end{array}$
    & 
    $\begin{array}{c}
     \fbox{1, 0,  -1, 0, 0, 0}  \\  
     \fbox{-1, 0,  1, 0, 0, 0}  \\  
         ~\\
         ~\\
         ~\\
         ~\\
         ~\\
         ~\\
         ~\\
         ~\\
   \end{array}$ \\
   \hline
\end{tabular}


  \begin{tabular} {|c|c|c|}
   \hline  
     & $\begin{array}{c} \text{Orbit $9$}   \end{array}$
     & $\begin{array}{c} \text{Orbit $10$}  \end{array}$ \cr
  \hline  
   $\begin{array}{c}
     \gamma^{\vee} \\
     \InvDualSC\gamma^{\vee} \\
    \end{array}$
    & 
    $\begin{array}{c}
    \fbox{0, -1, 1, 1, 0, 0} \\  
    \fbox{0, 1, -1, -1, 0,  0}  \\  
   \end{array}$
    & 
   $\begin{array}{c}
      \fbox{1, -1, 0, 1, 0, 0}  \\ 
      \fbox{-1, 1,  0, -1, 0, 0} \\  
   \end{array}$ \\
   \hline
\end{tabular}
  \vspace{2mm}
  \caption{\hspace{3mm} \footnotesize ${\bf D_6}$, there exist $10$ semi-Coxeter orbits:
  $7$ of length $10$, and three orbit of length $2$, all orbits contain unicolored linkage diagrams.
  All orbits are self-opposite}
  \label{D6pure_orbits}
\end{table}
~\\
~\\
~\\
~\\


\subsection{Semi-Coxeter orbits for $D_7(a_i)$, $E_7(a_i)$, $D_7$, $E_7$}
~\\


\begin{table}[H]
  \centering
  \scriptsize
  \renewcommand{\arraystretch}{1.5}  
  \begin{tabular} {|c|c|c|c|c|}
  \hline
    $\begin{array}{c} ~ \\ \includegraphics[scale=0.5]{E7a1.eps} \\ {\bf E_7(a_1)} \end{array}$
     & $\begin{array}{c} \text{Orbit $1$} \\  ~ \end{array}$     
     & $\begin{array}{c} \text{Orbit $2$} \\  ~ \end{array}$     
     & $\begin{array}{c} \text{Orbit $3$} \\  ~ \end{array}$     
     & $\begin{array}{c} \text{Orbit $4$} \\  ~ \end{array}$ \cr 
  \hline
   $\begin{array}{c}
     \gamma^{\vee} \\
     \InvDualSC\gamma^{\vee} \\
     (\InvDualSC)^2\gamma^{\vee} \\
     (\InvDualSC)^3\gamma^{\vee}\\
     (\InvDualSC)^4\gamma^{\vee}\\
     (\InvDualSC)^5\gamma^{\vee}\\
     (\InvDualSC)^6\gamma^{\vee}\\
     (\InvDualSC)^7\gamma^{\vee}\\
     (\InvDualSC)^8\gamma^{\vee}\\
     (\InvDualSC)^9\gamma^{\vee}\\
     (\InvDualSC)^{10}\gamma^{\vee}\\
     (\InvDualSC)^{11}\gamma^{\vee}\\
     (\InvDualSC)^{12}\gamma^{\vee}\\
     (\InvDualSC)^{13}\gamma^{\vee}\\
    \end{array}$
    & 
   $\begin{array}{c}
      \fbox{0, 0, 0, 0, 1, -1, 0}  \\ 
      0, 1,  0, -1, 0, 1, 0  \\
      1, 0,  0, 0, 0, -1, -1 \\
      0, 0, 1, -1,  -1, 0, 1 \\
      0, 0, 1, 0, 0, -1, -1  \\
      1, 0, 0, -1, 0, 1, 0 \\  
      0, 1,  0, 0, 1, -1, 0  \\
      \fbox{0, 0,  0, 0, -1,  1, 0} \\
      0, -1,  0, 1, 0, -1, 0  \\
      -1, 0,  0, 0,  0, 1, 1 \\
      0, 0,  -1, 1, 1, 0, -1  \\  
      0, 0,  -1, 0,  0, 1, 1 \\
      -1, 0,  0, 1, 0, -1, 0  \\
      0, -1,  0, 0,  -1, 1, 0 \\
   \end{array}$
    & 
   $\begin{array}{c}
      \fbox{0, 0, 0, 0, 0, 0, -1} \\
      1, 0, 0, -1, 0, 0, 0 \\  
      0, 1, 1, -1, 0, -1, 0 \\
      1, 0,  1, -1, -1, 0, -1 \\
      1, 0,  1, -1, 0, -1, 0 \\ 
      0, 1, 1, -1, 0, 0, 0 \\
      1, 0,  0, 0, 0, 0, -1 \\ 
     \fbox{0, 0, 0, 0, 0, 0, 1} \\ 
      -1, 0,  0, 1,  0, 0, 0 \\
      0, -1, -1, 1, 0, 1, 0  \\ 
      -1, 0,  -1, 1,  1, 0, 1 \\
      -1, 0,  -1, 1, 0, 1, 0  \\ 
      0, -1,  -1, 1,  0, 0, 0 \\
      -1, 0,  0, 0,  0, 0, 1 \\  
   \end{array}$
    & 
   $\begin{array}{c}
   \fbox{0, 0, 0, 0, -1, 0, 0} \\  
       0, -1, 1, 0, -1, -1, 0 \\
       0, 0, 1, -1, 0, 0, 0 \\
       1, 1, 0, -1, 1, 0, -1 \\ 
       1, 1, 0, -1, 0, 0, 0 \\
       0, 0, 1, 0, -1, -1, 0 \\ 
       0, -1, 1, 0, -1, 0, 0 \\
       \fbox{0, 0, 0, 0, 1, 0, 0} \\ 
       0, 1, -1, 0, 1, 1, 0 \\
       0, 0, -1, 1, 0, 0, 0 \\  
       -1, -1, 0, 1, -1, 0, 1 \\
       -1, -1, 0, 1, 0, 0, 0 \\  
       0, 0, -1,  0, 1, 1, 0 \\
       0, 1, -1, 0, 1, 0, 0 \\  
    \end{array}$
    &  
    $\begin{array}{c}
      \fbox{0, 0,  0, 0, 0, -1, 1}  \\
      -1, 0,  1, 0, -1, 0, 0  \\
      -1, 1,  0, 0,  0, 0, -1 \\ 
       0, 1, 0, -1,  1, 0, 1  \\
       0, 1, 0, 0, 0, 0, -1 \\
       1, -1, 0,  0, -1, 0, 0 \\  
      -1, 0,  1, 0, 0, -1, 1  \\
      \fbox{0, 0,  0, 0, 0, 1, -1}  \\ 
      1, 0,  -1, 0, 1, 0, 0  \\
      -1, 1,  0, 0, 0, 0, 1  \\ 
      0, -1,  0, 1, -1, 0, -1  \\
      0, -1,  0, 0, 0, 0, 1  \\ 
      -1, 1,  0, 0, 1, 0, 0  \\
      1, 0,  -1, 0, 0, 1, -1  \\ 
   \end{array}$ \\
   \hline
\end{tabular}
  \vspace{2mm}
  \caption{\hspace{3mm} \footnotesize ${\bf E_7(a_1)}$, there exist $4$ semi-Coxeter orbits,
  each of which of length $14$, all orbits are self-opposite}
  \label{E7a1_orbits}
\end{table}



\begin{table} [H]
   \centering
    \scriptsize
  \renewcommand{\arraystretch}{1.7}  
  \begin{tabular} {|c|c|c|c|c|}
  \hline
      $\begin{array}{c} ~\\ \includegraphics[scale=0.5]{E7a2.eps} \\ {\bf E_7(a_2)} \end{array}$
    & $\begin{array}{c} \text{Orbit $1$} \\  ~ \\ ~ \end{array}$     
    & $\begin{array}{c} \text{Orbit $2$} \\  ~ \\ ~ \end{array}$     
    & $\begin{array}{c} \text{Orbit $3$} \\  \text{(orbit $5$ is opposite)} \\ ~ \end{array}$     
    & $\begin{array}{c} \text{Orbit $4$} \\  \text{(orbit $6$ is opposite)} \\
                        \text{(no unicolored)} \end{array}$  \cr   
  \hline
   $\begin{array}{c}
     \gamma^{\vee} \\
      \InvDualSC\gamma^{\vee} \\
     (\InvDualSC)^2\gamma^{\vee} \\
     (\InvDualSC)^3\gamma^{\vee}\\
     (\InvDualSC)^4\gamma^{\vee}\\
     (\InvDualSC)^5\gamma^{\vee}\\
     (\InvDualSC)^6\gamma^{\vee}\\
     (\InvDualSC)^7\gamma^{\vee}\\
     (\InvDualSC)^8\gamma^{\vee}\\
     (\InvDualSC)^9\gamma^{\vee}\\
     (\InvDualSC)^{10}\gamma^{\vee}\\
     (\InvDualSC)^{11}\gamma^{\vee}\\
    \end{array}$
    & 
   $\begin{array}{c}
    0, -1, 1, 0, -1, 0, 0 \\
   \fbox{0, 0, 0, 0, 1, 0, 0} \\
    0, 1, -1, 0, 1, 1, -1 \\  
    0, 1, -1, 0, 1, 0, 0 \\
   \fbox{0, 0, 0, 0, -1, 0, 0} \\
    0, -1, 1, 0, -1, -1, 1 \\ 
      ~ \\
      ~ \\
      ~ \\
      ~ \\
      ~ \\
      ~ \\
    \end{array}$
   & 
   $\begin{array}{c}
      \fbox{0, 0, 0, 0, 1, -1, 1} \\
      \fbox{0, 0, 0, 0, -1, 1, -1} \\
      ~ \\
      ~ \\
      ~ \\
      ~ \\
      ~ \\
      ~ \\
      ~ \\
      ~ \\
      ~ \\
      ~ \\
    \end{array}$
      & 
    $\begin{array}{c}
    \fbox{0, 0, 0, 0, 0, -1, 0} \\ 
     0, 0, 1, -1, -1, 0, 0 \\  
     1, 0, 1, -1, 0, -1, 0 \\  
     0, 1, 1, -1, 0, 0, -1 \\  
     1, 1, 0, -1, 1, 0, 0  \\  
     0, 1, 0,  0, 0, 0, -1 \\  
     \fbox{0, 0, 0, 0, 0, 0, 1} \\ 
     0, -1, 0, 1, -1, 0, 0  \\ 
     -1, -1, 0, 1, 0, 0, 1 \\  
     0, -1, -1, 1, 0, 1, 0  \\ 
     -1, 0, -1, 1, 1, 0, 0  \\ 
     0, 0, -1, 0, 0, 1, 0  \\ 
    \end{array}$
   & 
     $\begin{array}{c}
      -1, 0, 0, 0, 0, 1, 0 \\ 
      1, 0, -1, 0, 1, 0, 0 \\ 
      -1, 1, 0, 0, 0, 0, -1 \\ 
      1, 0, 0, -1, 0, 0, 1 \\ 
      0, 0, 1, 0, -1, -1, -1 \\ 
      0, 0, 1, -1, 0, 0, 1 \\ 
      1, 0, 0, 0, 0, 0, -1 \\ 
      -1, 1, 0, 0, 1, 0, 0 \\ 
      1, 0, -1, 0, 0, 1, 0 \\ 
     -1, 0, 0, 1, 0, -1, 0 \\ 
      0, -1, 0, 0, -1, 1, 1 \\ 
      0, -1, 0, 1, 0, -1, 0 \\ 
    \end{array}$ \\
  \hline
\end{tabular}
\vspace{2mm}
  \caption{\hspace{3mm} ${\bf E_7(a_2)}$,
  \footnotesize $6$ semi-Coxeter orbits. Orbits $1$ and $2$ are self-opposite. Orbits $3$ and $5$ (resp. $4$ and $6$) are opposite}
  \label{E7a2_linkages}
  \end{table}



\begin{table} 
  \centering
  \scriptsize
  \renewcommand{\arraystretch}{1.9}
  \begin{tabular} {|c|c|c|c|c|c|} 
  \hline
    $\begin{array}{c} ~\\ \includegraphics[scale=0.5]{E7a3.eps} \\ {\bf E_7(a_3)} \end{array}$
       & $\begin{array}{c} \text{Orbit $1$} \\  ~ \end{array}$     
       &
       & $\begin{array}{c} \text{Orbit $1$ (cont.)} \\  ~ \end{array}$   \cr  
  \hline
   $\begin{array}{c}
     \gamma^{\vee} \\
     \InvDualSC\gamma^{\vee} \\
     (\InvDualSC)^2\gamma^{\vee} \\
     (\InvDualSC)^3\gamma^{\vee}\\
     (\InvDualSC)^4\gamma^{\vee}\\
     (\InvDualSC)^5\gamma^{\vee}\\
     (\InvDualSC)^6\gamma^{\vee}\\
     (\InvDualSC)^7\gamma^{\vee}\\
     (\InvDualSC)^8\gamma^{\vee}\\
     (\InvDualSC)^9\gamma^{\vee}\\
     (\InvDualSC)^{10}\gamma^{\vee}\\
     (\InvDualSC)^{11}\gamma^{\vee}\\
     (\InvDualSC)^{12}\gamma^{\vee} \\
     (\InvDualSC)^{13}\gamma^{\vee}\\
     (\InvDualSC)^{14}\gamma^{\vee}\\
    \end{array}$
    & 
   $\begin{array}{c}
   \fbox{0, 0, 0, 0, 0, 1, 0} \\  
     1, 0, -1, 0, 1, 1, 0 \\  
     0, 1, -1, 0, 1, 0, -1 \\ 
     0, 1, 0, -1, 0, 0, 0 \\ 
     1, 0, 1, -1, -1, 0, 0 \\ 
     0, 0, 1, 0, 0, -1, 0 \\  
     -1, 0, 0, 1, 0, 0, 0 \\  
     0, -1, -1, 1, 0, 1, 1 \\ 
     0, -1, -1, 1, 0, 0, 0 \\ 
     -1, 0, 0, 0, 0, -1, 0 \\ 
     0, 0, 1, -1, -1, 0, 0 \\ 
     1, 0, 1, -1, 0, 0, 0 \\  
     0, 1, 0, 0, 1, 0, -1 \\  
     0, 1, -1, 0, 1, 1, 0 \\  
     1, 0, -1, 0, 0, 1, 0 \\ 
    \end{array}$
    &
   $\begin{array}{c}
     (\InvDualSC)^{15}\gamma^{\vee}\\
     (\InvDualSC)^{16}\gamma^{\vee}\\
     (\InvDualSC)^{17}\gamma^{\vee}\\
     (\InvDualSC)^{18}\gamma^{\vee}\\
     (\InvDualSC)^{19}\gamma^{\vee}\\
     (\InvDualSC)^{20}\gamma^{\vee}\\
     (\InvDualSC)^{21}\gamma^{\vee}\\
     (\InvDualSC)^{22}\gamma^{\vee} \\
     (\InvDualSC)^{23}\gamma^{\vee}\\
     (\InvDualSC)^{24}\gamma^{\vee}\\
     (\InvDualSC)^{25}\gamma^{\vee}\\
     (\InvDualSC)^{26}\gamma^{\vee}\\
     (\InvDualSC)^{27}\gamma^{\vee}\\
     (\InvDualSC)^{28}\gamma^{\vee}\\
     (\InvDualSC)^{29}\gamma^{\vee}\\
    \end{array}$
    &  
   $\begin{array}{c}
   \fbox{0, 0, 0, 0, 0, -1, 0} \\  
     -1, 0, 1, 0, -1, -1, 0 \\ 
     0, -1, 1, 0, -1, 0, 1 \\  
     0, -1, 0, 1, 0, 0, 0 \\  
     -1, 0, -1, 1, 1, 0, 0 \\ 
      0, 0, -1, 0, 0, 1, 0 \\ 
     1, 0, 0, -1, 0, 0, 0 \\  
     0, 1, 1, -1, 0, -1, -1 \\ 
     0, 1, 1, -1, 0, 0, 0 \\  
     1, 0, 0, 0, 0, 1, 0 \\   
     0, 0, -1, 1, 1, 0, 0 \\  
     -1, 0, -1, 1, 0, 0, 0 \\ 
      0, -1, 0, 0, -1, 0, 1 \\ 
      0, -1, 1, 0, -1, -1, 0 \\ 
     -1, 0, 1, 0, 0, -1, 0 \\  
    \end{array}$ \\
  \hline
\end{tabular}

  \begin{tabular} {|c|c|c|c|c|c|c|}
   \hline
       & $\begin{array}{c} \text{Orbit $2$} \\  ~ \end{array}$     
       & $\begin{array}{c} \text{Orbit $3$} \\  ~ \end{array}$     
       & $\begin{array}{c} \text{Orbit $4$} \\  ~ \end{array}$  \cr   
  \hline
   $\begin{array}{c}
     \gamma^{\vee} \\
     \InvDualSC\gamma^{\vee} \\
     (\InvDualSC)^2\gamma^{\vee} \\
     (\InvDualSC)^3\gamma^{\vee}\\
     (\InvDualSC)^4\gamma^{\vee}\\
     (\InvDualSC)^5\gamma^{\vee}\\
     (\InvDualSC)^6\gamma^{\vee}\\
     (\InvDualSC)^7\gamma^{\vee}\\
     (\InvDualSC)^8\gamma^{\vee}\\
     (\InvDualSC)^9\gamma^{\vee}\\
    \end{array}$
   & 
   $\begin{array}{c}
      \fbox{0, 0, 0, 0, 0, 0, -1} \\ 
      0, 1, 0, -1, 1, 0, 0 \\  
      1, 1, 0, -1, 0, 1, -1 \\ 
      1, 1, 0, -1, 1, 0, 0 \\  
      0, 1, 0, 0, 0, 0, -1 \\  
     \fbox{0, 0, 0, 0, 0, 0, 1} \\ 
      0, -1, 0, 1, -1, 0, 0 \\  
      -1, -1, 0, 1, 0, -1, 1 \\ 
      -1, -1, 0, 1, -1, 0, 1 \\ 
       0, -1, 0, 0, 0, 0, 1 \\  
    \end{array}$
   & 
    $\begin{array}{c}
    \fbox{0, 0, 0, 0, 1, -1, 0} \\  
    -1, 1, 0, 0, 0, 0, -1  \\ 
    1, 0, 0, -1, 0, 1, 1   \\ 
    1, 0, 0, 0, 0, 0, -1   \\ 
    -1, 1, 0, 0, 1, -1, 0  \\ 
    \fbox{0, 0, 0, 0, -1, 1, 0} \\  
    1, -1, 0, 0, 0, 0, 1   \\ 
    -1, 0, 0, 1, 0, -1, -1 \\ 
    -1, 0, 0, 0, 0, 0, 1   \\  
     1, -1, 0, 0, -1, 1, 0 \\ 
    \end{array}$
   & 
     $\begin{array}{c}
     \fbox{0, 0, 0, 0, -1, 0, -1} \\ 
      0, 0, 1, -1, 0, -1, 1 \\ 
      0, 0, 1, 0, -1, 0, -1 \\ 
     \fbox{0, 0, 0, 0, 1, 0, 1} \\  
      0, 0, -1, 1, 0, 1, -1 \\ 
      0, 0, -1, 0, 1, 0, 1  \\ 
      ~ \\
      ~ \\
      ~ \\
      ~ \\
    \end{array}$ \\
  \hline
\end{tabular}
  \vspace{2mm}
  \caption{\hspace{3mm} \Small ${\bf E_7(a_3)}$, there exist $4$ semi-Coxeter orbits, one of length $30$,
  two of length $10$ and one of length $6$.  All orbits are self-opposite}
  \label{E7a3_linkages}
  \end{table}



\begin{table}[H]
  \centering
  \scriptsize
  \renewcommand{\arraystretch}{1.5}  
  \begin{tabular} {|c|c|c|c|c|c|c|}
  \hline
    $\begin{array}{c} ~\\ \includegraphics[scale=0.5]{E7a4_upd_28feb2011.eps} \\  {\bf E_7(a_4)} \end{array}$
       & $\begin{array}{c} \text{Orbit $1$} \\ ~ \end{array}$     
       & $\begin{array}{c} \text{Orbit $2$} \\ ~ \end{array}$     
       & $\begin{array}{c} \text{Orbit $3$} \\ ~ \end{array}$     
       & $\begin{array}{c} \text{Orbit $4$} \\ \text{(no unicolored)} \end{array}$ \cr
  \hline
   $\begin{array}{c}
     \gamma^{\vee} \\
     \InvDualSC\gamma^{\vee} \\
     (\InvDualSC)^2\gamma^{\vee} \\
     (\InvDualSC)^3\gamma^{\vee}\\
     (\InvDualSC)^4\gamma^{\vee}\\
     (\InvDualSC)^5\gamma^{\vee}\\
    \end{array}$
    & 
   $\begin{array}{c}
      \fbox{0, 0, 0, 0, 0, 0, 1}  \\ 
      1, -1,  0, 0, -1, -1, 1 \\  
      1, -1,  0, 0,  0,  0, 1 \\  
      \fbox{0, 0, 0, 0, 0, 0, -1}  \\ 
      -1, 1,  0, 0,  1, 1, -1 \\  
      -1, 1,  0, 0,  0, 0, -1 \\  
   \end{array}$
    & 
   $\begin{array}{c}
   \fbox{0, 0, 0, 0, 0, -1, 0} \\  
       1, 0, -1, 0, 1, -1, 1 \\ 
       1, 0, -1, 0, 0, -1, 0 \\ 
   \fbox{0, 0, 0, 0, 0, 1, 0} \\  
       -1, 0, 1, 0, -1, 1, -1 \\ 
       -1, 0, 1, 0, 0, 1, 0 \\ 
    \end{array}$
   &  
     $\begin{array}{c}
     \fbox{0, 0, 0, 0, -1, 0, 0}  \\ 
      0, -1,  1, 0, -1, 1,  1 \\ 
      0, -1,  1, 0, -1, 0,  0 \\ 
     \fbox{0, 0, 0, 0, 1, 0, 0}  \\ 
      0,  1, -1, 0,  1, -1, -1 \\ 
      0,  1, -1, 0,  1,  0,  0 \\ 
   \end{array}$
    & 
   $\begin{array}{c}
       -1, 0, 0, 0, 0, 1, 0 \\  
       0, 0, 1, -1, -1, 0, 0 \\ 
       1, 0, 1, -1, 0, 0, 1 \\  
       1, 0, 0, 0, 0, -1, 0 \\  
       0, 0, -1, 1, 1, 0, 0 \\  
       -1, 0, -1, 1, 0, 0, -1 \\ 
    \end{array}$ \\
   \hline
\end{tabular}


  \begin{tabular} {|c|c|c|c|c|c|c|}
   \hline  
       \hspace{2.35cm}
       & $\begin{array}{c} \text{Orbit $5$} \\ \text{(no unicolored)} \end{array}$
       & $\begin{array}{c} \text{Orbit $6$} \\ \text{(no unicolored)} \end{array}$
       & $\begin{array}{c} \text{Orbit $7$} \\ \text{(no unicolored)} \end{array}$
       & $\begin{array}{c} \text{Orbit $8$} \\ \text{(no unicolored)} \end{array}$ \cr
  \hline
   $\begin{array}{c}
     \gamma^{\vee} \\
     \InvDualSC\gamma^{\vee} \\
     (\InvDualSC)^2\gamma^{\vee} \\
     (\InvDualSC)^3\gamma^{\vee}\\
     (\InvDualSC)^4\gamma^{\vee}\\
     (\InvDualSC)^5\gamma^{\vee}\\
    \end{array}$
    &  
    $\begin{array}{c}
       \hspace{1.25mm} 0, 0, -1, 0, 0, -1, 0  \hspace{1.25mm} \\  
       1, 0, 0, -1, 0, 0, 1 \\  
       1, 0, 1, -1, -1, 0, 0 \\  
       0, 0, 1, 0, 0, 1, 0  \\  
       -1, 0, 0, 1, 0, 0, -1 \\  
       -1, 0, -1, 1, 1, 0, 0 \\ 
   \end{array}$
    &  
   $\begin{array}{c}
      \hspace{2.25mm}  0, 1, 1, -1, 0, 1, 0 \hspace{2.25mm}  \\ 
       0, 0, 1, 0, -1, 0, 0 \\  
       0, -1, 0, 1, 0, 0, 1 \\  
       0, -1, -1, 1, 0, -1, 0 \\ 
       0, 0, -1, 0, 1, 0, 0 \\ 
       0, 1, 0, -1, 0, 0, -1 \\  
   \end{array}$
    &  
    $\begin{array}{c}
      \hspace{1.4mm} 0, -1, 0, 0, -1, 0, 0  \hspace{1.4mm} \\ 
       0, 0, 1, -1, 0, 1, 0   \\ 
       0, 1, 1, -1, 0, 0, -1   \\ 
       0, 1, 0, 0, 1, 0, 0 \\ 
       0, 0, -1, 1, 0, -1, 0 \\ 
       0, -1, -1, 1, 0, 0, 1 \\ 
   \end{array}$
    &  
    $\begin{array}{c}
      -1, -1, 0, 1, -1, 0, 0  \\ 
       0, -1, 0, 0, 0, 0, 1   \\ 
       1, 0, 0, -1, 0, -1, 0   \\ 
       1, 1, 0, -1, 1, 0, 0 \\ 
       0, 1, 0, 0, 0, 0, -1 \\ 
      -1, 0, 0, 1, 0, 1, 0 \\ 
   \end{array}$ \\
   \hline
\end{tabular}


  \begin{tabular} {|c|c|c|c|c|c|c|}
   \hline  
     \hspace{2.35cm}
       & $\begin{array}{c} \text{Orbit $9$}  \\ \text{(no unicolored)} \end{array}$
       & $\begin{array}{c} \text{Orbit $10$} \\ ~  \end{array}$ \cr
  \hline
   $\begin{array}{c}
     \gamma^{\vee} \\
     \InvDualSC\gamma^{\vee} \\
     (\InvDualSC)^2\gamma^{\vee} \\
     (\InvDualSC)^3\gamma^{\vee}\\
     (\InvDualSC)^4\gamma^{\vee}\\
     (\InvDualSC)^5\gamma^{\vee}\\
    \end{array}$
    & 
   $\begin{array}{c}
      \hspace{2.3mm} 0, 1, 0, -1, 1, 0, 0 \hspace{2.3mm} \\  
       1, 1, 0, -1, 0, -1, 0 \\ 
       1, 0, 0, 0, 0, 0, 1   \\ 
       0, -1, 0, 1, -1, 0, 0 \\ 
       -1, -1, 0, 1, 0, 1, 0 \\ 
       -1, 0, 0, 0, 0, 0, -1 \\ 
    \end{array}$
   &
   $\begin{array}{c}
       \hspace{4.6mm}
       \fbox{0, 0, 0, 0, 1, 1, 1} \hspace{4.6mm} \\ 
       \fbox{0, 0, 0, 0, -1, -1, -1}  \\  
       ~ \\
       ~ \\
       ~ \\
       ~ \\
    \end{array}$ \\
    \hline
\end{tabular}
  \vspace{2mm}
  \caption{\hspace{3mm}  ${\bf E_7(a_4)}$, $10$ semi-Coxeter orbits, nine of length $6$,
  one of length $2$}
  \label{E7a4_linkages}
  \end{table}



\newpage
~\\
\begin{table}[t]
  \centering
  \scriptsize
  \renewcommand{\arraystretch}{1.5}  
  \begin{tabular} {|c|c|c|c|c|c|c|}
  \hline   
    $\begin{array}{c} ~\\ \includegraphics[scale=0.3]{E7_pure.eps} \\ {\bf E_7} \end{array}$
       & $\begin{array}{c} \text{Orbit $1$} \\ ~ \\ ~ \end{array}$ 
       & $\begin{array}{c} \text{Orbit $2$} \\ ~ \\ ~ \end{array}$ 
       & $\begin{array}{c} \text{Orbit $3$} \\ ~ \\ ~ \end{array}$  
       & $\begin{array}{c} \text{Orbit $4$} \\ ~ \\ ~ \end{array}$ \cr 
  \hline   
   $\begin{array}{c}
     \gamma^{\vee} \\
     \InvDualSC\gamma^{\vee} \\
     (\InvDualSC)^2\gamma^{\vee} \\
     (\InvDualSC)^3\gamma^{\vee}\\
     (\InvDualSC)^4\gamma^{\vee}\\
     (\InvDualSC)^5\gamma^{\vee}\\
     (\InvDualSC)^6\gamma^{\vee} \\
     (\InvDualSC)^7\gamma^{\vee}\\
     (\InvDualSC)^8\gamma^{\vee}\\
     (\InvDualSC)^9\gamma^{\vee}\\
     (\InvDualSC)^{10}\gamma^{\vee} \\
     (\InvDualSC)^{11}\gamma^{\vee}\\
     (\InvDualSC)^{12}\gamma^{\vee}\\
     (\InvDualSC)^{13}\gamma^{\vee}\\
     (\InvDualSC)^{14}\gamma^{\vee} \\
     (\InvDualSC)^{15}\gamma^{\vee}\\
     (\InvDualSC)^{16}\gamma^{\vee}\\
     (\InvDualSC)^{17}\gamma^{\vee}\\
    \end{array}$
    & 
   $\begin{array}{c}
  \fbox{0, 0, 0, -1, 0, 0, 0}  \\ 
       0, 0, 0, 1, 0, 0, -1 \\  
       0, 1,  0, 0,  -1,  0, 0 \\  
       1, 0, 1, 0, -1, -1, 0  \\ 
       1, 1,  0, 0,  -1, 0, -1 \\  
       0, 1,  1, 1,  -1, 0, -1 \\  
       1, 1,  0, 0,  -1, -1, 0 \\  
       1, 0,  1, 0,  -1, 0, 0 \\  
       0, 1,  0, 0,  0, 0, -1 \\  
   \fbox{0, 0,  0, 1,  0, 0,  0} \\  
       0, 0,  0, -1,  0, 0,  1 \\  
       0, -1,  0, 0,  1, 0,  0 \\  
      -1, 0,  -1, 0,  1, 1,  0 \\  
      -1, -1,  0, 0,  1, 0,  1 \\  
       0, -1, -1, -1, 1, 0,  1 \\  
      -1, -1,  0, 0,  1, 1, 0 \\  
      -1, 0,  -1, 0,  1, 0, 0 \\  
       0, -1,  0, 0,  0, 0, 1 \\  
   \end{array}$
    & 
   $\begin{array}{c}
  \fbox{-1, 0, 1, 0, 0, 0, 0} \\ 
       1, 0,  -1, 0, 0, -1, 0 \\  
       0, 0,  1, 0,  -1,  1, 0 \\ 
       0, 1, 0, 0, 0, -1, -1 \\   
       1, 0,  0, 1, -1, 0, 0 \\   
       0, 1,  1, -1,  -1, 0, 0 \\ 
       1, 0,  0, 1,  0, -1, -1 \\  
       0, 1,  0, 0,  -1, 1,  0 \\  
       0, 0,  1, 0,  0, -1, 0 \\  
  \fbox{1, 0,  -1, 0,  0, 0,  0} \\  
       -1, 0,  1, 0,  0, 1,  0 \\  
       0, 0,  -1, 0,  1, -1, 0 \\  
       0, -1,  0, 0,  0, 1,  1 \\  
       -1, 0,  0, -1,  1, 0, 0 \\  
       0, -1, -1, 1,  1, 0, 0 \\  
      -1, 0,  0, -1,  0, 1, 1 \\  
       0, -1,  0, 0,  1, -1, 0 \\  
       0, 0,  -1, 0,  0, 1, 0 \\  
    \end{array}$
   &  
     $\begin{array}{c}
    \fbox{1, -1, 0, 0, 0, 0, 0} \\  
      -1, 1, 0, 0, 0, 1, -1 \\ 
      0, 0, 0, 1, 0, -1, 0  \\ 
      1, 0,  0, -1, -1, 0,  1 \\ 
      0, 0,  1, 0, 0, 0, -1 \\ 
      0, 1,  -1, 1, 0, 0, -1 \\ 
      0, 0,  1, 0, -1, 0, 1 \\ 
      1, 0,  0, -1,  0, -1, 0 \\ 
      0, 0, 0, 1, 0, 1, -1 \\  
    \fbox{-1, 1, 0, 0, 0, 0, 0} \\  
      1, -1, 0, 0, 0, -1, 1 \\ 
      0, 0, 0, -1, 0,  1, 0 \\  
      -1,  0, 0, 1,  1, 0, -1 \\ 
      0,  0, -1, 0,  0,  0,  1 \\ 
      0, -1, 1, -1,  0, 0, 1  \\ 
      0,  0, -1, 0,  1,  0, -1 \\ 
      -1, 0, 0, 1, 0, 1, 0 \\ 
      0,  0, 0, -1,  0,  -1, 1 \\ 
    \end{array}$
    &  
   $\begin{array}{c}
      \fbox{0, 1, -1, -1, 0, 0, 0}  \\ 
       \fbox{0, -1, 1, 1, 0, 0, 0} \\  
       ~ \\
       ~ \\
       ~ \\
       ~ \\
       ~ \\
       ~ \\
       ~ \\
       ~ \\
       ~ \\
       ~ \\
       ~ \\
       ~ \\
       ~ \\
       ~ \\
       ~ \\
       ~ \\
   \end{array}$ \\
   \hline
\end{tabular}
  \vspace{2mm}
  \caption{\hspace{3mm}  ${\bf E_7}$, $4$ semi-Coxeter orbits, three of length $18$
  and one of length $2$.  All orbits are self-opposite.
  All orbits contain $\beta$-unicolored linkage diagrams}
  \label{E7pure_linkages}
   \vspace{20mm}
  \end{table}


\begin{table}[H]
  \centering
  \scriptsize
  \renewcommand{\arraystretch}{1.5}  
  \begin{tabular} {|c|c|c|c|c|c|c|}
  \hline   
    $\begin{array}{c} ~\\ \includegraphics[scale=0.3]{D7a1.eps} \\ {\bf D_7(a_1)} \end{array}$
       & $\begin{array}{c} \text{Orbit $1$} \\ ~ \end{array}$ 
       & $\begin{array}{c} \text{Orbit $2$} \\ ~ \end{array}$ 
       & $\begin{array}{c} \text{Orbit $3$} \\ ~ \end{array}$  
       & $\begin{array}{c} \text{Orbit $4$} \\ ~ \end{array}$ \cr 
  \hline   
   $\begin{array}{c}
     \gamma^{\vee} \\
     \InvDualSC\gamma^{\vee} \\
     (\InvDualSC)^2\gamma^{\vee} \\
     (\InvDualSC)^3\gamma^{\vee}\\
     (\InvDualSC)^4\gamma^{\vee}\\
     (\InvDualSC)^5\gamma^{\vee}\\
     (\InvDualSC)^6\gamma^{\vee} \\
     (\InvDualSC)^7\gamma^{\vee}\\
     (\InvDualSC)^8\gamma^{\vee}\\
     (\InvDualSC)^9\gamma^{\vee}\\
     (\InvDualSC)^{10}\gamma^{\vee} \\
     (\InvDualSC)^{11}\gamma^{\vee}\\
     (\InvDualSC)^{12}\gamma^{\vee}\\
     (\InvDualSC)^{13}\gamma^{\vee}\\
     (\InvDualSC)^{14}\gamma^{\vee} \\
     (\InvDualSC)^{15}\gamma^{\vee}\\
     (\InvDualSC)^{16}\gamma^{\vee}\\
     (\InvDualSC)^{17}\gamma^{\vee}\\
     (\InvDualSC)^{18}\gamma^{\vee}\\
     (\InvDualSC)^{19}\gamma^{\vee}\\
    \end{array}$
    & 
   $\begin{array}{c}
  \fbox{0, 0, 1, 0, 0, 0, 0}  \\ 
       0, 0, -1, 0, 1, 1, 0 \\  
      -1, 0, -1, 0, 1, 0, 1 \\  
      -1, -1, 0, -1, 1, -1, 1  \\ 
      -1, -1,  0, 0,  1, 0, 0 \\  
       0, 0,  -1, 0,  0, 1, 0 \\  
       0, 1,  0, 0,  -1, 0, 0 \\  
       1, 0,  1, 0,  -1, -1, -1 \\  
       1, 0,  1, 1,  -1, 0, -1 \\  
       1, 1,  0, 0,  -1, 1,  0 \\  
   \fbox{0, 1,  0, 0,  0, 0,  0} \\  
       0, -1,  0, 0,  1, -1,  0 \\  
      -1, -1,  0, 0,  1, 0,  1 \\  
      -1,  0, -1, -1,  1, 1,  1 \\  
      -1, 0, -1, 0, 1, 0,  0 \\  
       0, -1,  0, 0,  0, -1, 0 \\  
       0,  0,  1, 0,  -1, 0, 0 \\  
       1,  1,  0, 0,  -1, 1, -1 \\  
       1,  1,  0, 1,  -1, 0, -1 \\  
       1, 0,  1, 0,  -1, -1, 0 \\  
   \end{array}$
    & 
   $\begin{array}{c}
  \fbox{1, -1, 0, 0, 0, 0, 0} \\ 
       -1, 1,  0, 0, 0, 1, 1 \\  
       0, 0,  -1, -1,  1,  0, 0 \\ 
       -1, -1, 0, 1, 1, -1, 0 \\   
       0, -1,  0, -1, 0, 0, 1 \\   
       -1, 1,  0, 0,  0, 1, 0 \\ 
       1, 0,  -1, 0,  0, 0, -1 \\  
       0, 0,  1, 1,  -1, -1,  0 \\  
       1, 0,  1, -1,  -1, 0, 0 \\  
       0, 1,  0, 1,  0, 1,  -1 \\  
  \fbox{1, 0,  -1, 0,  0, 0,  0} \\  
       -1, 0,  1, 0,  0, -1, 1 \\  
       0, -1,  0, -1,  1, 0,  0 \\  
       -1, 0,  -1, 1,  1, 1, 0 \\  
       0, 0, -1, -1,  0, 0, 1 \\  
      -1, 0,  1, 0,  0, -1, 0 \\  
       1, -1,  0, 0,  0, 0, -1 \\  
       0, 1,  0, 1,  -1, 1, 0 \\  
       1, 1,  0, -1,  -1, 0, 0 \\  
       0, 0,  1, 1,  0, -1, -1 \\  
    \end{array}$
   &  
     $\begin{array}{c}
    \fbox{0, -1, 0, 1, 0, 0, 0} \\  
      0, 1, 0, -1, -1, 1, 1 \\ 
      0, 1, 0, 0, 0, 0, -1  \\ 
      1, -1,  0, 1, 0, -1,  -1 \\ 
      0, 0,  1, 0, -1, 0, 1 \\ 
      0, 1,  0, -1, 0, 1, 0 \\ 
      0, 0,  -1, 1, 1, 0, -1 \\ 
      0, -1,  0, 0,  0, -1, 1 \\ 
      -1, 0, 1, -1, 0, 0, 1 \\  
      0, 0, -1, 0, 1, 1, -1 \\  
    \fbox{0, 0, -1, 1, 0, 0, 0} \\ 
      0, 0, 1, -1, -1,  -1, 1 \\  
      0,  0, 1, 0,  0, 0, -1 \\ 
      1,  0, -1, 1,  0,  1, -1 \\ 
      0, 1, 0, 0,  -1, 0, 1  \\ 
      0,  0, 1, -1,  0,  -1, 0 \\ 
      0, -1, 0, 1, 1, 0, -1 \\ 
      0,  0, -1, 0,  0,  1, 1 \\ 
      -1,  1, 0, -1,  0,  0, 1 \\ 
      0,  -1, 0, 0,  1,  -1, -1 \\ 
    \end{array}$
    &  
   $\begin{array}{c}
      \fbox{-1, 0, 1, 1, 0, 0, 0}  \\ 
       1, 0, -1, -1, 0, 1, 0 \\  
      \fbox{-1, 1, 0, 1, 0 , 0, 0} \\
       1, -1, 0, -1, 0, -1, 0 \\
       ~ \\
       ~ \\
       ~ \\
       ~ \\
       ~ \\
       ~ \\
       ~ \\
       ~ \\
       ~ \\
       ~ \\
       ~ \\
       ~ \\
       ~ \\
       ~ \\
       ~ \\
       ~ \\
   \end{array}$ \\
   \hline
\end{tabular}
  \begin{tabular} {|c|c|c|c|c|c|c|}
  \hline   
       & $\begin{array}{c} \text{Orbit $5$} \end{array}$ 
       & $\begin{array}{c} \text{Orbit $6$}  \end{array}$ \cr 
  \hline   
   $\begin{array}{c}
     \gamma^{\vee} \\
     \InvDualSC\gamma^{\vee} \\
     (\InvDualSC)^2\gamma^{\vee} \\
     (\InvDualSC)^3\gamma^{\vee}\\
     (\InvDualSC)^4\gamma^{\vee}\\
     (\InvDualSC)^5\gamma^{\vee}\\
     (\InvDualSC)^6\gamma^{\vee} \\
     (\InvDualSC)^7\gamma^{\vee}\\
     (\InvDualSC)^8\gamma^{\vee}\\
     (\InvDualSC)^9\gamma^{\vee}\\
    \end{array}$
    & 
   $\begin{array}{c}
  \fbox{0, 0, 0, -1, 0, 0, 0}  \\ 
       0, 0, 0, 1, 0, 0, -1 \\  
       1, 0, 0, 0, -1, 0, 0 \\  
       0, 1, 1, 0, -1, 0, 0  \\ 
       1, 0,  0, 0,  0, 0, -1 \\  
  \fbox{0, 0,  0, 1,  0, 0, 0} \\  
       0, 0,  0, -1,  0, 0, 1 \\  
       -1, 0,  0, 0,  1, 0, 0 \\  
       0, -1,  -1, 0,  1, 0, 0 \\  
       -1, 0,  0, 0,  0, 0,  1 \\  
   \end{array}$
    &  
   $\begin{array}{c}
      \fbox{0, 0, 0, 0, 0, -1, 0}  \\ 
       0, -1, 1, 0, 0, -1, 0 \\  
       \fbox{0, 0, 0, 0, 0, 1, 0} \\
       0, 1, -1, 0, 0, 1, 0 \\
       ~ \\
       ~ \\
       ~ \\
       ~ \\
       ~ \\
       ~ \\
   \end{array}$ \\
   \hline
\end{tabular}
  \vspace{2mm}
  \caption{\hspace{3mm}  ${\bf D_7(a_1)}$, there exist $10$ semi-Coxeter orbits.
  Semi-Coxeter orbits $1$-$4$ (three of length $20$
  and one of length $4$) belong to the first $E$-type component. Every orbit out of $1$-$4$ has
  the opposite one ($7$-$10$) lying in the second $E$-type component. Orbits $5$ and $6$ are self-opposite.
  All orbits contain $\alpha$-unicolored or  $\beta$-unicolored linkage diagrams}
  \label{D7a1_linkages}
   \vspace{20mm}
  \end{table}

\begin{table}[H]
  \centering
  \scriptsize
  \renewcommand{\arraystretch}{1.5}  
  \begin{tabular} {|c|c|c|c|c|c|c|}
  \hline   
    $\begin{array}{c} ~\\ \includegraphics[scale=0.3]{D7a1.eps} \\ {\bf D_7(a_2)} \end{array}$
       & $\begin{array}{c} \text{Orbit $1$} \\ ~ \end{array}$ 
       & $\begin{array}{c} \text{Orbit $2$} \\ ~ \end{array}$ 
       & $\begin{array}{c} \text{Orbit $3$} \\ ~ \end{array}$  
       & $\begin{array}{c} \text{Orbit $4$} \\ ~ \end{array}$ \cr 
  \hline   
   $\begin{array}{c}
     \gamma^{\vee} \\
     \InvDualSC\gamma^{\vee} \\
     (\InvDualSC)^2\gamma^{\vee} \\
     (\InvDualSC)^3\gamma^{\vee}\\
     (\InvDualSC)^4\gamma^{\vee}\\
     (\InvDualSC)^5\gamma^{\vee}\\
     (\InvDualSC)^6\gamma^{\vee} \\
     (\InvDualSC)^7\gamma^{\vee}\\
     (\InvDualSC)^8\gamma^{\vee}\\
     (\InvDualSC)^9\gamma^{\vee}\\
     (\InvDualSC)^{10}\gamma^{\vee} \\
     (\InvDualSC)^{11}\gamma^{\vee}\\
     (\InvDualSC)^{13}\gamma^{\vee} \\
     (\InvDualSC)^{14}\gamma^{\vee} \\
     (\InvDualSC)^{15}\gamma^{\vee}\\
     (\InvDualSC)^{16}\gamma^{\vee}\\
     (\InvDualSC)^{17}\gamma^{\vee}\\
     (\InvDualSC)^{18}\gamma^{\vee} \\
     (\InvDualSC)^{19}\gamma^{\vee}\\
     (\InvDualSC)^{20}\gamma^{\vee}\\
     (\InvDualSC)^{21}\gamma^{\vee}\\
     (\InvDualSC)^{22}\gamma^{\vee} \\
     (\InvDualSC)^{23}\gamma^{\vee}\\
   \end{array}$
    & 
   $\begin{array}{c}
  \fbox{0, 0, 1, 0, 0, 0, 0}  \\ 
       0, 0, -1, 0, 1, 1, 0 \\  
      -1, 0, -1, -1, 1, 1, 1 \\  
      -1, 0, -1, 0, 1, 0, 0  \\ 
       0, -1, 0,  0, 0,  -1, 0 \\  
       0, 0,  1, 1,  -1, -1, 0 \\  
       1, 0,  1, 0,  -1, 0, -1 \\  
       1, 1,  0, 0,  -1, 1, 0 \\  
       0, 1,  0, -1,  0, 1, 0 \\  
       0, 0,  -1, 0,  1, 0,  0 \\  
       -1, -1,  0, 0,  1, -1,  1 \\  
       -1, -1,  0, 1,  1, -1,  0 \\  
   \fbox{0, -1,   0, 0,  0,  0, 0}  \\ 
       0, 1, 0, 0, -1, 1, 0 \\  
       1, 1, 0, -1, -1, 1, -1 \\ 
       1, 1, 0, 0, -1, 0, 0  \\ 
       0, 0, 1,  0, 0,  -1, 0 \\  
       0, -1,  0, 1,  1, -1, 0 \\  
       -1, -1,  0, 0,  1, 0, 1 \\  
       -1, 0,  -1, 0,  1, 1, 0 \\  
       0, 0,  -1, -1,  0, 1, 0 \\  
       0, 1,  0, 0,  -1, 0,  0 \\  
       1, 0,  1, 0,  -1, -1,  -1 \\  
       1, 0,  1, 1,  -1, -1,  0 \\  
   \end{array}$
    & 
   $\begin{array}{c}
  \fbox{0, 1, 0, 1, 0, 0, 0} \\ 
       0, -1,  0, -1, 1, 0, 0 \\  
       -1, 0,  -1, 1,  1,  0, 1 \\ 
       -1, -1, 0, -1, 1, 0, 0 \\  
   \fbox{0, 0, -1, 1, 0, 0, 0} \\   
       0, 0,  1, -1,  -1, 0, 0 \\ 
       1, 1,  0, 1,  -1, 0, -1 \\  
       1, 0,  1, -1,  -1, 0,  0 \\  
       ~ \\
       ~ \\
       ~ \\
       ~ \\
       ~ \\
       ~ \\
       ~ \\
       ~ \\
       ~ \\
       ~ \\
       ~ \\
       ~ \\
       ~ \\
       ~ \\
       ~ \\
       ~ \\
    \end{array}$
   &  
     $\begin{array}{c}
    \fbox{1, -1, 0, 0, 0, 0, 0} \\  
      -1, 1, 0, 0, 0, 1, 1 \\ 
      0, 0, -1, -1, 1, 1, -1  \\ 
      0, 0,  -1, 0, 0, 0,  1 \\ 
      -1, 0,  1, 0, 0, -1, 0 \\ 
      1, -1,  0, 1, 0, -1, -1 \\ 
      0, 0,   1, 0, -1, 0, 1 \\ 
      0,  1,  0, 0,  0, 1, -1 \\ 
      1, 0, -1, -1, 0, 1, 0 \\  
      -1, 1, 0, 0, 0, 0, 1 \\  
      0, -1, 0, 0, 1, -1, -1 \\ 
      0, -1, 0,  1, 0,  -1, 1 \\  
    \fbox{-1, 0, 1, 0, 0, 0, 0} \\  
      1, 0, -1, 0, 0, 1, -1 \\ 
      0, 1, 0, -1, -1, 1, 1  \\ 
      0, 1,  0, 0, 0, 0,  -1 \\ 
      1, -1,  0, 0, 0, -1, 0 \\ 
     -1, 0,  1, 1, 0, -1,  1 \\ 
      0, -1,  0, 0, 1, 0, -1 \\ 
      0,  0,  -1, 0,  0, 1, 1 \\ 
      -1, 1, 0, -1, 0, 1, 0 \\  
      1, 0, -1, 0, 0, 0, -1 \\  
      0, 0, 1, 0, -1, -1, 1 \\ 
      0, 0, 1,  1, 0,  -1, -1 \\  
    \end{array}$
    &  
     $\begin{array}{c}
    \fbox{1, 0, -1, 1, 0, 0, 0} \\  
      -1, 0, 1, -1, 0, 0, 1 \\ 
      0, 0, -1, 1, 1, 0, -1  \\ 
      0, -1,  0, -1, 0, 0,  1 \\ 
    \fbox{-1, 1,  0, 1, 0, 0, 0} \\ 
      1, -1,  0, -1, 0, 0, -1 \\ 
      0, 1,  0, 1, -1, 0, 1 \\ 
      0, 0,  1, -1,  0, 0, -1 \\ 
       ~ \\
       ~ \\
       ~ \\
       ~ \\
       ~ \\
       ~ \\
       ~ \\
       ~ \\
       ~ \\
       ~ \\
       ~ \\
       ~ \\
       ~ \\
       ~ \\
       ~ \\
       ~ \\
   \end{array}$ \\
   \hline
\end{tabular}
 \begin{tabular} {|c|c|c|c|c|c|c|}
  \hline   
       & $\begin{array}{c} \text{Orbit $5$} \end{array}$ 
       & $\begin{array}{c} \text{Orbit $6$}  \end{array}$ \cr 
  \hline   
   $\begin{array}{c}
     \gamma^{\vee} \\
     \InvDualSC\gamma^{\vee} \\
     (\InvDualSC)^2\gamma^{\vee} \\
     (\InvDualSC)^3\gamma^{\vee}\\
     (\InvDualSC)^4\gamma^{\vee}\\
     (\InvDualSC)^5\gamma^{\vee}\\
     (\InvDualSC)^6\gamma^{\vee} \\
     (\InvDualSC)^7\gamma^{\vee}\\
    \end{array}$
    & 
   $\begin{array}{c}
  \fbox{0, 0, 0, -1, 0, 0, 0}  \\ 
       0, 0, 0, 1, 0, -1, 0 \\  
       0, -1, 1, 0, 0, -1, 0 \\  
  \fbox{0, 0, 0, 1, 0, 0, 0}  \\ 
       0, 0,  0, -1,  0, 1, 0 \\  
       0, 1,  -1, 0,  0, 1, 0 \\  
       ~ \\
       ~ \\
   \end{array}$
    &  
   $\begin{array}{c}
      \fbox{0, 0, 0, 0, 0, 0, 1} \\ 
       -1, 0, 0, 0, 1, 0, 0 \\   
       0, -1, -1, 0, 1, 0, 0 \\  
       -1, 0, 0, 0, 0, 0, 1 \\   
       \fbox{0, 0, 0, 0, 0, 0, -1} \\ 
       1, 0, 0, 0, -1, 0, 0 \\   
       0, 1, 1, 0, -1, 0, 0 \\  
       1, 0, 0, 0, 0, 0, -1 \\   
   \end{array}$ \\
   \hline
\end{tabular}
  \vspace{2mm}
  \caption{\hspace{3mm}  ${\bf D_7(a_2)}$, there exist $10$ semi-Coxeter orbits.
  Semi-Coxeter orbits $1$-$4$ (two of length $24$
  and two of length $8$) belong to the first $E$-type component. Every orbit out of $1$-$4$ has
  the opposite one ($7$-$10$) lying in the second $E$-type component.
  Orbits $5$ and $6$ are self-opposite.
  All orbits contain $\alpha$-unicolored or  $\beta$-unicolored linkage diagrams}
  \label{D7a2_linkages}
   \vspace{20mm}
  \end{table}

\begin{table}[H]
  \centering
  \scriptsize
  \renewcommand{\arraystretch}{1.8}  
  \begin{tabular} {|c|c|c|c|c|c|c|}
  \hline   
         $\begin{array}{c} ~\\ \includegraphics[scale=0.3]{D7_pure.eps} \\ {\bf D_7} \end{array}$
       & $\begin{array}{c} \text{Orbit $1$} \\ \text{(no unicolored)} \end{array}$ 
       & $\begin{array}{c} \text{Orbit $2$} \\ ~ \end{array}$ 
       & $\begin{array}{c} \text{Orbit $3$} \\ \text{(no unicolored)} \end{array}$ 
       & $\begin{array}{c} \text{Orbit $4$} \\ ~ \end{array}$ \cr 
  \hline   
   $\begin{array}{c}
     \gamma^{\vee} \\
     \InvDualSC\gamma^{\vee} \\
     (\InvDualSC)^2\gamma^{\vee} \\
     (\InvDualSC)^3\gamma^{\vee}\\
     (\InvDualSC)^4\gamma^{\vee}\\
     (\InvDualSC)^5\gamma^{\vee}\\
     (\InvDualSC)^6\gamma^{\vee} \\
     (\InvDualSC)^7\gamma^{\vee}\\
     (\InvDualSC)^8\gamma^{\vee}\\
     (\InvDualSC)^9\gamma^{\vee}\\
     (\InvDualSC)^{10}\gamma^{\vee} \\
     (\InvDualSC)^{11}\gamma^{\vee}\\
   \end{array}$
    & 
   $\begin{array}{c}
       0, 0, -1, 0, 0, 1, 0  \\ 
       0, -1, 1, -1, 0, 1, 1 \\  
       0, 0, -1, -1, 1, 0, 0 \\  
      -1, -1, 0, 1, 1, 0, -1  \\ 
       0, 0, -1,  0, 0,  0, 1 \\  
       0, 0,  1, -1,  -1, 1, 0 \\  
       1, 0,  0, 0,  0, -1, 0 \\  
      -1, 1,  0, 1,  0, -1, -1 \\  
       1, 0,  0, 1,  -1, 0, 0 \\  
       0, 1,  1, -1,  -1, 0,  1 \\  
       1, 0,  0, 0,  0, 0,  -1 \\  
       -1, 0,  0, 1,  1, -1,  0 \\  
   \end{array}$
    & 
   $\begin{array}{c}
       0, -1, -1, -1, 1, 1, 1 \\ 
      -1, -1,  0, -1, 1, 1, 0 \\  
       0, -1,  -1, 0,  1,  0, 0 \\ 
  \fbox{-1, 0, 0, 0, 0, 0, 0} \\  
       1, 0, 0, 0, -1, 0, 0 \\   
       0, 1,  1, 0, -1, -1, 0 \\ 
       1, 1,  0, 1, -1, -1, -1 \\  
       0, 1,  1,  1,  -1, -1,  0 \\  
       1, 1,  0,  0, -1, 0, 0 \\ 
  \fbox{0, 0, 1, 0, 0, 0, 0}  \\ 
       0, 0, -1, 0, 1, 0, 0  \\ 
      -1, -1,  0, 0, 1, 1, 0 \\ 
    \end{array}$
   &  
     $\begin{array}{c}
      -1, 0, 0, -1, 0, 1, 1 \\  
       1, -1, 0, -1, 0, 1, 0 \\ 
      -1, 0, 0, 0, 1, -1, 0  \\ 
       0, 0, -1, 1, 0, 0, -1 \\ 
       0, 0,  1, 0, -1, 0, 1 \\ 
       1, 1,  0, -1, -1, 0, 0 \\ 
       0, 0,  1, 1, 0, -1, -1 \\ 
       0, 1, -1, 1,  0, -1, 0 \\ 
       0, 0,  1, 0, -1, 1, 0 \\  
       1, 0, 0, -1, 0, 0, 1 \\  
      -1, 0, 0, 0, 1, 0, -1 \\ 
       0, -1, -1,  1, 1,  0, 0 \\  
    \end{array}$
    &  
     $\begin{array}{c}
    \fbox{0, -1, 1, 1, 0, 0, 0} \\ 
      0, 1, -1, -1, 0, 0, 1   \\ 
      0, -1, 1, 0, 0, 1, -1   \\ 
      0, 0, -1, 0, 1, -1, 1 \\   
     -1, 0,  0, 0, 0, 1, -1 \\ 
      1, -1,  0, -0, 0, 0, 1 \\ 
    \fbox{-1, 1,  0, -1, 0, 0, 0} \\ 
      1, -1,  0, 1,  0, 0, -1 \\ 
     -1,  1,  0, 0,  0, -1, 1 \\  
      1,  0,  0, 0,  -1, 1, -1 \\ 
      0,  0,  1, 0,  0, -1, 1 \\  
      0,  1,  -1, 0,  0, 0, -1  \\ 
   \end{array}$ \\
   \hline
\end{tabular}
 \begin{tabular} {|c|c|c|c|c|c|c|}
  \hline   
       & $\begin{array}{c} \text{Orbit $5$} \\ ~ \end{array}$ 
       & $\begin{array}{c} \text{Orbit $6$} \\ ~ \end{array}$ 
       & $\begin{array}{c} \text{Orbit $7$} \\ ~ \end{array}$ 
       & $\begin{array}{c} \text{Orbit $8$} \\ ~ \end{array}$ \cr 
  \hline   
   $\begin{array}{c}
     \gamma^{\vee} \\
     \InvDualSC\gamma^{\vee} \\
     (\InvDualSC)^2\gamma^{\vee} \\
     (\InvDualSC)^3\gamma^{\vee}\\
     (\InvDualSC)^4\gamma^{\vee}\\
     (\InvDualSC)^5\gamma^{\vee}\\
     (\InvDualSC)^6\gamma^{\vee} \\
     (\InvDualSC)^7\gamma^{\vee}\\
     (\InvDualSC)^8\gamma^{\vee}\\
     (\InvDualSC)^9\gamma^{\vee}\\
     (\InvDualSC)^{10}\gamma^{\vee} \\
     (\InvDualSC)^{11}\gamma^{\vee}\\
   \end{array}$
    & 
   $\begin{array}{c}
   \fbox{-1, 1, 0, 0, 0, 0, 0}  \\ 
       1, -1, 0, 0, 0, 1, 0 \\   
      -1, 0,  0, -1, 1, 0, 1 \\  
       0, -1, -1, 0, 1, 1, -1  \\ 
      -1, -1, 0,  0, 1, 0, 1 \\  
       0, 0, -1, -1,  0, 1, 0 \\  
   \fbox{0, -1, 1,  0, 0,  0, 0} \\  
       0, 1,  -1, 0,  0, -1, 0 \\  
       0, 0,  1, 1,  -1, 0, -1 \\  
       1, 1,  0, 0,  -1, -1,  1 \\  
       0, 1,  1, 0,  -1, 0,  -1 \\  
       1, 0,  0, 1,  0, -1,  0 \\  
   \end{array}$
    & 
   $\begin{array}{c}
   \fbox{-1,  0,  0,  1, 0, 0, 0} \\ 
       1,  0,  0, -1, -1, 1, 1 \\  
   \fbox{0, 0, 1, -1, 0, 0, 0} \\  
       0, 0, -1, 1,  1, -1, -1 \\ 
       ~ \\   
       ~ \\   
       ~ \\   
       ~ \\   
       ~ \\   
       ~ \\   
       ~ \\   
       ~ \\   
    \end{array}$
   &  
     $\begin{array}{c}
     \fbox{0, 0, 0,  0, 0, 0, -1} \\  
       0, 0, 0, 1, 0, -1, 0 \\ 
       0, 1, 0, 0, -1, 0, 0 \\ 
       1, 0, 1, 0, -1, 0, 0 \\ 
       0, 1, 0, 0, 0, -1, 0 \\ 
       0, 0,  0, 1, 0, 0, -1 \\ 
     \fbox{0, 0,  0, 0, 0, 0, 1} \\ 
       0, 0, 0, -1, 0, 1, 0 \\ 
       0, -1, 0, 0, 1, 0, 0 \\ 
      -1, 0, -1, 0, 1, 0, 0 \\ 
       0, -1, 0, 0, 0, 1, 0 \\ 
       0, 0,  0, -1, 0, 0,  1 \\ 
     \end{array}$
    &  
     $\begin{array}{c}
    \fbox{1, 0, -1, 0, 0, 0, 0} \\ 
    \fbox{-1, 0, 1, 0, 0, 0, 0} \\ 
      ~\\
      ~\\
      ~\\
      ~\\
      ~\\
      ~\\
      ~\\
      ~\\
      ~\\
      ~\\
   \end{array}$ \\
   \hline
\end{tabular}
  \vspace{2mm}
  \caption{\hspace{3mm}  ${\bf D_7}$, there exist $14$ semi-Coxeter orbits.
    Semi-Coxeter orbits $1$-$6$ ($5$ of length $12$ and one of length $4$ belong
    to the first $E$-type component. Every orbit out of $1$-$6$ has the opposite one
    ($9$-$14$) lying in the second $E$-type component. Orbits $7$ and $8$ are
    self-opposite. Orbits $1$, $3$ (and opposite to these orbits, i.e., $9$, $11$)
    do not contain unicolored linkage diagrams}
  \label{D7pure_linkages}
   \vspace{20mm}
  \end{table}